%% file: thesis.tex
\documentclass{cththesis}
\usepackage{rotating}
\usepackage{proof}
\usepackage[dvips]{xypic} 
\usepackage{ae}
\usepackage{mathrsfs}
\usepackage{amssymb,amsmath,amsthm}
\usepackage{braket}

\DeclareRobustCommand{\inferforall}[2]{\ldots \ #2 \ \ldots \ _{#1}}

\DeclareRobustCommand{\inference}[2]{\raisebox{-1.3ex}{\infer{#1}{#2}}}

\newtheorem{thm}{Theorem}[chapter]
\newtheorem{prop}[thm]{Proposition}
\newtheorem{lem}[thm]{Lemma}
\newtheorem{cor}[thm]{Corollary}
\newtheorem{conj}[thm]{Conjecture}
\theoremstyle{definition}
\newtheorem{defin}[thm]{Definition}
\newtheorem{que}[thm]{Question}

\input{head}

\graphicspath{{./}}

\newlength{\overmarginal}
\newlength{\undermarginal}
\setulmarginsandblock{12mm}{12mm}{*}
\checkandfixthelayout

\begin{document}

\renewcommand{\title}{Expansions, omitting types, and standard systems}
\renewcommand{\author}{FREDRIK ENGSTR\"OM}
\newcommand{\isbn}{91-7291-524-2}

\pagestyle{empty} 

%
%
%
%
%
%

\setulmarginsandblock{\overmarginal}{\undermarginal}{*}
\checkandfixthelayout

\begin{center}
\noindent THESIS FOR THE DEGREE OF DOCTOR OF PHILOSOPHY
\end{center}

\vspace{20mm}

\begin{center}
\bfseries \noindent
\Huge\title
\end{center}

\vspace{5mm}

\begin{center}
\author
\end{center}

\vspace{30mm}


\vfill\relax

\begin{center}
\noindent
Mathematical Sciences \\
CHALMERS UNIVERSITY OF TECHNOLOGY\\
AND G\"OTEBORG UNIVERSITY\\
G\"oteborg, Sweden 2004
\end{center}

\clearpage

\null
\vfill

\noindent \title \\
\author \\
ISBN \isbn

\vspace{15mm}

\noindent \copyright \author, 2004

\vspace{10mm}

\noindent Doktorsavhandlingar vid Chalmers Tekniska H\"ogskola \\
Ny serie nr 2206 \\
ISSN 0346-718x

\vspace{10mm}

\noindent Mathematical Sciences \\
Chalmers University of Technology and G\"oteborg University \\
SE--412 96 G\"oteborg \\
Sweden \\
Telephone $+$46--(0)--31--772 100

\vspace{10mm}


\vfill

\noindent Matematiskt Centrum \\
G\"oteborg, Sweden 2004

\clearpage

\noindent \textsc{\title\\
\noindent Fredrik Engstr\"om }

\vskip 2mm

\noindent Mathematical Sciences\\
Chalmers University of Technology and G\"oteborg University

\vskip 15mm

\textsc{\Large Abstract}

\vskip 5mm

\noindent Recursive saturation and resplendence are two important notions in models of arithmetic. Kaye, Kossak, and
Kotlarski introduced the notion of arithmetic saturation and argued that recursive saturation
might not be as rigid as first assumed. 

In this thesis we give further examples of variations of recursive saturation, all of which are connected with
expandability properties similar to resplendence. However, the expandability properties are stronger than 
resplendence and implies, in one
way or another, that the expansion not only satisfies a theory, but also omits a type. We conjecture that a special
version of this expandability is in fact equivalent to arithmetic saturation. We prove that another of these properties
is equivalent to $\beta$-saturation. We also introduce a variant on recursive saturation which makes sense in the context
of a standard predicate, and which is equivalent to a certain amount of ordinary saturation.

The theory of all models which omit a certain type $p(\bx)$ is also investigated. We define a proof system, which proves
a sentence if and only if it is true in all models omitting the type $p(\bx)$. The complexity of such proof systems are discussed and some
explicit examples of theories and types with high complexity, in a special sense, are given.

We end the thesis by a small comment on Scott's problem. The problem 
is to characterise standard systems of models of arithmetic.
We prove that, under the assumption of Martin's axiom,
        every Scott set of cardinality $<2^{\aleph_0}$ closed under
        arithmetic comprehension which has the
        countable chain condition is the standard system of some model of
        $\PA$. However, we do not know if there
        exists any such uncountable Scott sets.

\vfill


\vspace{1ex}


\cleardoublepage

\begin{center}
\textsc{\Large \noindent Acknowledgement}
\end{center}

\noindent My supervisor, Richard Kaye, deserves all my thanks, for taking me on and
for being such an inspiring teacher. Without his help this thesis would never have been. 

Jan Smith and Thierry Coquand made it possible for me to study models of arithmetic; thank you for
encouraging me.

I would also like to thank my friends at Chalmers; Julf, Leif, Bom, Kenny, Bettan, and Claesson for 
making my everyday life so pleasant.

At the University of Birmingham I would like to thank everyone for making
my stays there very enjoyable. My special thanks go to
John, David, Andrey, Stuart, Simon, John, Richard, Kamila, Panos, and Lylah.

At the philosophy department of G\"oteborg University, I thank Per Lindstr\"om, Christian Bennet, and
Dag Westerst\aa hl for general discussions on logic, interesting seminars, and for their interest in my
work.

Joakim Pirinen was kind to let me use two of 
his marvellous drawings in this thesis, thank you.

I would also like to thank Hilda for always encouraging and supporting me.

\vskip 10mm

\noindent\textbf{Thank you.}

\cleardoublepage






\pagestyle{ruled} 

\frontmatter
\maxtocdepth{subsection}
\tableofcontents

\mainmatter

\chapter*{Introduction}\addcontentsline{toc}{chapter}{Introduction}

The story of recursive saturation and resplendence, of which this thesis is a part, sprung from the 
useful notion of saturated models in classic model theory. However, as it turned out, there are set theoretic universes
in which some theories, with infinite models, do not have saturated models. The situation is even worse if you want countable
saturated models. By joining computability with saturation, Barwise and Schlipf, and independently Ressayre,\footnote{For a more
detailed account on the history of recursive saturation see \cite[\S 3]{Barwise.Schlipf:76}.} 
came up with the notion of recursive saturation in the seventies.

Every theory, in a recursive language, with an infinite model has a lot of recursively saturated model. More importantly,
any such theory has a countable recursively saturated model.
It turned out that countable recursively saturated models behave very much like saturated models. A lot of techniques used
in classic model theory could now be adopted to, for example, first-order arithmetic, or Peano arithmetic as it is often called. 

However, recursive saturation is mostly a useful notion only for countable models. The slightly stronger notion of resplendence 
seems to work better with uncountable models. A model is resplendent if any $\Sigma_1^1$-formula which is consistent with
the theory of the model is in fact true in the model. By a theorem of Kleene this is equivalent to that if $T$ is a theory
in a recursive extension of the language of the model which is consistent with the theory of the model, then there
is an expansion of the model satisfying the theory. For countable models resplendence and recursive saturation coincide.
The notion of resplendence was, also, introduced by Barwise and Schlipf.

For long it seemed that recursive saturation was a very rigid notion. For example, recursively enumerable saturation coincide
with recursive saturation. Also resplendent and recursive saturation coincide for the most interesting case of countable models. 
This view has changed since the work of Kaye, Kossak and Kotlarski, \cite{Kaye.Kossak.ea:91}, where they find an interesting
variation which is strictly stronger than recursive saturation. They call it arithmetic saturation and they characterise
all countable models of arithmetic which are arithmetic saturated in terms of properties on the automorphism group of the
model, both as a permutation group and as a topological group. Later Lascar proved that the automorphism group of such a model of arithmetic
has the small index property, thus reducing the property of the automorphism group as a topological group to a property of the 
group as an abstract group.\footnote{Thus there are countable recursively saturated models of arithmetic with different automorphisms
groups. This was unknown before the paper of Kaye, et. al. and the paper by Lascar.} 

It is in this context the second chapter of this thesis should be read. There we introduce several new variations of
recursive saturation and resplendence. 

The notion of recursive saturation, and its cousins, is very tightly intervened with coding properties. The standard system
of a model is the collection of all sets of natural numbers coded in the model. It can be seen as a measure on the degree
of saturation of the model.
For countable models Scott characterised all algebras of sets of natural numbers occurring as standard systems. By
doing a limit construction this actually also work for models of cardinality $\aleph_1$, however the construction 
can not, to our knowledge, be taken further. Very little is known about which algebras of cardinality larger than $\aleph_1$ 
are realized as standard systems. There will be more on this in the last chapter of this thesis.

The thesis consists of four chapters. The first one is an introduction to some background material we will need in the other three
chapters. In it we present an overview of the literature and most of the results are not ours. Some small, and easy, remarks and 
are, however, to our knowledge, new.
The next three chapters are new, except possibly for the first part of the third chapter.

\section*{First chapter, background}\addcontentsline{toc}{section}{First chapter, background}
 
Various background material needed for the rest of the thesis is covered in the first chapter. Most of it comes
from models of arithmetic, but there are also some material from descriptive set theory and second-order arithmetic. Some proofs
are presented, and others are not. We have chosen to include proofs which are important for the rest of the thesis or which are particularly nice. 

Most of the material is known, but some small comments seems to be new. For example the notion of low saturation and the easy
propositions \ref{prop:arit.imp.xsat} and \ref{prop:low.imp.xsat} are new. Also, we have not been able to
find the proof of Theorem \ref{thm:scott.sets.are.ssy} explicitly in the literature even though we understand it is well-known.

\section*{Second chapter, expansions omitting types}\addcontentsline{toc}{section}{Second chapter, expansions omitting types}

The second chapter is the main one of this thesis. It constitutes of several proposed generalisations of
the notion of resplendence, an expandability property of structures. 
A model is resplendent if for all recursive theories $T$, in a language recursively extending the language of the model, which
is consistent with the theory of the model, there is an expansion of the model satisfying it.

The proposed notions try to generalise this to, not only satisfying a theory, 
but also omitting a type. The most naive generalisation would be, for all recursive $T$ and $p(\bar x)$ if
$T+p\om+\Th(M)$, where $p\om$ is a non first-order sentence expressing that $p(\bx)$ is omitted, has a model
then there is an expansion of $M$ satisfying $T+p\om$. Trivially, this generalisation is too strong, since if
$p(\bx)$ is a type realized in $M$ such that $\Th(M)+p\om$ is consistent, then no expansions of $M$ could ever satisfy
$p\om$.  We propose two different ways of weakening this naive notion:
\begin{enumerate}
\item\label{weak1} Strengthening the consistency assumption on $T+p\om$.
\item\label{weak2} Weakening the expandability conclusion of $M$.
\end{enumerate}

In the first section of this chapter the first proposed solution is discussed. The resulting notion is called \emph{transcendence}, and
we prove that, for countable models, enough saturation implies transcendence; and that, for models of $\PA$, transcendence implies 
quite a lot saturation. However, we have not been able to tie up the loose ends completely and prove an equivalence between
transcendence and a saturation property.

The second section discusses what happens if we restrict the types we can omit in expansions to limit types, i.e., types
which are not isolated. It turns out that arithmetic saturation is enough to prove this form of transcendence, at least for
countable models; and we conjecture that, for models of arithmetic, it is also necessary.

We go on and have a look at expansions to theories $T+p\om$ which are categorical, in the sense that for a given $T_0$ for
every $M \models T_0$ there is at most one expansion satisfying $T+p\om$. The main example of such a theory is
the theory $\Kom$ which expresses that a unary predicate $K$ is the predicate of standard numbers, it is categorical 
over $\PA$. Some general results about such theories are proven with the help of transcendence; these results inspired
the next section.

The fourth section discusses the special case when $T+p\om$ extends $\Kom$ and only one constant symbol is added to 
the language; this notion only applies to models of $\PA$.
We prove that the resulting notion, called \emph{recursive standard saturation} is in fact equivalent, for countable models of $\PA$,
to a saturation property.

Lastly we discuss the second way of weakening the naive notion of expandability. The conclusion is now, not that there is
an expansion satisfying a theory and omitting a type, but that there is an elementary submodel having such an expansion. It turns out that
this notion, called \emph{subtranscendence}, is equivalent, for models of $\PA$, to $\beta$-saturation. A model of $\PA$ is
$\beta$-saturated if its standard
system is a $\beta$-model, i.e., if $\SSy(M) \embin_{\Sigma^1_1} \power(\omega)$, where $\SSy(M)$ and $\power(\omega)$ 
are interpreted as $\omega$-models of second-order arithmetic.

All results presented in this chapter are new, and we think they show that there are other saturation properties
apart from recursive and arithmetic saturation that are interesting.

\section*{Third chapter, the theory of omitting types}\addcontentsline{toc}{section}{Third chapter, the theory of omitting types}

Given a theory $T$ and a type $p(\bx)$ what is the theory, $\Th(T+p\om)$, of all models of $T$ omitting $p(\bx)$? This question is
answered in the third chapter by defining a new inference rule that, schematically, looks like this:

\begin{equation}
  \tag{$p$-rule}
  \inference{\lnot \exists \bar x \varphi(\bar x)}{\inferforall{i \in \omega}{\forall \bar x (\varphi(\bar x) \imp
  p_i(\bar x))}}
\end{equation}
where $p(\bx)=\set{p_i(\bx) | i \in \omega}$, i.e., we may deduce $\lnot \exists \bar x \varphi(\bar x)$ if we can deduce $\forall \bar x (\varphi(\bar x) \imp
p_i(\bar x))$ for all $i \in \omega$.

We also apply some of the theory developed in the second chapter to prove that there is a type and a theory with rank
$\omega_1^\mathrm{CK}$, where the rank of a theory and a type is a measure on how isolated the type is by the theory. 
If the type $p(\bx)$ is isolated by $T$ then the rank of $p(\bx)$ over $T$ is $0$, if not then the rank measures how
many times the $p$-rule has to be used to get the theory $\Th(T+p\om)$.

The proof
system discussed in this chapter is already implicit in a paper
by Casanovas and Farr\'e, \cite{Casanovas.Farre:96}. However, our approach is somewhat different, and we think
that most of results are new.

\section*{Fourth chapter, standard systems}\addcontentsline{toc}{section}{Fourth chapter, standard systems}

Scott's problem is to characterise the standard systems of models of $\PA$. For countable models, and for models of cardinality 
$\aleph_1$, this has been done. In the case of the continuum hypothesis this settles the problem. However, if the continuum 
hypothesis fails very little is known about the problem for models of cardinalities greater than $\aleph_1$. 

We prove that, under the assumption of Martin's axiom,
every Scott set of cardinality $<2^{\aleph_0}$ closed under
arithmetic comprehension which has the
countable chain condition is the standard system of some model of
$\PA$. However, we do not know if there
exists any such uncountable Scott sets.

The ultraproduct construction in this chapter is strongly 
inspired by one of Kanovei's papers, \cite{Kanovei:96}, where he, given a 
countable arithmetically closed set $\scott{X}$,
constructs a model $M$ of true arithmetic with $\SSy(M)=\scott{X}$ and such that a set $A \subseteq \omega$ is representable (without parameters) over $(M,\omega)$ by a $\Sigma_k$-formula iff it is definable (without parameters) 
by a $\Sigma^1_k$ formula over $\scott{X}$.

\chapter{Background}

We start this chapter by presenting some notation and definitions. Then we go on with some words on how to 
arithmetise logic inside first-order arithmetic, this yields the arithmetised completeness theorem. The important notions
of recursive saturation and resplendence are presented in sections \ref{sec:rec.sat} and \ref{sec:resp}. Scott sets 
and $\scott{X}\!$-saturation are also presented together with arithmetic saturation. We also introduce the notion
of low saturation. We end the chapter with some results from second-order arithmetic and descriptive set theory. 

\section{Notation and preliminaries}

Most of the notation and definition used are taken from \cite{Kaye:91}. 
For clarity we repeat most of them here.

Languages will mostly be countable and recursive and denoted by $\La$. 
The theory of a model $M$ is
\[
\Th(M)=\set{\varphi | M \models \varphi \text{ and $\varphi$ is a sentence}};
\]
if $T$ is a first-order theory then 
\[
\Th(T)=\set{ \varphi | T \prf \varphi \text{ and $\varphi$ is a sentence}}.
\] 
The underlying language will, hopefully, be clear from the context.
If $\bar a, \bar b \in M$ then
\begin{multline*}
\tp_M(\bar a/\bar b)=\{\varphi(\bar x,\bar b) \mathbin| M \models \varphi(\bar 
a,\bar b),\text{and $\varphi(\bar x,\bar y)$ is a formula} \\
\text{ with all free variables shown}\}.
\end{multline*}
We write $\tp_M(\bar a)$ to mean $\tp_M(\bar a/\emptyset)$, and $\tp(\bar a /\bar 
b)$ if the model $M$ is understood from the context. Observe that $\tp_M(\bar 
a/\bar b)$ is the same set of formulas as $\tp_{(M,\bar b)}(\bar a)$, and that 
$\tp_M(\emptyset/\bar b)$ is $\Th(M,\bar b)$. 

Given 
any first-order theory $T$, a \emph{complete type over $T$} is a set $\tp_M(\bar a)$ 
for some $\bar a \in M \models T$, and a \emph{type over $T$} is a subset of a complete 
type over $T$. The set of all complete types over $T$ with $k$ free variables 
is denoted $S_k(T)$, i.e., 
\[
S_k(M)=\set{\tp_M(\bar a/\bar b) | \bar a,\bar b \in M, |\bar a|=k, |\bar b|<\omega}.
\]

A complete type over a model $M$ is a set $\tp_N(\bar 
a/\bar b)$ for some $\bar b \in M \embin N$ where $\bar a \in N$, and a type 
over $M$ is a subset of a complete type over $M$. The set of all complete types 
over $M$ with $k$ free variables is denoted $S_k(M)$. 

The three sets 
$S_k(\Th(M))$, $S_k(M)$, and $S_k(\ThM)$ differs only in how many parameters we 
allow the types to have; in the first no parameters are allowed, in the second 
we allow finitely many parameters in each type, and in the third arbitrarily many 
parameters are allowed. 

$\PA$ is full first-order arithmetic in the ordinary language of arithmetic, $\La_A= \set{<,+,\cdot,0,1}$. 
$\PA^-$ is $\PA$ but without the induction scheme.


\section{Arithmetising logic}

Let us fix some standard G\"odel numbering of formulas and terms, and identify a syntactic object with its 
G\"odel number. Thus; it makes sense saying that a theory is, for example, recursive. Sometimes we will, however,
write the G\"odel number of a formula $\varphi(\bar x)$ as $\godel{\varphi(\bar x)}$ if clarity is gained.

We will use some machinery for coding sequences of 
elements by a single element in a model of $\PA$. The details for constructing such a coding will not
be carried out, see for example \cite{Kaye:91}. The standard notation $(x)_y$ 
will be used for the $y$th element coded by $x$. A set $A \subseteq \omega$ is 
said to be coded in a model $M$ of $\PA$ if there is $a \in M$ such that $A 
=\codedset_M(a)$, where
\[
\codedset_M(a)= \set{ k \in \omega  |  M \models (a)_k \neq 0 }.
\]
The standard system of a model $M$ of $\PA$ is
\[ 
\SSy(M) = \set{ \mathrm{set}_M(a)  |  a \in M}.
\]
Given a finite set $S \subset \omega$, the least
standard natural number coding $S$ will be denoted $[S]$.

A set $X \subseteq M$ is also said to be coded in $M$ if there is $a \in M$ such that
$X = \set{ m \in M | M \models (a)_m \neq 0}$. This terminology might be slightly confusing. Since
a set $A \subseteq \omega$ is coded, in this second sense, iff $A$ is finite, we hope that
the reader accepts this abuse of terminology. All coded sets, again in this second sense, are definable
and bounded, and, by using the induction axiom, any definable bounded set is coded.

We will assume the existence of formulas 
$\Prf_\varphi^z(x,y)$ which, in the standard way, enumerates the relation 
of
\begin{quote}
``$x$ is a proof of the formula $y$ using the non logical axioms  
$\set{z | \varphi(z)}$.''
\end{quote}
See any good textbook on the arithmetisation of logic for the details.
The formula $\Prf_\varphi^z(y)$ is a convenient 
short-hand for $\exists x\Prf_\varphi^z(x,y)$, and $\con_\varphi^z$ is the 
formula $\exists x \lnot \Prf_\varphi^z (x)$. When the free variable of $\varphi$ 
is easily understood from the context we will omit the superscript $z$. If 
$a$ is an element of some model $M$ the formula $\con_{(a)_x \neq 0}^x$ is denoted by $\con_a$.

\begin{defin}
Given a model $M \models \PA$ and another model $N$ in some recursive language
we say that \emph{$N$ is strongly interpreted in $M$} if the domain of $N$ is a subset 
of $M$, and there are formulas $\mathrm{dom}(x)$ and $\mathrm{sat}(x,y)$, with parameters from $M$,
such that the domain of $N$ is 
\[
\set{a \in M | M \models \mathrm{dom}(a)}
\]
and for any 
$\bar a \in N$, any $b \in M$ coding the sequence $\bar a$, and any 
$\varphi(\bar x)$ in the language of $N$ the following holds
\[
N \models \varphi(\bar a) \quad \text{iff} \quad M \models 
\mathrm{sat}(\varphi,b).
\]
\end{defin}

By proving the completeness theorem inside $\PA$ you get the following:

\begin{thm}[Arithmetised Completeness Theorem]\label{thm:arit.comp}
If $M \models \PA$ and $T\in\SSy(M)$ is a consistent theory in a 
recursive language such that $M \models 
\con_\tau$ for some $\tau(x,a)$, $a \in M$, enumerating $T$ in $M$; 
then there exists a model of $T$ strongly interpreted in $M$.
\end{thm}
\begin{proof}
Do the ordinary Henkin construction of the completeness theorem, but
this time inside $M$. See \cite[Section 13.2]{Kaye:91} for the details.
\end{proof}

In fact, by an overspill argument, we can replace the assumption $M \models 
\con_\tau$ in the theorem by the assumption that 
\begin{equation}\label{eq:con}
 \text{for all finite } S\subseteq T \text{ we have } M \models \con_{[S]},
\end{equation}
where $[S]$ is the least standard natural number coding the finite set $S$:
If \eqref{eq:con} holds and $\tau(x,a)$ enumerates $T$ 
then 
\[
M \models \con_{\tau(x,a) \land x<k}
\] 
for all $k\in\omega$, and so 
by overspill 
\[
M \models \con_{\tau(x,a) \land x<b}
\]
 for some nonstandard $b \in 
M$.

This last argument gives a formula, with parameters from $M$, enumerating $T$ 
such that $M$ thinks the theory defined by the formula is consistent. 
If we assume $T$ to be recursive there is a 
more uniform way of doing this: the Feferman construction \cite[Theorem 
5.9]{Feferman:60}. It gives a formula (without parameters) that works in any 
model of \[\PA + \set{\con_{[S]} | S \subseteq T \text{ and $S$ is finite}},\] 
i.e., a formula $\tau(x)$ such that
\[
\PA + \set{\con_{[S]} | S \subseteq T  \\ \text{ and $S$ is finite}} \models \con_\tau.
\]
It should be noted that even though it is assumed that $T$ is recursive
the formula $\tau(x)$ is, in general, not $\Sigma_1$.

\section{Recursive saturation and homogeneity}\label{sec:rec.sat}

Throughout this section, and most of the thesis, we will assume the 
base language $\La$ to be recursive. When not specifying anything else $M$ 
will be a structure in $\La$. A type $p(\bar x,\bar a)$ over $M$ is said to be 
recursive if the set
\[
\set{ \godel{\varphi(\bar x,\bar y)} | \varphi(\bar x,\bar a) \in p(\bar x,\bar 
a)} \subseteq \omega
\]
is recursive.

\begin{defin}
A model $M$ is \emph{recursively saturated} if all recursive types over $M$ are 
realized in $M$.
\end{defin}

Recursively saturated models behave very well as we will see; and they exist in 
abundance:
\begin{prop}
Any model $M$ (in a recursive language) has an elementary extension $N$ which is 
recursively saturated and such that $|M| = |N|$.
\end{prop}
\begin{proof}
The proof is by a union of chains argument.
\end{proof}

A closely related notion is homogeneity.

\begin{defin}
A model $M$ is called \emph{$\omega$-homogeneous} if for all $\bar{a},\bar{b} 
\in M$ satisfying $\tp(\bar a)=\tp(\bar b)$, and for all $c \in M$ there is a 
$d \in M$ such that $\tp(\bar a,c)=\tp(\bar b,d)$. $M$ is said to be \emph{strongly homogeneous} 
if for every $\bar{a},\bar{b} \in M$ such that $\tp(\bar 
a)=\tp(\bar b)$ there is an automorphism of $M$ taking $\bar a$ 
to $\bar b$.
\end{defin}

For countable models these two notions of homogeneity 
coincide; which is proved by a straight forward back-and-forth argument. 

Recursive  saturation implies $\omega$-homogeneity:
\begin{prop}
Any recursively saturated model is $\omega$-homogeneous.
\end{prop}
\begin{proof}
Let $\bar a,\bar b \in M$ be such that $\tp_M(\bar a)=\tp_M(\bar b)$, 
and let $c \in M$. Define
\[
p(x) = \set{\varphi(\bar a,c) \imp \varphi(\bar b,x) | \varphi(\bar y,x) \text{ 
a formula}};
\]
then $p(x)$ is a recursive type over $M$ since for any $\varphi(\bar y,x)$, if $M 
\models \varphi(\bar a,c)$ then $M \models \exists x \varphi(\bar a,x)$ and so 
$M \models \exists x \varphi(\bar b,x)$, which, by compactness,  proves that 
$p(x)$ is a type over $M$. Any $d \in M$ realizing $p(x)$ has the property that 
 $\tp_M(\bar a,c)=\tp_M(\bar b,d)$.
 \end{proof}
 
Combining the last proposition 
with the fact that any countable $\omega$-homogeneous model is strongly homogeneous
we get that any countable recursively saturated model 
is strongly homogeneous.

\section{Scott sets}

\begin{defin}
A set $\scott{X} \subseteq \power(\omega)$ is a \emph{Scott set} if it is
non-empty and
\begin{itemize}
\item if $A,B \in \scott{X}$ then $A \cup B$, $\omega \setminus A \in
  \scott{X}$,
\item if $A \in \scott{X}$ and $B \leq_T A$, i.e., $B$ is recursive in
  $A$, then $B \in \scott{X}$,
\item if $T \in \scott{X}$ is a theory in some recursive language then
  there is a complete consistent theory $T_c \in \scott{X}$  extending $T$.
\end{itemize}
\end{defin}

The definition says that a Scott set is a boolean algebra of sets of
natural number which is closed under  relative
recursiveness, and completing theories.

The last clause in the definition of a Scott set could be replaced by
saying that if $T \in \scott{X}$ is an infinite binary tree (coded as a set of 
natural numbers), then there is an infinite path $P \in \scott{X}$ through $T$. \footnote{The alternative 
definition is the more common way of defining a Scott set; but we think that our definition makes more sense in this
setting.}

There are recursive theories $T$ without recursive completions (for example 
$\PA$), so the set of recursive sets is \emph{not} a Scott set.
By using the low basis theorem we can find a low completion $T_c$ of any
theory $T$, meaning that ${T_c}'\leq_T T'$, i.e., that the jump of the
completion is recursive in the jump of the theory. Therefore; the set 
$\LOW$ of all low sets, i.e., $A \in \LOW$ iff $A' \leq_T \emptyset'$, 
is a Scott set. 
It is important to note that there are Scott sets not including $\LOW$ as a subset; in fact,
for each non-recursive $A \subseteq 
\omega$ there is a Scott set $\scott{X}$ such that $A \notin \scott{X}$.

\begin{defin}
A model $M$ is said to be \emph{$\scott{X}$-saturated}, where $\scott{X}$ is a Scott set, 
if for every complete type $p(\bar x, \bar a)$ over $M$ we have
\[
\exists \bar b \in M \models p(\bar b, \bar a) \quad\text{iff}\quad p(\bar x,\bar y) \in \scott{X}.
\]
\end{defin}

We will sometimes be a bit sloppy and write $p(\bar x,\bar a) \in \scott{X}$ 
instead of $p(\bar x,\bar y) \in \scott{X}$. In fact; this sloppiness is quite alright 
since, if we extend the G\"odel numbering to the parameters $\bar a$, 
we have that $p(\bar x,\bar a)$ and $p(\bar x, \bar y)$ are Turing equivalent;
thus, one of them is in a given Scott set iff the other one is. 

Since models of $\PA$ admits coding of finite sequences we can replace 
finite sequences of variables in types by a single variable in the following \emph{
recursive} way: Given a type $p(\bar x,\bar a)$ define
\begin{multline*}
\leftfix p'(x,b) = \{ \forall \bar y, \bar z \bigl(y_0=(x)_0 \land \ldots \land y_{k-1}=(x)_{k-1} \land {} \\ 
z_0=(b)_0 \land \ldots \land z_{l-1}=(b)_{l-1} \imp \varphi(\bar y,\bar z) \bigr) 
\mathbin| \varphi(\bar x,\bar a) \in p(\bar x,\bar a)\}
\end{multline*}
where $b$  codes the finite sequence $\bar a$, and $k$ and $l$ are the lengths of $\bar x$ 
and $\bar a$ respectively. For any model of $\PA$, $p(\bar x,\bar a)$ is 
realized exactly when $p'(x,b)$ is realized.

\begin{prop}
If $M$ is  $\scott{X}$-saturated and $p(\bar x,\bar a) \in \scott{X}$ is a, not necessarily complete, 
type over $M$, then $p(\bar x,\bar a)$ is realized in $M$.
\end{prop}
\begin{proof}
Let $q(\bar c,\bar a) \in \scott{X}$ be a consistent completion of $p(\bar c,\bar a)$, where $\bar c$ are new constant symbols. 
Then $q(\bar x,\bar y) \in 
\scott{X}$ and so the type $q(\bar x,\bar a)$ is realized which clearly implies that 
$p(\bar x,\bar a)$ is realized.
\end{proof}

Since any Scott set includes all 
recursive sets we get that any $\scott{X}\!$-saturated model is recursively 
saturated.

\begin{prop}\label{prop:standard.system.is.scott}
If $M \models \PA$ is nonstandard then $\SSy(M)$ is a 
Scott set. Also, if $M$ is recursively saturated then $M$ is
$\scott{X}\!$-saturated iff $\scott{X}=\SSy(M)$.
\end{prop}
\begin{proof}
The proof that $\SSy(M)$ is a Scott set is straight forward but a bit lengthy, so we skip 
it here.

Let $M \models \PA$ be recursively saturated and $p(x,a) \in \SSy(M)$ a 
complete type over $M$. Then there is $b \in M$ such that 
$\codedset_M(b)=p(x,y)$. Let $q(x,a,b)$ be the recursive set
\[ 
\set{\varphi(x,a) \imp (b)_n \neq 0 | \varphi(x,y) \text{ a formula and } n=\godel{\varphi(x,y)}}. 
\]
It is a type since $p(x,a)$ is a type; thus, by the recursive saturation, it is 
realized by some $m \in M$. Now, $M \models p(m,a)$ since if $\varphi(x,a) \in 
p(x,a)$ then $M \models (b)_n = 0$, where $n=\godel{\neg \varphi(x,y)}$, and so $M \nmodels \neg 
\varphi(m,a)$, i.e., $M \models \varphi(m,a)$. This proves that $M$ realizes all 
complete types in $\SSy(M)$.

To prove that it does not realize anything outside 
of $\SSy(M)$ let $a,b \in M$. Then
\[ 
p(x,a,b)=\set{\varphi(a,b) \ekv 
(x)_n \neq 0 | \varphi(x,y)\text{ a formula and } n=\godel{\varphi(x,y)}}
\]
is a recursive type and so realized. Any $m \in M$ realizing $p(x,a,b)$ 
satisfies
\[
\codedset_M(m)=\set{\varphi(x,y) | \varphi(x,b) \in \tp_M(a/b)}.
\] 
Thus $tp_M(a/b) \in \SSy(M)$ and we conclude that $M$ is $\SSy(M)$-saturated.

For the other direction suppose $M$ is $\scott{X}\!$-saturated. Let $A \in \scott{X}$ and 
\[ 
p(x) = \set{(x)_k \neq 0 | k \in A} \cup \set{(x)_k=0 | k \notin A}. 
\]
The type $p(x)$ is recursive in $A$ and so is in $\scott{X}$ which implies that 
$p(x)$ is realized. Any $m \in M$ 
realizing $p(x)$ satisfies $\codedset_M(m)=A$. Thus; $\scott{X} \subseteq 
\SSy(M)$. 

If $A \in \SSy(M)$ let $a \in M$ be such that $\codedset_M(a)=A$. Since 
$M$ is $\scott{X}\!$-saturated $\tp_M(a) \in \scott{X}$ and $A$ is recursive in 
$\tp_M(a)$, so $A \in \scott{X}$. All this proves that if $M$ is 
$\scott{X}\!$-saturated then $\scott{X}=\SSy(M)$.
\end{proof}

\begin{thm}[\cite{Scott:62}]\label{thm:scott}
For any consistent recursive theory $T$ (in the language of arithmetic) 
extending $\PA$ and any $\scott{X} \subseteq \power(\omega)$ the following are equivalent:
\begin{itemize}
\item[$(i)$] $\scott{X}$ is a countable Scott set.
\item[$(ii)$] There is a countable nonstandard model $M \models T$ with  $\SSy(M)=\scott{X}$.
\end{itemize}
\end{thm}
\begin{proof}
The implication  $(ii) \Rightarrow (i)$ is Proposition \ref{prop:standard.system.is.scott}. The other
implication is a Henkin construction.
\end{proof}

This theorem could be taken a bit further: Given a Scott set $\scott{X}$ of cardinality 
$\aleph_1$ there is a nonstandard model of $T$ with standard system $\scott{X}$. 
To prove this strengthening we need a definition and a lemma:

\begin{defin}
Given a theory $T$, in the language of arithmetic, 
we say that a set $A \subseteq \omega$ is \emph{represented in $T$} (sometimes called strongly represented)
if there is a formula $\varphi(x)$ such that for all $n \in A$ we have $T \prf \varphi(n)$ and for all $n \in \omega \setminus A$ we have
$T \prf \lnot \varphi(n)$. By $\rep(T)$ we denote the collection of all sets represented in $T$.
\end{defin}

\begin{lem}
Let $T$ be a consistent complete theory in the language of arithmetic 
extended by countably many constant symbols.
If $\rep(T) \subseteq \scott{X}$, where $\scott{X}$ is a countable Scott set,
then there is a countable model $M \models T$ such that $\SSy(M) = \scott{X}$.
\end{lem}
\begin{proof}
We use a Henkin style argument to construct the model $M \models T$. Let $\La= \La_A \cup \set{c_i}_{i \in \omega}$ 
be the language of $T$, $D=\set{d}_{i \in \omega}$ be new constants,  $\varphi_k(x)$ be an enumeration of all $\La(D)$-formulas with
one free variable, and $\set{X_i}_{i \in \omega}$ an enumeration of $\scott{X}$. 

We construct consistent theories $T_i \subseteq T_{i+1}$ 
such that $X_i \in \rep(T_{i+1})$, $T_i \setminus T \in \scott{X}$, 
and if $T_{i+1} \prf \exists x \varphi_i(x)$
then $T_{i+1} \prf \varphi(d_j)$ for some $j \in \omega$. Let $T_0=T$ and, given $T_i$, 
if $T_i \prf \exists x \varphi_i(x)$ let $T_{i+1}$ be
\begin{multline*}
T_i + \varphi_i(d_j) +
\set{ (d_j)_k \neq 0 | k \in X} + \set{(d_j)_k =0 | k \notin X} \\
+ \set{ (d_l)_k \neq 0 | k \in X_i} + \set{(d_l)_k =0 | k \notin X_i}
\end{multline*}
where $d_j$ and $d_l$ are constants not occurring in $T_i$ or in $\varphi_i(x)$, and $X \in \scott{X}$. To find $X$ such that
$T_{i+1}$ is consistent let $\Th^m_n(T)$ be the set of all $\La_A \cup \set{c_0,\ldots, c_m}$ sentences 
which are $\Sigma_n$ and provable by $T$. Even though $T$ might not be in $\scott{X}$ the set $\Th^m_n(T) \in \scott{X}$
for all $m,n \in \omega$ since it is represented in $T$ (by a truth-definition for $\Sigma_n$-formulas).
Let $n$ be such that all sentences in $\set{\varphi_i(x)} \cup T_i \setminus T$ is 
$\Sigma_n$ and let $m$ be large enough so that no $c_i$, where $i \geq m$,
occurs in $T_i \setminus T$ or in $\varphi_i(x)$. By construction, the set $T_i \setminus T$ is in $\scott{X}$ so the theory
$T_i \setminus T + \Th^m_n(T) + \varphi(d_j)$ 
is in $\scott{X}$. Let $S \in \scott{X}$ be any consistent completion of that theory, and let $X \in \scott{X}$ be 
such that
\[
n \in X \quad\text{iff}\quad (d_j)_n\neq 0 \in S.
\]
It should be clear that this $X$ makes the theory $T_{i+1}$ consistent.

If $T_i \nprf \exists x \varphi_i(x)$ let $T_{i+1}$ be
\[
T_i + \lnot \exists x \varphi_i(x) +
\set{ (d_l)_k \neq 0 | k \in X_i} + \set{(d_l)_k =0 | k \notin X_i},
\]
where $d_l$ does not occur in $T_i$ or $\varphi_i(x)$. 

The union of all $T_i$ is a consistent complete Henkin theory; its term model has standard system $\scott{X}$.
\end{proof}

To prove that any Scott set of cardinality $\aleph_1$ is realized as a standard system we use a 
union of chains argument. Let $\scott{X}=\set{X_\alpha}_{\alpha < \omega_1}$ and $\scott{X}_0 \subset \scott{X}$ be
any countable Scott set. Moreover; let $M_0 \models T$ be a countable model with standard system $\scott{X}_0$.
Define $M_\alpha$ and $\scott{X}_\alpha$, for $\alpha \leq \omega_1$, by induction: Given $\scott{X}_\alpha$ and 
$M_\alpha$, let $\scott{X}_{\alpha+1}$ be a countable Scott set such that 
\[
\scott{X}_\alpha \cup \set{X_\alpha} \subseteq \scott{X}_{\alpha+1} \subseteq \scott{X}.
\]
Since $\SSy(M_\alpha) \subseteq \scott{X}_{\alpha+1}$ we have that $\rep(\Th(M_\alpha,a)_{a \in M_\alpha})
\subseteq \scott{X}_{\alpha+1}$ and so, by the lemma, there is a model $M_{\alpha+1}$ 
of $\Th(M_\alpha,a)_{a \in M_\alpha}$ with standard system
$\scott{X}_{\alpha+1}$. Clearly $M_\alpha \embin M_{\alpha+1}$. 
For limit ordinals $\lambda \leq \omega_1$ let 
$M_\lambda$ be the union of $\set{M_\alpha}_{\alpha<\lambda}$ and $\scott{X}_\lambda$ the union 
of $\set{\scott{X}_\alpha}_{\alpha<\lambda}$.
It is should be clear that $M_{\omega_1} \models T$ and $\SSy(M_{\omega_1})=\scott{X}_{\omega_1}=\scott{X}$.
We have proved the following theorem.

\begin{thm}\label{thm:scott.sets.are.ssy}
Given a consistent completion $T$ of $\PA$ and a Scott set $\scott{X}$
of cardinality $\aleph_1$ or $\aleph_0$, there is a model of $T$ with
$\scott{X}$ as its standard system.
\end{thm}

The problem of characterising standard systems in general is known as Scott's problem. 
The theorem above solves it completely if we assume the CH to hold. 
Very little is known about it otherwise, but see
Chapter \ref{ch:scott} for a small result in this direction.

Let us now use the arithmetised completeness theorem to give a sufficient condition for when
a theory $T$ has a $\scott{X}\!$-saturated model, for some given Scott set $\scott{X}$ of cardinality $\aleph_0$ or
$\aleph_1$.

\begin{thm}[{\cite[Theorem 2.29]{Wilmers:75}}]\label{thm:scott.model.existence}
If $\scott{X}$ is a Scott set, $|\scott{X}|\leq \aleph_1$ and $T \in \scott{X}$ is a
consistent theory then there is an $\scott{X}\!$-saturated
model of $T$.
\end{thm}
\begin{proof}
Let $T' = \PA + \set{\con_{[S]} | S\subseteq T \text{ is finite}}$.
By Theorem \ref{thm:scott.sets.are.ssy}
there is a nonstandard model $M \models T'$ such that $\SSy(M) =
\scott{X}$. Let $t \in M$ code $T$ and
by overspill let $a \in M \setminus \omega$ be such
that 
\[
M \models \con_{(t)_x \neq 0 \land x<a}.
\]
By the arithmetised completeness theorem, Theorem \ref{thm:arit.comp}, there is a model 
$N \models T$ strongly interpreted in $M$. It should now be easy to see that $N$ is 
$\scott{X}\!$-saturated.
\end{proof}

If $T$ is complete this condition is also necessary: If $M$ is $\scott{X}\!$-saturated then
$\Th(M) \in \scott{X}$. 

To sum up we have;
\begin{itemize}
\item any $\SSy(M)$, where $M \models \PA$, is a Scott set,
\item any Scott set $\scott{X}\!$, where $|\scott{X}|\leq \aleph_1$, is equal to $\SSy(M)$ for some $M \models \PA$, and
\item if $M \models \PA$ and $T \in \SSy(M)$ is a consistent theory, then there is a $\SSy(M)$-saturated
model of $T$. 
\end{itemize}

We will later need the next proposition which says that there is only one (up to isomorphism) countable
$\scott{X}\!$-saturated model.

\begin{prop}\label{prop:scott.uniqueness}
If $M$ and $N$ are two countable models of the same complete theory
which both are $\scott{X}\!$-saturated, then they are isomorphic.
\end{prop}
\begin{proof}
Easy back-and-forth.
\end{proof}

\section{Arithmetic saturation}

\begin{defin}
A model $M$ is \emph{arithmetically saturated} if every type $p(\bar x,\bar a)$ 
over $M$ which is arithmetic in some realized type $\tp_M(\bar b)$ is realized in $M$.
\end{defin}

\begin{thm}[{\cite[Proposition 5.2]{Kaye.Kossak.ea:91}}]\label{thm:k.k.k}
For countable recursively saturated models $M \models\PA$ the following are 
equivalent:
\begin{enumerate}
\item $M$ is arithmetically saturated,
\item For any $f \in M$ there is $c \in M$ such that $M \models f(k) > n$ for 
all $n \in \omega$ iff $M \models f(k) > c$, $k \in 
\omega$.
\item\label{it:maxaut} There exists $g \in \Aut(M)$ such that 
$\mathrm{fix}(g)=\set{a \in M | g(a)=a}=\dcl(\emptyset)$.
\end{enumerate}
\end{thm}

Observe that \eqref{it:maxaut} could be expressed as realizing a theory and 
omitting a type in a bigger language: Let $T$ say that $g$ is an automorphism 
and let $p(x)$ say that $x$ is not definable and $g(x) \neq x$.

We could extend the list with properties of the automorphism group of $M$, either as a permutation group, a 
topological group, or an abstract group. Therefore, there are two 
countable recursively saturated models of $\PA$ with non-isomorphic automorphism groups.

Arithmetically saturated models are easier to handle than recursively saturated models as the 
following theorem shows, which is unknown to be true if arithmetic saturation is replaced by recursive saturation. Observe 
that for models of $\PA$ this theorem reduces to Proposition \ref{prop:standard.system.is.scott}.

\begin{prop}\label{prop:arit.imp.xsat}
Every arithmetically saturated model is $\scott{X}\!$-saturated for some Scott set 
$\scott{X}$.
\end{prop}
\begin{proof}
Let $\scott{X}$ be the arithmetic closure of $\set{\tp(\bar a) | \bar a \in 
M}$. Clearly $\scott{X}$ is a Scott set (since it is the arithmetic 
closure of something) and for all $\bar a, \bar b \in M$, 
the type $\tp(\bar a/\bar b)$ is in $\scott{X}$ since 
$\tp(\bar a/\bar b) \in \scott{X}$ by definition means $\tp(\bar a,\bar b) \in \scott{X}$. Let
$p(\bar x,\bar y) \in \scott{X}$. Since if $A$ is arithmetic in $B$ which is
arithmetic in $C$ then $A$ is arithmetic in $C$ we have that $p(\bar x,\bar y)$
is arithmetic in some $\tp(\bar b)$. Thus, by the definition of arithmetic
saturation, if $p(\bar x,\bar a)$ is a type then it is realized.
\end{proof}

The same is true for recursively saturated models of cardinality at most $\aleph_1$
as we will prove in the next section, but the proof is more complicated and it is
unknown if it holds for models of greater cardinalities. 

In fact, we do not need
full arithmetic saturation to prove the last theorem.
\begin{defin}
A model $M$ is \emph{low saturated} if every type $p(\bar a,\bar x)$ over $M$ which is low in
some $\tp_M(\bar b)$ is realized in $M$.
\end{defin}

\begin{prop}\label{prop:low.imp.xsat}
Any low saturated
model is $\scott{X}\!$-saturated for some Scott set
$\scott{X}$.
\end{prop}
\begin{proof}
Let $\scott{X}$ be the set of sets low in some
$\tp_M(\bar a)$ where $\bar a \in M$. If $A$ is low in $B$, which, in turn, 
is low in $C$; then $A$ is low in $C$. Thus; the proof above also works for low
saturated models.
\end{proof}

\section{Resplendence}\label{sec:resp}

In some sense recursive saturation is a sort of expandability property: For any
recursive theory $T$ in a language expanded by finitely many constants and finitely many
parameters $\bar a \in M$ which is
consistent with $\Th(M,\bar a)$ there is an expansion of $M$ satisfying $T$. Here is a
version of recursive saturation, introduced by Barwise and Schlipf,\footnote{Their
definition is slightly different and needs a theorem of Kleene on the expressibility of $\Sigma_1^1$-formulas 
to prove that it implies recursive saturation.} along those lines:

\begin{defin}
An $\La$-model $M$ is
resplendent if for any $\Ba \in M$, any recursive extension $\La^+$ of $\La
(\Ba)$ and any recursive $T$ in $\La^+$ consistent with $\Th(M,\bar a)$ 
there exists an expansion $M^+$ of $M$ satisfying $T$.
\end{defin}

The next theorem shows that, in the countable case, recursive saturation implies 
resplendence. However, first we would like to prove something much easier: 
Any countable \emph{low} saturated model is resplendent:

Let $M$ be 
countable and $\scott{X}\!$-saturated, where $\scott{X}$ is a Scott set such that if
$A$ is low in $B \in \scott{X}$ then $A \in \scott{X}$. Let $\bar a \in M$ and $\La^+$ be a recursive
extension of $\La(\bar a)$, and let $T$ be a recursive 
theory consistent with $\Th(M,\bar a)$. By Theorem 
\ref{thm:scott.model.existence} there is a countable $\scott{X}\!$-saturated
model $N$ of $\Th(M,\bar a)+T$ since $\Th(M,\bar a)+T \in \scott{X}$. 
Since both $M$ and $N$ are countable and $\scott{X}\!$-saturated, by using Theorem 
\ref{prop:scott.uniqueness}, the restriction of $N$ to the language of $M$ is 
isomorphic to $M$ with an isomorphism taking the interpretation of $\bar a$ in $N$ to $\bar a$ in $M$. 
Thus; $M$ has an expansion satisfying $T$.

\begin{thm}\label{thm:rec.ekv.resp}
A countable 
$\La$-model is resplendent iff it is recursively saturated. In fact, if $M$ is 
countable and recursively saturated and $T$ is as in the definition of 
resplendence, then there is a recursively saturated expansion of $M$ satisfying
$T$.
\end{thm}
\begin{proof}
The easy direction is to prove that any resplendent 
model is recursively saturated:

If $M$ is resplendent, then given a recursive
type $p(\bx,\bar a)$, $\bar a \in M$, let $T=p(\bar c,\bar a)$ where $\bar c$ 
are new constants. Clearly $\Th(M,\bar a) +T$ is consistent so by resplendence 
there is an expansion $M \redof M^+\models T$. The elements of $M$ 
interpreting $\bar c$ realizes the type $p(\bar x,\bar a)$. Since $p(\bar x,\bar 
a)$ was chosen arbitrarily, any recursive type is satisfied, i.e., $M$ is 
recursively saturated.

For the other implication let $\bar a \in M$ and $T$ be a 
recursive theory in a recursive extension $\La^+$ of $\La (\bar a)$ 
such that $\Th(M,\bar a) + T$ is consistent. Also; let $\set{\varphi_i(x)}$ be an 
enumeration of all $\La^+(M)$ formulas with exactly one free variable. We will 
build finite $\La^+(M)$-theories $T_k$ such that $T+T_k+\ThM$ is consistent and 
if $T_{i+1} \prf \exists x \varphi_i(x)$ then for some $m \in M$ we have 
$T_{i+1} \prf \varphi_i(m)$. 

We start off by letting $T_0=\emptyset$. Assume $T_i$ 
has been defined and satisfies the above properties; define $T_{i+1}$ in the
following way: Let $\bar b \in M$ be all parameters occurring in $T_i, \varphi_i(x)$ and
$\bar a$.
If
\[
T+T_i+\Th(M,\bar b) \nprf \exists x \varphi_i(x)
\]
let $T_{i+1} = T_i + \lnot \exists x \varphi_i(x)$, else, if there is $b \in \bar b$ such that
\[
T+T_i+\Th(M,\bar b) \nprf \neg \varphi_i(b),
\]
let
$T_{i+1}=T_i+\varphi(b)+\exists x \varphi_i(x)$.

Otherwise; let
\[
p(x,\bar b) = \set{\theta(x,\bar b) \in \La(\bar b) | T+T_i
\prf \forall x \bigl(\varphi_i(x) \land x \neq \bar b \imp \theta(x,\bar b)\bigr)},
\]
where
$x \neq \bar b$ is a shorthand for $\bigwedge_{b \in \bar b} x \neq b$.
The set $p(x,\bar b)$ is
recursively enumerable, so by Craig's trick (see \cite[p. 150]{Kaye:91}) it is
logically equivalent to a recursive set. To prove that it is a type let
$\theta(x,\bar b) \in p(x,\bar b)$ then, since
\[
T+T_i+\Th(M,\bar b) \prf \exists x\bigl(\varphi_i(x) \land x \neq \bar b\bigr) \land 
\forall x \bigl(\varphi_i(x) \land x \neq \bar
b \imp \theta(x,\bar b)\bigr),
\]
we have
\[
T+T_i+\Th(M,\bar b) \prf \exists x \theta(x,\bar b)
\]
and so $M \models \exists x \theta(x,\bar b)$. This shows that $p(x,\bar b)$ is
a type over $M$ since $p(x,\bar b)$ is closed under conjunctions.
Let $m \in M \models p(m,\bar b)$; if
\[
\ThM+T+T_i+\varphi_i(m) \prf \perp
\]
then there is $M \models \theta(m,\bar b)$ such that
\[
T+T_i \prf \varphi_i(m) \imp \lnot\theta(m,\bar b).
\]
Since $m$ does not appear in $T$ or $T_i$ (since $x \neq \bar b \in p(x,\bar
b)$) we have
\[
T+T_i \prf \forall x \bigl(\varphi(x) \imp \lnot\theta(x,\bar b)\bigr),
\]
and so $\lnot\theta(x,\bar b) \in p(x,\bar b)$ which contradicts the fact that
$m$ satisfies $p(x,\bar b)$.

Therefore $T_{i+1}=T_i+\varphi_i(m)$ makes the theory
$\ThM+T+T_{i+1}$ consistent.

Any Henkin model of $\ThM+T+\cup_i T_i$ is a model of $T$ whose $\La$-reduct is
isomorphic to $M$.

To prove that the model can be taken to be recursively saturated we need to introduce
satisfaction classes. We will not do it here; see \cite[Theorem
15.8]{Kaye:91} for the details, and \cite{Engstrom:02} for more on satisfaction classes.
\end{proof}

The theorem is not true in the uncountable case since any
$\omega_1$-like, i.e., of cardinality $\aleph_1$ but with all its
proper initial segments countable, recursively saturated model $M$ of $\PA$ is
not resplendent: Let $T$ state that $f$ is an isomorphism from $M$ onto a proper
initial segment of $M$, $T$ is consistent with $\Th(M)$ by Friedman's theorem
\cite[Section 4]{Friedman:73} (or \cite[Theorem 12.4]{Kaye:91}), but clearly
there is no expansion of $M$ satisfying $T$. A recursively saturated $\omega_1$-like 
model of $\PA$ can be found by using Friedman's theorem ones more, see 
\cite[Page 247]{Kaye:91} for the details.

The next proposition says that in the countable case all 
recursively saturated models are  $\scott{X}\!$-saturated for some Scott set 
$\scott{X}$. It is not known if this is true for higher cardinalities, except 
for the case of models of cardinality $\aleph_1$, in which a union of chains 
argument will prove that any recursively saturated model is 
$\scott{X}\!$-saturated for some $\scott{X}$. If, in addition, the continuum 
hypothesis holds it is not very hard to prove that any recursively saturated 
model $M$ is $\scott{X}\!$-saturated for some Scott set $\scott{X}$: 
By a 
downward L\"owenheim-Skolem argument let $M_0 \embin M$ be recursively saturated 
and realizing exactly the same complete types as $M$ does. $M_0$ can be chosen 
of cardinality $\leq \aleph_1$ and so there is a Scott set $\scott{X}$ for which 
$M_0$ is $\scott{X}\!$-saturated. Clearly, $M$ is also 
$\scott{X}\!$-saturated.
See Theorem \ref{prop:arit.imp.xsat} for more in this direction.

\begin{prop}
Every countable and recursively saturated model is 
$\scott{X}\!$-saturated for some Scott set $\scott{X}$.
\end{prop}
\begin{proof}
Let $M$ be a countable recursively saturated model in $\La$,  by Theorem 
\ref{thm:rec.ekv.resp}, it is resplendent and so let $M \redof M^+ \models \PA$, 
where the arithmetic language of $\PA$ is disjoint with $\La$, be 
recursively saturated. Let $\scott{X}=\SSy(M^+)$; then $M^+$ is 
$\scott{X}\!$-saturated and so for any complete type $p(\bar x,\bar a)$ over 
$M$, $p$ is realized in $M$ iff $p(\bar x,\bar y) \in \scott{X}$. This means 
that also $M$ is $\scott{X}\!$-saturated.
\end{proof}

We can now prove something slightly stronger than Theorem \ref{thm:rec.ekv.resp}.

\begin{prop}
Let $M$ be a 
countable recursively saturated model in the language $\La$, $\scott{X}$  
such that $M$ is $\scott{X}\!$-saturated, $\bar a \in M$, and $T \in \scott{X}$ 
 a theory in a recursive extension of $\La(\bar a)$ such that 
$\Th(M,\bar a) + T$ is consistent. Then there is an \emph{$\scott{X}\!$-saturated}
expansion $M \redof M^+ \models T$.
\end{prop}
\begin{proof}
Since $\scott{X}$ is 
a Scott set we have $\Th(M,\bar a) + T\in \scott{X}$. By Theorem 
\ref{thm:scott.model.existence} there is a countable $\scott{X}\!$-saturated 
model of $\Th(M,\bar a) + T$ and by Theorem \ref{prop:scott.uniqueness} the 
restriction of this model to $\La$ is isomorphic to $M$.
\end{proof}

We finish off 
this section by proving that resplendent models are very symmetric in the sense 
of homogeneity.

\begin{prop}
Any resplendent model is strongly homogeneous.
\end{prop}
\begin{proof}
Let $\bar a,\bar b \in M$ be such that 
$\tp_M(\bar a)=\tp_M(\bar b)$. Let $\La^+$ be $\La \cup \set{g}$ where $g$ is a 
unary function symbol. Let $T$ be
\[ 
\set{ \forall \bar x\bigl(\varphi(\bar x) \ekv 
\varphi(g(\bar x))\bigr) | \varphi(\bar x) \in \La}  +  \forall 
x\exists y \bigl(x=g(y)\bigr) + g(\bar a)=\bar b,
\]
where $g(\bar x)$ is 
$g(x_0),g(x_1),\ldots,g(x_k)$ if $\bar x$ is $x_0,x_1,\ldots,x_k$. We have to 
prove that $\Th(M,\bar a,\bar b) + T$ is consistent. Let $N$ be any countable 
recursively saturated model of $\Th(M,\bar a,\bar b)$, then $N$ is 
$\omega$-homogeneous and countable so it is strongly homogeneous. Also 
$\tp_N(\bar a)=\tp_N(\bar b)$ and so there is an automorphism $f\in \Aut(N)$ of 
$N$ taking $a$ to $b$, and so $(N,f) \models T$. By resplendence there 
is an expansion of $M$ satisfying $T$. Any realization of the function symbol $g$ is an 
automorphism of $M$ taking $\bar a$ to $\bar b$.
\end{proof}

\section{Second-order arithmetic}

We will now shortly discuss second-order arithmetic theories, i.e,  theories 
with two types of variables, one for numbers and one for sets of numbers. The set variables 
will be written with capital letters and the number variables with lower case
letters. We add a new type of atomic formula: $x \in X$, which should be interpreted as 
`the number $x$ is a member of the set $X$'. Models of second-order arithmetic have two domains, 
the number domain is the ordinary domain and the set domain which is a subset of 
the power set of the number domain. The semantics for such models are evident.
It should be noted that these theories
are expressible 
in first-order logic. Second-order logic would have implied that the set domain was the full power set 
of the number domain.

Full second-order arithmetic, denoted $Z_2$, is 
$\PA^-$ ($\PA$ without the induction axioms) plus the second-order induction axiom
\begin{equation}\label{eq:ind}
\forall X \bigl(0 \in X \land \forall x (x \in X \imp x+1 \in X) \imp \forall x (x \in X)\bigr),
\end{equation}
and 
the comprehension scheme
\[
\exists X \forall x \bigl(x \in X \ekv \theta(x)\bigr),
\]
where $\theta(x)$ is any formula in which $X$ does not occur freely 
($\theta(x)$ may have other free variables, both set 
and number variables).

A formula is called arithmetic, or first-order, 
if it has no bound set variables (it may still have free set variables). The 
ordinary arithmetic hierarchy of formulas ($\Sigma_k$, $\Pi_k$ and $\Delta_k$) 
extends to set quantifiers ($\Sigma^1_k$, $\Pi^1_k$ and $\Delta^1_k$). For example, 
a formula
\[
\exists X_1 \forall X_2 \ldots Q X_k \varphi(X_1,X_2,\ldots,X_k,\bar Y,\bar x)
\]
is $\Sigma_k^1$ if $Q$ is $\exists$ when $k$ is odd and $\forall$ otherwise, and
$\varphi$ a first-order formula. 
Hence; the set of arithmetic formulas can be written as $\Delta^1_0$.

An $\omega$-model of $Z_2$ (or some of its subsystems) is a model
\[ 
(\omega,\scott{X},<,+,\cdot,0,1) 
\]
where $\scott{X} \subseteq \power(\omega)$. 
We often specify an $\omega$-model only by the set domain $\scott{X}$. $\Nat_2$ is the full 
$\omega$-model of $Z_2$, i.e., $\power(\omega)$.

$\RCA_0$ is the subsystem of 
$Z_2$ consisting of $I\Sigma_1$ ($\PA^-$ plus the first-order induction scheme 
for $\Sigma_1$-formulas), the second-order induction axiom \eqref{eq:ind}, and the 
$\Delta_1$-comprehension scheme, i.e.,
\[
\exists X \forall x (x \in X \ekv \theta(x))
\]
for all $\Delta_1$-formulas $\theta(x)$, where $X$ does not occur in $\theta$ and 
$\theta(x)$ may have other free variables.

$\WKL_0$ is $\RCA_0$ plus an axiom 
saying that any coded infinite tree has a coded infinite path. The 
$\omega$-models of $\WKL_0$ are precisely the Scott sets, or to be more precise
the set domains of $\omega$-models of $\WKL_0$ are exactly the Scott sets.

$\ACA_0$ is $Z_2$ but 
with the comprehension scheme restricted to arithmetic formulas.

We say that an 
$\omega$-model $\scott{X} \subseteq \power(\omega)$ is a \emph{$\beta_k$-model} if 
for any $\Sigma^1_k$-sentence $\Theta$ with set-parameters from $\scott{X}$, we 
have $\scott{X} \models\Theta$ iff $\Nat_2 \models \Theta$.  A \emph{$\beta$-model} 
is a $\beta_1$-model ($\beta$-models were first studied by 
Mostowski in \cite{Mostowski:61}).

It could be worth mentioning that if we are 
working over full second-order arithmetic $Z_2$ a $\beta$-model is just a model 
$M$ for which well-orderings are absolute. To be more precise, if $\prec \in M$ 
is any non well-ordering on $|M|$ then there exists a witness for the 
non well-ordering of $\prec$ in $M$, i.e.,the following holds
\[ 
M \models \exists f \forall x (f(x+1) \prec f(x)). 
\]
Clearly well-orderings are absolute in any $\beta$-model. To prove the other implication 
we need that any $\Sigma_1^1$-formula is equivalent to saying that a certain tree has an infinite 
branch. This is, in turn, equivalent to saying that the Brouwer-Kleene ordering 
of the specific tree is a well-ordering. For the details of the proof see \cite[Theorem 
1.11]{Apt.Marek:73}.

A set $\scott{X}$ is a \emph{$\beta_\omega$-model} if it is
a $\beta_k$-model for all $k \in \omega$. Observe that $\scott{X}$ is a 
$\beta_\omega$-model iff $\scott{X} \embin \power(\omega)$ as second-order $\omega$-models.

\begin{defin}
If $\Delta$ is a collection of sets of natural numbers and $\Gamma$ a collection of subsets of $\power(\omega)$ we 
say that $\Delta$ is a \emph{basis for $\Gamma$} if for any $\gamma \in \Gamma$ we have $\gamma \cap \Delta \neq \emptyset$.  
\end{defin}

That $\scott{X}$ is a $\beta_\omega$-model is then equivalent to that $\scott{X}$ is a basis for the collection 
of subsets of $\power(\omega)$ definable by a second-order formula with parameters from $\scott{X}$.

If $\bar A \subseteq \omega$, then $\Sigma_k^{1,\bar A}$  will denote both the collection of sets of 
natural numbers definable in $\Nat_2$
by an $\Sigma^1_k$-formula $\theta(x,\bar A)$, and the collection of subsets of $\power(\omega)$ definable by 
a $\Sigma^1_k$-formula $\theta(X,\bar A)$. The sets $\Pi_k^{1,\bar A}$ and $\Delta_k^{1,\bar A}$ are defined 
similarly.

Under G\"odels set-theoretic assumption $\VisL$, saying that every set is constructible, 
see \cite[Chapter 13]{Jech:03}, we have the following:

\begin{thm}[{\cite[Corollary V.2.7]{Hinman:78}}]\label{thm:basis.for.beta}
Assume $\VisL$ holds; then $\Delta_k^{1,A}$ is a basis for 
$\Sigma_k^{1,A}$ for all $k \geq 2$ and $A \subseteq \omega$.
\end{thm}
\begin{proof}
The proof is by defining a
well-ordering of $\power(\omega)$ which is $\Delta_2^1$ and then choose the 
smallest element according to that well-ordering. This would, at first sight, give us that $\Pi_{k+1}^{1,A}$ is a basis 
for $\Sigma_k^{1,A}$, but you can take this further and prove the theorem.
\end{proof}

Similarly, under the assumption of Projective Determinacy, or $\PD$ for short (see \cite[Chapter 33]{Jech:03}), 
we have:

\begin{thm}[{\cite[Corollary V.3.6]{Hinman:78}}]\label{thm:basis.for.beta.PD}
Assume $\PD$ holds; then $\Delta^{1,A}_k$ is a 
basis for $\Sigma^{1,A}_k$ for all even $k \geq 2$ and all $A \subseteq \omega$.
\end{thm}

Together these last two theorems give us:

\begin{cor} \label{cor:basis.for.beta}
Assume that either $\VisL$ or $\PD$ holds, then for all $A \subseteq \omega$ 
the set $\Delta^{1,A}_\infty$ is a basis for itself.
\end{cor}

If $\Gamma$ is any class of formulas then a set $\scott{X} \subseteq \power(\omega)$ 
is said to satisfy \emph{$\Gamma$-comprehension} if
\[ 
\scott{X} \models  \forall \bar 
Y\exists X \forall n \bigl(n \in X \ekv \theta(n,\bar Y)\bigr)
\]
for all formulas $\theta \in \Gamma$, where $X$ does not occur freely. 
The set $\scott{X}$ is said to satisfy \emph{true $\Gamma$-comprehension} if for all $\theta 
\in \Gamma$ and all $\bar{Y} \in \scott{X}$ there is an $X \in \scott{X}$ such that
\[
\Nat_2 \models \forall n \bigl(n \in X \ekv \theta(n,\bar{Y})\bigr).
\]

\begin{prop}
Let $\scott{X}$ be a $\beta_k$-model, $k \geq 1$, then:
\begin{enumerate}
\item If $\Theta(x,\bar X)$ is a $\Sigma_k^1$-formula, then there is a $\Sigma^1_k$-formula $\Psi(\bar X)$
such that both $\scott{X}$ and $\Nat_2$ satisfies
\[
\forall \bar X \bigl( \forall x (\Theta(x,\bar X) \ekv \Psi(\bar X)\bigr).
\]
\item $\scott{X}$ satisfies $\Sigma^1_k$-comprehension iff it satisfies true $\Sigma^1_k$-comprehension.
\item $\scott{X}$ satisfies $\Delta^1_k$-comprehension and true $\Delta^1_k$-comprehension.
\end{enumerate}
\end{prop}
\begin{proof}
The first statement is easily seen to be true for $\Nat_2$ (see for example \cite[Theorem 16.III]{Rogers:87}).
On the other hand assume $\scott{X} \models\forall x \Theta(x,\bar A)$ for some $\bar A \in \scott{X}$. Then 
$\Nat_2 \models \forall x \Theta(x,\bar A)$ and so $\Nat_2 \models \Psi(\bar A)$. We can conclude that 
$\scott{X} \models \Psi(\bar A)$. The other direction is equally easy.

For the second statement all we have to do is to observe that for any $\beta_k$-model $\scott{X}$
\[
\Nat_2 \models \theta(n,\bar A) \quad\text{iff}\quad \scott{X} \models \theta(n,\bar A)
\]
for all $\Sigma^1_k$-formulas $\theta(x,\bar A)$, where $\bar A \in \scott{X}$.

By statement two, to prove the third one, it is enough to prove that $\scott{X}$ satisfies
$\Delta^1_k$-comprehension.

We want to prove that for all $\Delta_k^1$-formulas $\Theta(\bar X)$ and all 
$\bar{A} \in \scott{X}$ there is a $\Sigma_{k}^1$-formula $\Phi(\bar X)$ such that
\[ 
\scott{X} \models \Phi(\bar A) \ekv \exists Y \forall n (n \in Y \ekv \Theta(n,\bar A)) 
\]
and $\Nat_2 \models \Phi(\bar A)$; in that case $\scott{X} \models \Phi(\bar A)$ and so
\[ 
\scott{X} \models \exists Y \forall n (n \in Y \ekv \Theta(n,\bar{A})). 
\]
Since $\Theta(\bar X)$ is $\Delta^1_k$ it is also $\Delta^1_k$ in $\scott{X}$ (remember that $\scott{X}$ is a
$\beta_k$-model). Therefore it is not hard to see that $x \in Y \ekv \Theta(x,\bar X)$ is $\Sigma^1_k$. By using statement one
we see that $\forall x (x \in Y \ekv \Theta(x,\bar X)$ is also $\Sigma^1_k$ proving the statement.
\end{proof}

\begin{cor}\label{cor:beta.comprehension}
Assume either $\VisL$ or $\PD$ and $k \geq 2$ is even, 
then a set $\scott{X} \subseteq \power(\omega)$ is a $\beta_k$-model iff
it satisfies true $\Delta^1_k$-comprehension.
\end{cor}
\begin{proof}
The right-to-left implication is just a combination of Theorem \ref{thm:basis.for.beta} and \ref{thm:basis.for.beta.PD}. For the
other direction let $\bar A \in \scott{X}$ and $\Nat_2 \models \Theta(x,\bar A) \ekv \Psi(x,\bar A)$, where
$\Theta$ is a $\Sigma^1_k$-formula and $\Psi$ a $\Pi^1_k$ ditto. Clearly 
\[
\scott{X} \models \Theta(x,\bar A) \ekv \Psi(x,\bar A),
\]
and so the formula $x \in X \ekv \Theta(x,\bar A)$ is equivalent to a $\Sigma^1_k$-formula. And thus, by the proposition,
\[
\exists X \forall x \bigl(x \in X \ekv \Theta(x,\bar A)\bigr)
\]
is also equivalent to a $\Sigma^1_k$-formula. Therefore it holds in $\scott{X}$, 
since it holds in $\Nat_2$. 

The $\omega$-model $\scott{X}$ therefore satisfies $\Delta^1_k$-comprehension; by the proposition it also satisfies true $\Delta^1_k$-comprehension.
\end{proof}

\begin{cor}
Any $\beta_\omega$-model satisfies both 
$\Sigma_\infty^1$-comprehension and true $\Sigma_\infty^1$-comprehension.
\end{cor}

We end this chapter by a recent theorem of Mummert and Simpson, which we will need later.

\begin{thm}[{\cite[Corollary 3.7]{Mummert.Simpson:04}}]\label{thm:exists.beta.not.beta}
For any $k\in\omega$ 
there is a countable $\beta_k$-model which is not a $\beta_{k+1}$-model.
\end{thm}

The proof is by showing that for each recursively enumerable second-order arithmetic theory $T$, if there
is a $\beta_k$-model $T$ then there is a $\beta_k$-model of $T+\text{``there is no $\beta_k$-model of $T$''}$,
a sort of incompleteness theorem for $\beta_k$-models, and observe that any $\beta_{k+1}$-model of $T$ satisfies
``there is a $\beta_k$-model of $T$''.


\chapter{Expansions omitting types}\label{ch:tr}

In this chapter we will investigate when a model can be expanded to satisfy certain 
non first-order theories in bigger languages. First, we discuss the notion of \emph{transcendence}. 
A model is transcendent
if for all theories of the form $T+p\om$, where $T$ is a first-order theory, $p(\bar x)$ is a type, and
$p\om$ is the (non first-order) sentence expressing that $p(\bar x)$ is omitted, is consistent, in a rather strong special sense,
then there is an expansion of the model satisfying $T+p\om$. Transcendence has some connections with 
saturation properties of the model, for example it implies some very strong saturation.

We will then go on and investigate what happens if we restrict ourselves to non-isolated, or limit, types $p(\bar x)$.
It turns out that that notion is much weaker than transcendence.

Some theories of the form $T+p\om$ has the property that for any $M \models T_0$ there is at most one expansion of $M$ satisfying $T+p\om$. Such theories $T+p\om$ are said to be \emph{categorical over $T_0$}. For those theories we can say quite a bit more about the expansions satisfying it. This is dealt with
in the section following. As a special case we get the standard predicate, which tells you if an element of a model of arithmetic is standard or not. If we expand a model of arithmetic with the standard predicate the resulting structure is clearly not 
recursively saturated. We define a version of recursive saturation that works for such structures. We also characterise all such models
in terms of ordinary saturation. 

We end the chapter with the notion of \emph{subtranscendence}, which is like transcendence, but instead of concluding that there
is an expansion of the model in question, it predicts that the model has an elementary submodel with such an expansion. It turns out
that a model of arithmetic is subtranscendent iff it is $\beta$-saturated. The last proposition of the
chapter is a characterisation of $\beta$-models in terms of closure under completing certain non first-order theories.

\section{Transcendence}

Given a set $p(\bar{x})$ of first-order formulas in a countable language
$\La$ we  denote the $\La_{\omega_1 \omega}$-sentence
\[
\forall \bar{x} \bigvee_{\psi(\bar{x}) \in p(\bar{x})} \neg \psi(\bar{x})
\]
by $p\om$, where the arrow binds all free variables of $p(\bar x)$. The sentence expressing that 
$p(\bar x)$ is realized, i.e., $\lnot p\om$ is denoted by $p\re$. Observe that even though $p\re$ is not a first-order
sentence it is equivalent to realizing a first-order theory.

A transcendent model is a model which is resplendent in, not only first-order logic, but first-order logic
plus the option of omitting a type. More precisely, but see Definition \ref{def:trans} for the exact definition, 
a strongly homogeneous $\scott{X}\!$-saturated model 
$M$ is transcendent if the language $\La$ 
of $M$ is recursive and for every recursive extension $\La^+$ of $\La$ and every $T,p(\bar x) \in \scott{X}$ 
such that there exists a model $N \models \Th(M)+T+p\om$, where $N \restrictedto \La$ is $\omega$-saturated; there is an expansion of $M$ satisfying $T+p\om$.

The property is of the schematic form:
\begin{quote}
for all $T,p(\bx) \in \scott{X}$ such that $\con(\Th(M),T,p(\bar x))$ there is an 
expansion of $M$ satisfying $T+p\om$. 
\end{quote}
In this case $\con(T_0,T,p(x))$ is taken to be $\satcon(T+p\om/T_0)$ which means 
that there exists a model of $T_0+T+p\om$ whose $\La$-part 
is $\omega$-saturated, where $\La$ is the language of $T_0$. 
This is a rather strong condition but none weaker \emph{seems} to work. 
In particular the notion of ordinary consistency is too weak as the following shows:

Let $M$ be a non-standard model of true arithmetic, i.e., $M \equiv \Nat$. 
There exists a realized recursive type $p(x)$ such that
$\Th(M) +p\om$ is consistent; let $p(x)$ be $\set{x>0,x>1,\ldots}$.  
Clearly, $M$ realizes $p(x)$ and so no expansion of $M$ satisfies $p\om$. 

Another possibility would be to take $\con(T_0,T,p(x))$ to hold iff 
for all types $q(\bar x)$ over $T_0$ the theory $T_0+T+p\om + q\re$ is consistent.
Later in this section we will see why this does not work.

A natural question to ask is if it is possible 
to strengthen the express power of the logic to full $\La^+_{\omega_1\omega}$ 
to get a property 
like: if $\Theta$ is a $\La^+_{\omega_1\omega}$ sentence which is somehow 
consistent with $M$ then there exists an expansion of $M$ satisfying 
$\Theta$? The simple answer is no, but of course this depends on the consistency 
notion in use. The following very well-known theorem by Dana Scott will show that the answer is negative if 
the consistency notion is taken to be $\satcon$.

\begin{thm}[\cite{Scott:65}]
For any countable model $N$ in a 
countable language $\La$ there is an $\La_{\omega_1 \omega}$-sentence $\sigma_N$ 
such that any other model in the same language $\La$ is back-and-forth equivalent 
with $N$ iff it satisfies $\sigma_N$; in particular any countable model satisfying $\sigma_N$ is
isomorphic to $N$.
\end{thm}

Let $M \models T$ be any 
countable model of a non $\omega$-categoric complete theory $T$ (in a countable 
language $\La$) and let $N$ be any other countable model of $T$, i.e., $M \ncong N$. An 
$\omega$-saturated model $K$ of $T$ is back-and-forth equivalent with $N$ so $K 
\models \sigma_N$, and so $\satcon(\sigma_N/\Th(M))$. However, $M \nmodels 
\sigma_N$, since for countable models back-and-forth equivalence and isomorphism coincide.

\begin{que}
Is there a fragment $\La'$ of $\La_{\omega_1 \omega}$ such that any countable model has a 
countable elementary extension having expansions satisfying all $\La'$-theories $T$ such that 
$\satcon(T/\Th(M))$? 
\end{que}

\subsection{Basic definitions}

Since we are investigating the property of omitting
a type it is not important that the type is consistent with the base
theory. Therefore, when we say that $p(\bar x)$ is a type, we only mean that 
$p(\bar x)$ is a set of formulas with $\bar x$ as the only free variables. However,
that $p(\bar x,\bar a)$ is type \emph{over} a model $M$ (or a theory $T$)
still means that $p(\bar x)$ is consistent with $\Th(M,\bar a)$ (or $T$).

Consistency will always mean consistency in the semantical sense, i.e., 
a theory is consistent iff there is a model satisfying it. This is important since we are dealing
with theories which are not first-order. We will later on, in Chapter \ref{ch:proof}, deal with syntactic
notions of consistency. 

We write $M^+ \expof M$ or $M \redof M^+$ to mean that $M^+$ is an expansion of $M$.

Let $T_0$ be a first-order theory in a recursive language $\La$ and $S$ a theory
(not necessarily first-order) in a recursive extension $\La^+$ of $\La$. 
\begin{defin}
The property $\satcon(S/T_0)$ holds if there exists a model $N \models T_0+S$ such
that the $\La$-reduct of $N$ is $\omega$-saturated. 
\end{defin}
Strictly speaking the
smaller language $\La$ should be indicated when we write $\satcon(S/T_0)$, but it is
usually understood to be the language of $T_0$ (mostly $T_0$ is complete
in $\La$). For first-order theories $S$ $\satcon$ is the ordinary first-order consistency, i.e., 
$\satcon(S/T_0)$ holds iff $S+T_0$ is consistent.

Let us now define the main concept of this chapter.

\begin{defin}\label{def:trans}
A strongly homogeneous model $M$ is called \emph{transcendent} if there is a Scott
set $\scott{X}$ for which $M$ is $\scott{X}\!$-saturated and for all recursive
extensions $\La^+$ of the language of $M$, all first-order $\La^+$-theories $T \in
\scott{X}$, and all types $p(\bar x) \in \scott{X}$ also in $\La^+$ such
that $\satcon(T+p\om/\Th(M))$ there is $M \redof M^+ \models
T+p\om$.
\end{defin}

That we require the model to be strongly homogeneous is a way of getting rid of mentioning parameters 
in the definition. The following proposition shows that.

\begin{prop}\label{prop:trans.imp.param.trans}
If $M$ is transcendent and $\bar{a} \in M$ then
$(M,\bar{a})$ is transcendent.
\end{prop}
\begin{proof}
Given that
$\satcon(T+p\om/\Th(M,\bar a))$ it clearly follows that 
\[
\satcon(T+S+p\om/\Th(M)),
\]
where $S$ is $\Th(M,\bar a)$ but with the parameters $\bar a$ replaced by some new constants $\bar c$.
By the transcendence of $M$ let $M'$ be an expansion of $M$ satisfying $T+S+p\om$. 
The problem now is that $\bar c$ might not be interpreted in $M'$ as
$\bar a$, the homogeneity of $M$ helps us fix this.

Let $\bar b \in M$ be the
interpretation of $\bar c$ in $M'$. Clearly $\tp_M(\bar a)=\tp_M(\bar b)$ and
so by the strong homogeneity of $M$ there is an $f \in \Aut(M)$ such that
$f(\bar b )=\bar a$. Let $M^+$ be $f(M')$, i.e., 
\begin{align*}
M^+ \models P(\bar m) &\quad\text{iff}\quad M' \models P(f^{-1}(\bar m))\\
\intertext{for any predicate symbol $P$ and any $\bar m \in \M$,} 
M^+ \models g(\bar m)=n &\quad\text{iff}\quad M' \models g(f^{-1}(\bar m))=f^{-1}(n) 
\intertext{for any function symbol $g$ and any $\bar m, n \in M$, and}
M^+ \models c=m &\quad\text{iff}\quad M' \models c=f^{-1}(m)
\end{align*}
for any constant symbol $c$.
The mapping
$f:M' \to M^+$ is an isomorphism, so $M^+ \models T+p\om$ and the
interpretation of $\bar c$ in $M^+$ is $\bar a$. Furthermore; 
$M^+$ is an expansion of $M$ since $M^+\restrictedto \La$ is $M$.
\end{proof}

Transcendence implies resplendence which follows directly form the last proposition; still we formulate it as a proposition.

\begin{prop}
Any transcendent model $M$ is resplendent.
\end{prop}
\begin{proof}
Let $\bar a \in M$, $\La^+$ and $T$ be as in the definition
of resplendence, i.e., $\La^+$ is an extension of 
$\La(\bar a)$ and $T$ is a recursive first-order $\La^+$-theory  
consistent with $\Th(M,\bar a)$. Clearly
$\satcon(T/\Th(M,\bar a))$ and since $(M,\bar a)$ is transcendent there is an
expansion of $(M,\bar a)$ satisfying $T$.
\end{proof}

Clearly; if $M$ is transcendent and $T, p(\bar x), q(\bar x) \in \scott{X}$ 
are such that 
\[
\satcon(T+p\om+q\om/\Th(M))
\]
then there is an expansion of $M$ 
satisfying $T+p\om+q\om$. This is easily seen since 
\[
\models p \om \land q \om \ekv r \om
\]
where 
\[
r(\bar x) = \set{ \varphi(\bar x) \vee \psi(\bar x) | \varphi(\bar x) \in p(\bar x), 
\psi(\bar x) \in q(\bar x)}.
\]

\subsection{Countable transcendent models exist}

The next proposition gives us some easy facts about $\satcon$; it will help us to build
transcendent models.

\begin{prop}\label{prop:satcon.properties}
Let $\La^+$ be an extension of a language $\La$, $T$ a first-order theory in $\La^+$, 
$T_0$ a first-order theory in $\La$ and $p(\bar x)$ an $\La^+$-type. Assume $\satcon(T+p\om/T_0)$.
\begin{enumerate}

\item\label{item:ett} If $\sigma$ is a sentence in the language $\La^+$ then either \[\satcon(T+\sigma+p\om/T_0), \text{ or }\satcon(T+\lnot\sigma+p\om/T_0).\]

\item\label{item:tva} If $T_0$ is complete and $q(\bar x)$ a type over $T_0$ (i.e., consistent with $T_0$) then 
\[\satcon(T+p\om/T_0+q(\bar c))\] for any constant symbols $\bar c$ not in $\La^+$.

\item\label{item:tre} If $\bar c$ are constants in $\La^+$ then there is $\psi(\bar x) \in p(\bar x)$ such that
\[\satcon(T+\lnot \psi(\bar c)+p\om/T_0).\]

\item\label{item:fyra} If $T \prf \exists x \varphi(x)$ then there is a complete type $q(x)$ over
$T_0$ and a new constant symbol $c$ such that \[\satcon(T+\varphi(c)+p\om/T_0+q(c)).\]
\end{enumerate}
\end{prop}
\begin{proof}
Let $N$ witness that $\satcon(T+p\om/T_0)$, i.e., $N \restrictedto \La$ is $\omega$-saturated and
$N \models T+p\om + T_0$.
\begin{enumerate}
\item Either $N \models \sigma$ or $N \models \lnot \sigma$, and so either $\satcon(T+\sigma+p\om/T_0)$ or $\satcon(T+\lnot\sigma+p\om/T_0)$.
\item Since $q(\bar x)$ 
is consistent with $T_0$ and $T_0=
\Th(N)$ there is $\bar b \in N$ realizing $q(\bar x)$. The model $(N,\bar b)$ witnesses that $\satcon(T+p\om/T_0+q(\bar c))$. 
\item Let $\bar b \in N$ be the interpretation of $\bar c$. Since $N \models p\om$ there has to be $\psi(\bar x) \in p(\bar x)$ such that $N \models \lnot\psi(\bar b)$ and so $\satcon(T+\lnot \psi(\bar c) +
p\om /T_0)$. 
\item Let $b \in \N$ be such that $N \models \varphi(b)$, and let $q(x)=\tp_N(b)$. The model $(N,b)$ witnesses 
that $\satcon(T+\varphi(c)+p\om/T_0+q(c))$.\qedhere
\end{enumerate}
\end{proof}

These are the only properties of $\satcon$ we will be using. 
Observe that it is only \eqref{item:tva} which is not true for ordinary consistency, it is the property making 
$\satcon$ work for us.

Let us now find a saturation property which, for countable model, implies transcendence. The following definition says that
a Scott set $\scott{X}$ is $\satcon$-closed if we can complete theories, using the consistency notion $\satcon$, inside $\scott{X}$.

\begin{defin}
We say that a Scott set $\scott{X}$ is \emph{$\trans$-closed} if for every
language $\La$, any extension $\La^+$ of $\La$, and any $T,T_0,p(\bar x) \in \scott{X}$ 
such that $\satcon(T+p\om/T_0)$ 
there exists a complete $\La^+$-theory $T_c \in \scott{X}$ such that $T \subseteq T_c$ 
and $\satcon(T_c+p\om/T_0)$.
\end{defin}

Clearly, by a union of chains argument, any Scott set lies inside a $\trans$-closed Scott set 
of the same cardinality.

\begin{defin}
We say that a model is $\trans$-saturated if it is $\scott{X}\!$-saturated for a
$\trans$-closed Scott set $\scott{X}$.
\end{defin}

Thus; a model of arithmetic is $\satcon$-saturated iff it is recursively saturated and $\SSy(M)$ is $\satcon$-closed.
We formulate the existence of $\trans$-saturated models as a theorem.

\begin{thm}\label{thm:exp.sat.exist}
For any model $M$, countable or not, 
there exists an elementary extension $N \extof M$ which is $\trans$-saturated and such
that $|N|=|M|$.
\end{thm}

If $M$ is an $\scott{X}\!$-saturated model, where $\scott{X}$ is $\trans$-closed, and $T$, $p(\bar x) \in \scott{X}$
(as usual $T$ and $p(\bar x)$ are in a recursive extension of the recursive language of $M$) 
are such that $\satcon(T+p\om/\Th(M))$ we can choose
the type $q(x)$ in part \eqref{item:fyra} of Proposition \ref{prop:satcon.properties} 
to lie in $\scott{X}$ in the following sense.

Let $\bar a \in M$ and $T,p(\bar x) \in \scott{X}$ be in an extension of $\La(\bar a)$, where $\La$
is the language of $M$. 
If 
\[
\satcon(T+\exists x \varphi(x)+p\om/\Th(M,\bar a))
\]
then there is a complete $\La(\bar a)$-type $q(x) \in \scott{X}$ over $\Th(M,\bar a)$ and a new constant $c$ such that 
\[
\satcon(T+\varphi(c)+p\om/\Th(M,\bar a)+q(c)). 
\]
To see this let $S \in
\scott{X}$ be a completion of $T+\varphi(c)$ such that $\satcon(S+p\om/\Th(M,\bar
a))$; there is such a $S$ since $\satcon(T+\varphi(c)+p\om/\Th(M,\bar a))$. Let
\[
q(x)=\set{\psi(x) \in \La(\bar a) | \psi(c) \in S}.
\] 
As $q(c) \subseteq S$ we have $\satcon(S+p\om/\Th(M,\bar a)+q(c))$ and thus
\[
\satcon(T+\varphi(c)+p\om/\Th(M,\bar a)+q(c)).
\]

We can now show that countable transcendent models exist in abundance.

\begin{thm}\label{thm:sat.imp.exp}
Every countable $\trans$-saturated model is transcendent.
\end{thm}
\begin{proof}
Let $M$ be a $\trans$-saturated countable model and 
$\scott{X}$ a $\trans$-closed Scott set such that $M$ is  $\scott{X}\!$-saturated.  Let
$\La^+$, $T$ and $p(\bar x)$ be as in the definition of transcendence.

The expansion $M \redof M^+ \models T+p\om$ is constructed by a Henkin type
construction. Let
$\set{\varphi_k(x)}_{k \in \omega}$
be an enumeration of all formulas in the language $\La^+(M)$ with at most one free variable.

We define a sequence of $\La^+(M)$-sentences $\set{\sigma_{k}}_{k \in \omega}$ such that
\begin{enumerate}

\item\label{a:ett} $\sigma_{k+1} \prf \sigma_k$, and all parameters from $M$ occurring in some $\sigma_i, 
i \leq k$ occur in $\sigma_{k+1}$,

\item\label{a:tva} $\satcon(T+\sigma_k+p\om/\Th(M,\bar b))$, where $\bar b$ are all parameters from $M$ 
occurring in $\sigma_k$,

\item\label{a:tre} $\sigma_k \prf \exists x \varphi_k(x)$ or $\sigma_k \prf
\neg \exists x \varphi_k(x)$,

\item\label{a:fyra} if $\sigma_k \prf \exists x \varphi_k(x)$ then
$\sigma_{k} \prf \varphi_k(m)$ for some $m \in M$, and

\item\label{a:fem} if all elements of $\bar m \in M$ occur in $\sigma_k$ then there exists
$\psi(\bar x) \in p(\bar x)$ such that $\sigma_k \prf \neg \psi(\bar m)$,

\end{enumerate}
hold for all $k \in \omega$.

Given such a sequence of sentences define
$S = T+ \set{\sigma_k}_{k\in \omega}$. The theory $S$ is a complete Henkin theory in the language
$\La^+(M)$. Let $M^+$ be the term model of $S$. The domain of $M^+$ can be
identified with the domain of $M$.
If $\varphi(\bar a)$ is an $\La(M)$-sentence and $M \models \varphi(\bar a)$ then there is $k$ 
such that all elements of $\bar a$ occur in $\sigma_k$ and either $\sigma_k \prf \varphi(\bar a)$ or 
$\sigma_k \prf \lnot\varphi(\bar a)$; let $k$ be greater than the $l$ satisfying that $\varphi(\bar a)$ is
$\varphi_l$, then, by \eqref{a:tre}, either $\sigma_k \prf \exists x \varphi(\bar a)$ or $\sigma_k \prf \lnot \exists
x \varphi(\bar a)$. 
By \eqref{a:tva} we have $\satcon(T+\sigma_k+p\om/\Th(M,\bar a))$ and so  
the theory $\sigma_k+\Th(M,\bar a)$ is consistent. This implies that
$\sigma_k \prf \varphi(\bar a)$ and so $S \prf \varphi(\bar a)$. Thus
\[
\Th(M^+\restrictedto \La,a)_{a \in M^+} = \Th(M,a)_{a \in M},
\]
and so $M^+\restrictedto\La$ is $M$. Clearly $M^+ \models T$, and if $\bar a
\in M^+$ then there is $k$ such that all elements of $\bar a$ occur in $\sigma_k$;
 by \eqref{a:fem} there is some $\psi(\bar x) \in p(\bar x)$ such that $S \prf \neg \psi(\bar a)$, i.e., $M^+
\models \neg \psi(\bar a)$. Therefore $M^+ \models p\om$.

We have to construct such a sequence 
$\{\sigma_k\}$. For the construction to be as uniform as possible define $\sigma_{-1}$ to be  
$\exists x ( x = x)$ and assume $\varphi_0(x)$ to be $x\neq x$.

Suppose $\sigma_{k-1}$ has been constructed. Let $\bar b \in M$ be all parameters 
occurring in $\sigma_{k-1}$ or in $\varphi_k(x)$. 
If
\[
\satcon(T+\sigma_{k-1}+\neg \exists x\varphi_k(x)+p\om/\Th(M,\bar b))
\] 
let $\sigma$ be $\neg \exists x \varphi_k(x)$; this is the case when $k=0$ and so then $\sigma$ is $\lnot \exists x
(x \neq x)$.

Otherwise; if there exists a parameter $d \in \bar b$ such that
\[
\satcon(T+\sigma_{k-1} + \varphi_k(d)+p\om/\Th(M,\bar b))
\]
let $\sigma$ be $\varphi_k(d)$. 
 
For the last case we have
\[
\satcon(T+\sigma_{k-1} + \exists x(\varphi_k(x)\land x \neq \bar b)+p\om/\Th(M,\bar b)),
\]
where $x \neq \bar b$ means $\bigwedge_{b \in \bar b} x \neq b$.
Let $q(x) \in\scott{X}$ be a complete type over  $\Th(M,\bar b)$ 
such that for a new constant symbol $c$ we 
have 
\[
\satcon(T+\sigma_{k-1}+\varphi_k(c)+c \neq \bar b+p\om /\Th(M,\bar b)+q(c));
\]
such a $q(x)$ can be found by the argument preceding this theorem. 
Since $q(x) \in\scott{X}$ and $M$ is $\scott{X}\!$-saturated
there is $m \in M$ realizing $q(x)$. Let $\sigma$ be $\varphi_k(m)$ and expand $\bar b$ to include $m$.

In all three cases we have 
\[
\satcon(T+\sigma_{k-1}+\sigma+p\om/\Th(M,\bar b)).
\]
Let $N$ witness this, i.e., 
\[
N \models T+\sigma_{k-1}+ \sigma+p\om+\Th(M,\bar b)
\] 
is such that $N \restrictedto \La$ is $\omega$-saturated;
we may therefore assume that $M \embin N \restrictedto \La$. For all 
$\bar m \subseteq \bar b$ let $\psi_{\bar m}(\bar x) \in 
p(\bar x)$ be such that $N \models \neg \psi_{\bar m}(\bar m)$. 
Finally let $\sigma_k$ be the conjunction of $\sigma_{k-1}$, 
$\sigma$ and all sentences of the form $\neg \psi_{\bar m}(\bar m)$.

We have to check that $\sigma_k$ satisfies all 
the properties it is supposed to satisfy. \eqref{a:ett} is 
clear since $\sigma_{k-1}$ is one of the conjuncts 
of $\sigma_k$. Property \eqref{a:tva} is also easily 
seen to be true since all the conjuncts of $\sigma_k$ 
is true in the model $N$ above. The other three, \eqref{a:tre},
\eqref{a:fyra} and \eqref{a:fem} are all obviously true.
\end{proof}

In the previous proof full $\trans$-closedness of $\scott{X}$ is not needed; all we used was the
special case when $T_0$ is $\Th(M,\bar a)$ for some $\bar a \in M$, i.e., if $T,p(\bar x) \in \scott{X}$
is in an extension $\La^+$ of $\La(\bar a)$ for some $\bar a \in M$ and $\satcon(T+p\om/\Th(M,\bar a))$ 
then there is  a completion $T_c \in \scott{X}$ of $T$ such that $\satcon(T_c+p\om/\Th(M,\bar a))$. We could even take
this further and get rid of the parameters $\bar a \in M$:

\begin{cor}\label{por:sat.imp.exp}
If $M$ is a countable $\scott{X}\!$-saturated model, where $\scott{X}$ satisfies that for all 
$T,p(\bar x) \in \scott{X}$
in an extension $\La^+$ of the language $\La$ of $M$ satisfying $\satcon(T+p\om/\Th(M))$ 
there is a completion $T_c \in \scott{X}$ of $T$ such that $\satcon(T_c+p\om/\Th(M))$, then 
$M$ is transcendent.
\end{cor}
\begin{proof}
By the argument above we only need to prove that if $T,p(\bar x) \in \scott{X}$
is in an extension $\La^+$ of $\La(\bar a)$ for some $\bar a \in M$ and $\satcon(T+p\om/\Th(M,\bar a))$ 
then there is  a completion $T_c \in \scott{X}$ of $T$ such that 
\[
\satcon(T_c+p\om/\Th(M,\bar a)).
\] 

Let $T$ and $p(\bar x)$ be such. Then $\satcon(T+\Th(M,\bar a)+p\om /\Th(M))$, and so there is a completion $T_c \in \scott{X}$ 
of $T$ such that $\satcon(T_c+\Th(M,\bar a)+p\om/\Th(M))$. But then clearly $\satcon(T_c+p\om/\Th(M,\bar a))$.
\end{proof}

If we fix $p(\bar x)$ in the proof of the theorem we get the following.

\begin{cor}\label{por:sat.imp.exp.fixed.p}
Let $M$ be a countable $\scott{X}\!$-saturated model and $\La^+$ an extension of the language of $M$. Fix an $\La^+$-type 
$p(\bar x)$ and suppose that for all first-order theories $T \in \scott{X}$ such that $\satcon(T+p\om/\Th(M))$ there
is a completion $T_c \in \scott{X}$ of $T$ satisfying $\satcon(T_c+p\om/\Th(M))$. Then for any $\La^+$-theory 
$T \in \scott{X}$ such that $\satcon(T+p\om/\Th(M))$ there is an expansion of $M$ satisfying $T+p\om$. 
\end{cor}

\begin{cor}
Any countable model $M$ has a countable transcendent elementary extension.
\end{cor}
\begin{proof}
Combine Theorems \ref{thm:exp.sat.exist} and \ref{thm:sat.imp.exp}.
\end{proof}

As an easy application we get a joint consistency test for theories of the form $T+p\om$. It should 
be noted that this theorem is provable by a more direct argument as well.

\begin{thm}
Let $T_1$ and 
$p_1(\bar x)$ be a theory and a type in the language $\La_1$ and $T_2$ and $p_2(\bar x)$ in the 
language $\La_2$. Furthermore let $\La_0$ be $\La_1 \cap \La_2$ and $T_0$ a 
complete theory in $\La_0$. Assume also that $\La_0$ is recursive and $\La_1$ 
and $\La_2$ are recursive extensions of $\La_0$. If 
$\satcon(T_1+p_1\om/T_0)$ and $\satcon(T_2+p_2\om/T_0)$ then there exists a countable 
model of $T_0+T_1+T_2+p_1\om+p_2\om$.
\end{thm}
\begin{proof}
Let $\scott{X}$ be a countable $\trans$-closed Scott set 
including $T_0$, $T_1$, $T_2$, $p_1(\bar x)$, and $p_2(\bar x)$.
Assume $M \models T_0$ is countable and $\scott{X}\!$-saturated.
Since 
$\satcon(T_i+p_i\om/T_0)$ for $ i=1,2$ there exists expansions $M_1$ and $M_2$ of $M$ such 
that $M_i \models T_i+p_i\om$. Since $\La_0 = \La_1 \cap \La_2$ we can 
merge $M_1$ and $M_2$ together to one expansion $M^+ \models T_0+T_1+T_2+p_1\om+p_2\om$.
\end{proof}

\subsection{Not all recursively saturated models are transcendent}

Given a language $\La$, which is an extension of the arithmetic language
$\La_A=\set{0,1,+,\cdot,<}$, 
let $K$ be a new unary predicate symbol and define the $\La_{\omega_1 \omega}$-theory 
$\Kom$ to be $S+p\om$ where
\begin{align*}
S &= \set{ K(\num{n}) | n \in \omega }, \\
p(x) &= \set{ x \neq \num{n}  |  n \in \omega } \cup \set{K(x)},
\end{align*}
and $\num{n}$ is the $n$th numeral, i.e., $1+1+\ldots+1$ with $n$ ones.
It is easy to see that $S + p\om$ holds in an expansion of a
model of $\PA$ iff $K$ is interpreted as the set of natural numbers.

Define the $K$-translate, $\Theta^K$, of any second-order
arithmetic formula $\Theta$\footnote{For simplicity we
assume that the only logical symbols of $\Theta$ are $=$, $\vee$, $\neg$,
and~$\exists$.}, so that the following holds:
\begin{align*}
  (t=r)^K \; &\text{is} \; t=r, \\
  (X=Y)^K \; &\text{is} \; \bigl(\forall x (x \in X \ekv x \in Y)\bigr)^K,\\
  (t \in X)^K \; &\text{is} \; K(t) \land (x_X)_t \neq 0, \\
  (\Psi_1 \vee \Psi_2)^K  \; &\text{is} \; \Psi_1^K \vee \Psi_2^K, \\
  (\neg \Psi)^K  \; &\text{is} \; \neg \Psi^K, \\
  (\exists x \Psi)^K \; &\text{is} \; \exists x (K(x) \land \Psi^K), \text{ and} \\
  (\exists X \Psi)^K \; &\text{is} \; \exists x_X \Psi^K
\end{align*}
where $X$ and $Y$ are set variables, $t$ and $r$ are terms, and
$x_X$ is a first-order variable chosen in such way that $x_X$ does
not occur in $\Theta$ and if $X$ and $Y$ are two different second-order
variables then $x_X$ and $x_Y$ are different.

\begin{lem}\label{lem:k.translate}
For any $N \models \PA$, any second-order arithmetic formula 
\[
\Theta(x_0,\ldots,x_k,X_0,\ldots,X_l),
\]
any $n_0$, \dots, $n_k \in \omega$ and any $d_0$, \dots, $d_l \in N$ we have
\begin{align*}
(N,\omega) &\models \Theta^K(n_0,\ldots,n_k,d_0,\ldots,d_l) \quad\text{iff} \\ 
\SSy(N) &\models \Theta(n_0,\ldots,n_k,\codedset_N(d_0),\ldots,\codedset_N(d_l)).
\end{align*}
\end{lem}
\begin{proof}
The proof is by induction on the construction of $\Theta$. First assume $\Theta$ to be atomic. 
There are three cases.
\begin{itemize}
\item $\Theta$ is $t=r$ for some terms $t$ and $r$. 
Clearly $N \models t(\bar n)=r(\bar n)$ iff $t(\bar n)=r(\bar n)$.
\item $\Theta$ is $t \in X$ for some term $t$. $(N,\omega) \models (t \in X)^K(\bar n,d)$ iff 
$(N,\omega) \models K(t(\bar n)) \land (d)_{t(\bar n)} \neq 0$ iff $t(\bar n) \in \codedset_N(d)$.
\item $\Theta$ is $X=Y$. This case reduces to the other cases.
\end{itemize}
If $\Theta$ is not atomic, it is composite; there are three cases here as well.
\begin{itemize}
\item $\Theta$ is $\lnot \Psi$ or $\Psi_1 \vee \Psi_2$. This is obvious from the 
definition (since $K$-translate and $\lnot$/$\vee$ commutes ).
\item $\Theta$ is $\exists x \Psi(x_0,\ldots,x_k,x,X_0,\ldots,X_l)$. $(N,\omega) 
\models \exists x (K(x) \land \Psi^K)(\bar n,\bar d)$ 
iff there is $n \in \omega$ such that $(N,\omega)\models \Psi^K(\bar n,n,\bar d)$ iff there is $n \in \omega$ 
such that $\SSy(N) \models \Psi(\bar n,n,\bar D)$ iff $\SSy(N) \models  
\exists x \Psi(\bar n,x, \bar D)$, where $\bar D$ are the sets coded by the elements
$\bar d$.
\item $\Theta$ is $\exists X \Psi(x_0,\ldots,x_k,X_0,\ldots,X_l,X)$. We have
\[(N,\omega) 
\models \exists x_X \Psi^K(\bar n,\bar d)\] iff there is $d \in N$ 
such that \[(N,\omega) \models \Psi^K(\bar n,\bar d,d)\] 
iff there is $D \in \SSy(N)$ such that \[\SSy(N) \models
\Psi(\bar n,\bar D,D)\] iff 
\[\SSy(N) \models \exists X \Psi(\bar n,\bar D,X).\]
\end{itemize}
By induction the lemma holds for any second-order arithmetic formula $\Theta$.
\end{proof}

\begin{thm}\label{thm:exp.imp.ana}
If $M \models \PA$ is transcendent then $\SSy(M)$ is a $\beta_\omega$-model.
\end{thm}
\begin{proof}
Let $\Psi(\bar{A})$, where $\bar A \in \SSy(M)$, be a second-order sentence true in $\Nat_2$. 
Let $\bar{a} \in M$ code $\bar{A}$, i.e., 
$a_i$ codes $A_i$.  By
taking $N$ to be an $\omega$-saturated model of $\Th(M,\bar a)$ we see that
$(N,\omega) \models \Psi^K (\bar{a})$ since by the lemma this is
equivalent to $\SSy(N) \models \Psi(\bar{A})$ and $\SSy(N) =
\power(\omega)$. Therefore 
\[
\satcon(\Kom + \Psi^K(\bar{a})/\Th(M,\bar a))
\]
and so by the
assumption and Proposition \ref{prop:trans.imp.param.trans}
there is an expansion of $M$ satisfying $\Kom +
\Psi^K(\bar{a})$. There could only be one such expansion
and so we have
\[ (M,\omega) \models \Psi^K(\bar{a}).\]
By using the lemma once again we have that
\[ \SSy(M) \models \Psi(\bar{A})\]
and thus $\SSy(M)$ is a $\beta_\omega$-model.
\end{proof}

\begin{defin}
We say that $M$ is \emph{$\beta_\omega$-saturated} (\emph{$\beta$-saturated}) if it is $\scott{X}\!$-saturated
for some $\beta_\omega$-model ($\beta$-model) $\scott{X}$.
\end{defin}

For historical reasons we note that Jonathan Stavi\footnote{At least Smory{\'n}ski, in \cite{Smorynski:81*4}, claims
that it is due to Jonathan Stavi.} 
proved that a short cofinally expandable model has a standard system which is a $\beta$-model,
see \cite[Fact 3.13]{Smorynski:81*4}.\footnote{In fact, Stavi seems to have proved something weaker, but a
slight modification of his proofs gives the result.}
That is the only other place we found any notion closely related to $\beta$-saturation, though it should be noted that
a short model is never recursively saturated.
However, 
Robert Solovay later proved that no short cofinally expandable models exist, see \cite{Smorynski:82}.

\begin{cor}
If $M \models \PA$ is transcendent then $M$ is $\beta_\omega$-saturated.
\end{cor}
\begin{proof}
Since $M$ is a recursively saturated model of $\PA$ it is $\scott{X}$-saturated, where $\scott{X}=\SSy(M)$, 
and by the theorem above $\SSy(M)$ is a $\beta_\omega$-model.
\end{proof}

The predicate $\satcon$ is certainly not recursive, in fast it is not even arithmetic or analytic. 
The next corollary shows that it is not $\Sigma_k^1$ for any $k \in \omega$, i.e.,
that it is not in the analytical hierarchy of sets.

\begin{cor}\label{cor:satcon.not.analytic}
There is no second-order arithmetic formula $\Theta(X,Y)$ such that 
for all first-order theories $T_0$ and $T$
\[
\satcon(\Kom+T/T_0) \quad \text{iff}\quad \Nat_2 \models \Theta(T,T_0).
\]
\end{cor}
\begin{proof}
Assume, by contradiction, that $\Theta(X,Y)$ is such a formula which is, say, 
$\Sigma^1_k$. 
Let $\scott{X}$ be a countable $\beta_k$-model which is not a $\beta_{k+1}$-model, such a model 
exists by Theorem \ref{thm:exists.beta.not.beta}; and let $M\models\PA$ be countable and $\scott{X}\!$-saturated.
The model $M$ is $\beta_k$-saturated but not $\beta_{k+1}$-saturated, since if $M$ is $\scott{Y}$-saturated then
$\scott{Y}=\scott{X}$, and $\scott{X}$ is not a $\beta_{k+1}$-model.

Assume that $T,T_0 \in \SSy(M)$ and $\satcon(T+\Kom/T_0)$ then
\[
\Nat_2 \models \exists X \bigl(\Theta(X,T_0) \land T \subseteq T_0 \land X \text{ is a complete theory}\bigr).
\]
Since $\SSy(M)$ is a $\beta_k$-model and the sentence is $\Sigma^1_k$ it is also true in $\SSy(M)$ and
so there is a completion $T_c \in \SSy(M)$ of $T$ satisfying $\satcon(T_c+\Kom/T_0)$. 

We don't know if $M$ is $\trans$-saturated, but what we do know is that if $T,T_0 \in \SSy(M)$ are such
that $\satcon(T+\Kom/T_0)$ then there is a completion $T_c \in \SSy(M)$ of $T$ such that 
$\satcon(T_c+\Kom/T_0)$.

Therefore, by Corollary \ref{por:sat.imp.exp.fixed.p}, 
for any $T \in \SSy(M)$ satisfying \[\satcon(T+\Kom/\Th(M))\] there is an expansion of $M$ satisfying
$T+\Kom$. We might call this property $\Kom$-transcendence. In fact, since $M$ is strongly homogeneous, 
for all $\bar a \in M$ the model $(M,\bar a)$ is $\Kom$-transcendent by the same argument that 
proved Proposition \ref{prop:trans.imp.param.trans}.

The proof of Theorem \ref{thm:exp.imp.ana} only uses that $(M,\bar a)$ is $\Kom$-transcendent for all
$\bar a \in M$. Thus; that argument proves that $\SSy(M)$ is a $\beta_\omega$-model which contradicts
the assumption that $\SSy(M)$ is not a $\beta_{k+1}$-model.
\end{proof}

Let $\tp_{\Nat_2}(A)$, where $A \subseteq \omega$, be the type of $A$ in the standard second-order model of arithmetic, i.e.,
\[
\tp_{\Nat_2}(A)=\set{\Theta(X) \text{ second-order arithmetic formula}| \Nat_2 \models \Theta(A)}.
\]

\begin{thm}\label{thm:super.duper.jump}
Let $M \models \PA$ be transcendent; if $A \in \SSy(M)$ then 
\[
\tp_{\Nat_2}(A) \in \SSy(M).
\]
\end{thm}
\begin{proof}
Assume $M \models \PA$ is transcendent and $A \in \SSy(M)$ is coded in $M$ by $a \in M$.
Let $T+p\om$ be
\[
\Kom  +  \set{(c)_n \neq 0 \leftrightarrow \Theta^K(a) | \text{$\Theta(X)$ second-order, } n=\godel{\Theta(X)}}.
\]
If $N$ is an $\omega$-saturated model of $\Th(M)$ and $b \in N$ codes the type $\tp_{\Nat_2}(A)$ then
\[
(N,\omega,b) \models T+p\om
\]
since for all second order $\Theta(X)$ we have 
\[
\Nat_2 \models \Theta(A) \quad\text{iff}\quad (N,\omega,a) \models \Theta^K(a).
\] 
By the transcendence of $M$ there is $d \in M$ such that 
\[ 
(M,\omega,d) \models T+p\om.
\]
Thus, $d$ codes the theory of the second-order model $(\SSy(M),A)$ which is elementary equivalent to $(\Nat_2,A)$ 
since $\SSy(M)$ is a $\beta_\omega$-model.
\end{proof}

Under certain set-theoretic assumptions we have that if
a Scott set is closed under the operator $A \mapsto \tp_{\Nat_2}(A)$ 
then it is a $\beta_\omega$-model:

\begin{thm}
If $\VisL$ or $\PD$ hold then 
any Scott set $\scott{X}$ satisfying the property that if $A \in \scott{X}$ then $\tp_{\Nat_2}(A) \in \scott{X}\!$, 
is a $\beta_\omega$-model.
\end{thm}
\begin{proof}
By Corollary \ref{cor:basis.for.beta} it is enough to prove that $\scott{X}$ 
satisfies true $\Delta^1_\infty$-comprehension. Let $\theta(x,A)$ be any 
second-order arithmetic formula, where $A \in \scott{X}$. Clearly $\scott{X}$ satisfies 
arithmetic comprehension since it is closed under the jump operator. Let $B$ be 
such that
\[ 
\scott{X} \models \forall x (x \in B \ekv \godel{\theta(x,A)} \in \tp_{\Nat_2}(A)), 
\]
clearly
\[ 
\Nat_2 \models \forall x (x \in B \ekv \theta(x,A));
\]
and so $\scott{X}$ satisfies true $\Delta_\infty^1$-comprehension which shows that $\scott{X}$ is indeed
a $\beta_\omega$-model.
\end{proof}

We end this section with an open question and a conjectured partial answer.

\begin{que}
Is, for countable models, transcendence characterised by saturation properties? I.e.,  is there a property
of Scott sets such that any countable model is transcendent iff there is such a Scott set $\scott{X}$ for which 
the model is $\scott{X}\!$-saturated?
\end{que}

Our guess is that this is true for models of $\PA$, but not otherwise.
Clearly, for models of any complete $T \supseteq \PA$ it is true,
since any countable recursively saturated model of $\PA$ is
characterised by its theory and standard system.

\begin{conj}
There is a property of Scott sets such that a countable recursively saturated model $M 
\models \PA$ is transcendent iff $\SSy(M)$ has the property.
\end{conj}

\subsection{Are all transcendent models $\beta_\omega$-saturated?}

Observe that if $T_0$ is complete and has a countable saturated model then there is a $\Sigma^1_2$ formula
$\Theta(X,Y,Z)$ such that if $T$ and $p(\bar x)$ are a theory and a type respectively in an
extension $\La^+$ of the language $\La$ of $T_0$ then 
\[
\Nat_2 \models \Theta(T,p(\bar x),T_0) \quad\text{iff}\quad \satcon(T+p\om/T_0).
\]
To see this let $\Theta(X,Y,Z)$ say 
\begin{multline*}
\exists M \forall q(\bar x,\bar a) \bigl( \text{$M$ is a model of $T_0+T+p\om$} \land {}\\
\text{if $q(\bar x,\bar a)$ is a type over $M$ in $\La$ then $q(\bar x,\bar a)$ is realized in $M$} \bigl).
\end{multline*}
This formula has the desired property since if $\satcon(T+p\om/T_0)$ then there is a countable
$N^+ \models T+T_0+p\om$ such that $N^+\restrictedto \La$ is $\omega$-saturated. To see this let
$L^+ \models T+T_0+p\om$ be such that $L^+ \restrictedto \La$ is $\omega$-saturated and 
$N^+_0 \embin L^+$ a countable model. By a standard argument let $N^+_0 \embin N^+_1 \embin L^+$ be 
countable such that $N^+_1$ realizes all $\La$-types over $N^+_0$ (with parameters from $N^+_0$), such a model
can be constructed since there are only countably many $\La$-types over $N^+_0$. Continue in this fashion
to construct a sequence $\set{N_i^+}_{i \in \omega}$ such that $N_i^+ \embin N_{i+1}^+ \embin L^+$ and such that
$N^+_{i-1}$ realizes all $\La$-types over $N^+_i$. Let $N^+$ be the union of this sequence.

\begin{cor}\label{cor:trans.not.beta}
If $T_0$ has a countable saturated model then there is a transcendent model
of $T_0$ which is $\scott{X}\!$-saturated for a Scott set $\scott{X}$ which is not 
a $\beta_3$-model.
\end{cor}
\begin{proof}
Let $T_1$ be a completion of $T_0$ such that $T_1$ has a countable saturated model and, 
by Theorem \ref{thm:exists.beta.not.beta},
let $\scott{X}$ be a countable $\beta_2$-model which is not a $\beta_3$-model.
Assume $M$ is a countable $\scott{X}\!$-saturated model of $T_1$. 
If $T,p(\bar x) \in \scott{X}$ are such
that $\satcon(T+p\om/T_1)$ and $\Theta(X,Y,Z)$ is as above we have
\[
\Nat_2 \models \exists X \bigl(\Theta(X,p(\bar x),T_1) \land T \subseteq X \land \text{$X$ is complete} \bigr)
\]
and so, since $\scott{X}$ is a $\beta_2$-model, there is a completion $T_c \in \scott{X}$ of $T$ such
that $\satcon(T_c+p\om/T_1)$. By Corollary \ref{por:sat.imp.exp} this is enough for
the model $M$ to be transcendent.
\end{proof}

It is still, for us, not clear whether there is a model which is transcendent but not $\beta_\omega$-saturated. 
If we assume $\VisL$ or $\PD$, to find such a model it is enough to construct a theory $T \in \scott{X}_2$
with countably many complete types of which one is in  $\scott{X}_\omega \setminus \scott{X}_2$, where
$\scott{X}_2$ is the least $\beta_2$-model and $\scott{X}_\omega$ is the least $\beta_\omega$-model. That $\scott{X}_2$
and $\scott{X}_\omega$ exists follows from Corollary \ref{cor:beta.comprehension}. In this case any $\scott{X}_2$-saturated model
of $T$ omits $p(\bar x)$, and is therefore not $\beta_\omega$-saturated, but is transcendent by the discussion above. 

\begin{que}
Is there a transcendent model which is not $\beta_\omega$-saturated? 
\end{que}

%
%
%

\subsection{Alternative consistency notions}

Let us now discuss some possible alternative consistency notions: Are there weaker 
notions of consistency which can replace $\satcon$ in the definition of transcendence? 

We say that a property $\con(T_0,T,p(\bar x))$ is \emph{possible} if any countable model
has a countable elementary extension $M$ satisfying the definition of transcendence with 
\[
\satcon(T+p\om/\Th(M))
\] 
replaced by 
\[
\con(\Th(M),T,p(\bar x)).
\] 

If $\con$ also satisfies that $\con(T_0,T,p(\bar x))$ holds for every $T_0$, $T$, and $p(\bar x)$ such that
there is a model of $T_0+T$ and $\models p\om$; then we say
that $\con$ is \emph{good}, i.e., $\con$ is good if for first-order $T+p\om$ the predicate 
$\con(T_0,T,p(\bar x))$ coincide with 
the ordinary first-order consistency of $T_0+T+p\om$. 

Any model $M$ satisfying the definition of transcendence with 
$\satcon$ replaced by some good $\con$ has to be recursively saturated. 

As we have seen before, defining $\con(T_0,T,p(\bar x))$ to hold iff there is a model of $T_0+T+p\om$ makes
$\con$ not possible. We can do a bit more:

\begin{prop}
Define $\con(T_0,T,p\om)$ to hold iff for all types $q(\bar x)$ over $T_0$ there is a model of
$T_0+q\re+T+p\om$. Then $\con$ is not possible. 
\end{prop}
\begin{proof}
Let $T_0$ be any complete extension of $\PA$.
We will define a theory $T+p\om$ such that no recursively saturated model of $T_0$ has an expansion 
satisfying $T+p\om$, but $\con(T_0,T,p(x))$. If $\con$, defined in this way, was possible it would also be good so finding such a theory
would be enough. The theory $T+p\om$ will be $\Kom$ together with a formalisation of
\begin{equation}\label{eq:is.type}
\text{$\Sigma$ is a truth predicate} \land \text{there is an omitted coded complete type.} 
\end{equation}
Given a type $q(\bar x)$ over $T_0$, let $M$ be the prime model of $T_0+q(\bar c)$; $M$ is not recursively saturated, therefore
\[
(M,\omega,\ThM) \models T+p\om.
\]
Thus; $\con(T_0,T,p(\bar x))$ but no recursively saturated model of $T_0$ has an expansion satisfying $T+p\om$. 

We will give a hint on how to formalise \eqref{eq:is.type}. 
That $\Sigma$ is a truth predicate means that it includes the truth-definition for $\Delta_0$-formulas, which is
definable in $\PA$, that it respects Tarski's definition of truth, and that it is complete for all \emph{standard}
formulas. Let $\Psi(\Sigma,K)$ express this, i.e.,
\[ 
(M,\Sigma,\omega) \models \Psi(\Sigma,K) \quad\text{iff}\quad \Sigma=\ThM,
\]
for every $M \models \PA$. 

By the use of $\Sigma$ and the standard predicate $K$ it is easy to express that an element
$a \in M$ codes a complete type. It should also be clear how to formalise the statement that a type is omitted in
this context. This should, hopefully, convince the reader that \eqref{eq:is.type} is indeed formalisable.
\end{proof}

We end this section with an open question:

\begin{que}
Assume $\con(T_0,T,p(\bar x))$ holds iff there is a model of 
\[
T_0+T+p\om+\set{q\re | q(\bar x) \in S_k(T) \text{ for some $k \in \omega$}}.
\]
Is $\con$ possible?
\end{que}

%
%
%
%
%
%
%

\section{Limit types}

A type $p(\bar x)$ is \emph{isolated in $T$} if there exists $\varphi(\bar x)$ such that
$T+\exists \bar x\varphi(\bar x)$ is consistent and
\[
T \prf \forall \bar x\bigl(\varphi(\bar x) \imp \psi(\bar x)\bigr)
\] 
for all $\psi(\bar x) \in p(\bar x)$.
We say that $p(\bar x)$ is a \emph{limit in  $T$} if it is not isolated in $T$.\footnote{The terms
\emph{isolated} and \emph{limit} comes from the the fact that a complete type over $T$ 
is isolated iff it is a isolated point in the topological space $S(T)$ of all complete types
over $T$.} Observe that if $p(\bar x)$ 
is not consistent with $T$, i.e., it is not a type over $T$, or $T$ itself is inconsistent, then $p(\bar x)$ is 
a limit in $T$.

If $p(\bar x)$ is a limit in $T$ and $\sigma$ is a sentence such that there exists $\varphi(\bar x)$
satisfying 
\[
T+\sigma \prf \forall \bar x\bigl(\varphi(\bar x) \imp \psi(\bar x)\bigr) 
\]
for all $\psi(\bx) \in p(\bx)$,
then
\[
T \prf \forall \bar x\bigl(\sigma \land \varphi(\bar x) \imp \psi(\bar x)\bigr).
\]
for all $\psi(\bx) \in p(\bx)$.
Thus; if the type $p(\bar x)$ is a limit in $T$ then it is also a limit in $T+\sigma$.

Let us formulate this as a proposition:

\begin{prop}
If $p(\bar x)$ is a limit in $T$ and $\sigma$ is a sentence consistent with $T$ then
\begin{enumerate}
\item\label{c:ett} $p(\bar x)$ is a limit in $T+\sigma$, and
\item\label{c:tva} $T+\sigma+p\om$ is consistent iff $T+\sigma$ is consistent.
\end{enumerate}
\end{prop}
\begin{proof}
We have already proven \eqref{c:ett}. 
The omitting types theorem says that if $T$ is consistent then
so is $T+p\om$. The statement \eqref{c:tva} then follows from the omitting types theorem 
and \eqref{c:ett}.
\end{proof}

Furthermore; if $p(\bar x)$ is a limit in $T$ then there is a 
completion $T_c$ of $T$, arithmetic using $T$ and $p(\bar x)$ as oracles, 
such that $T_c+p\om$ is consistent. To see this do the ordinary proof of the omitting types theorem. That
gives you a theory $S$ in the language $\La \cup \set{c_i | i \in \omega}$, where the $c_i$s are new constant symbols and 
$\La$ is the language of $T$ and $p(\bar x)$. The Henkin theory $S$ is such that $T \subseteq S$ and the canonical model
of $S$ omits $p(\bar x)$. Furthermore; $S$ is recursive using $\Th(T)$ and $p(\bar x)$ as oracles, and so is arithmetic using $T$ and $p(\bar x)$ as oracles. The theory $T_c = S \cap \La$ is a completion of $T$ consistent with $p\om$.

Thus; if $\scott{X}$ is any arithmetically closed Scott
set and $T,p(\bx) \in \scott{X}$ are such that $p(\bar x)$ is a limit in $T$ and $T$ is consistent then there
exists a completion $T_c \in \scott{X}$ of $T$ such that $T_c +p\om$ is consistent.


\begin{thm}\label{thm:rec.imp.limittr}
Let $\scott{X}$ be an arithmetically closed Scott set and 
let $M$ be an $\scott{X}\!$-saturated countable model in a recursive language $\La$. 
Assume $T,p(\bar x) \in \scott{X}$ are in a recursive extension 
$\La^+$ of $\La$ and such that $p(\bar x)$ is a limit in $T+\Th(M,\bar m)$ for 
all $\bar m \in M$. If 
$T+\Th(M)$ is consistent then there exists $M^+ \expof M$ such that $M^+ \models
T+p\om$.
\end{thm}
\begin{proof}
The proof follows, more or less, the proof of Theorem \ref{thm:sat.imp.exp}. 
We will construct the sequence $\sigma_k$ as in that proof with the exception that \eqref{a:tva} is 
replaced by the weaker condition

\begin{itemize}
\item[(2')]\label{dd:tva} $T+\sigma_k+\Th(M,\bar b)+p\om$ is consistent, where $\bar b$ 
are all parameters occurring in $\sigma_k$.
\end{itemize}

Clearly this is enough for the term model of $\set{\sigma_k}_{k \in \omega}$
to be (isomorphic to) an expansion of $M$ satisfying $T+p\om$.

Before we construct $\sigma_k$ please observe that $T+\Th(M,\bar m)+\sigma$ is consistent iff
$T+\Th(M,\bar m)+\sigma+p\om$ is consistent, this is true since $p(\bar x)$ is a limit in $T+\Th(M,\bar m)$
and so in $T+\Th(M,\bar m)+\sigma$.

The construction of $\sigma_k$ given $\sigma_{k-1}$ is as follows (as before we let $\sigma_{-1}$ be 
$\exists x (x=x)$ and can forget about the base case in the construction).

Suppose $\sigma_{k-1}$ has been constructed. Let $\bar b \in M$ be all parameters 
occurring in $\sigma_{k-1}$ or in $\varphi_k(x)$. 
If
\[
T+\sigma_{k-1}+\neg \exists x\varphi_k(x)+\Th(M,\bar b)+p\om
\] 
is consistent let $\sigma$ be $\neg \exists x \varphi_k(x)$. 

Otherwise, if
\[
T+\sigma_{k-1} + \varphi_k(d)+\Th(M,\bar b)+p\om
\]
is consistent for some $d \in \bar b$ let $\sigma$ be $\varphi_k(d)$.

In the last case 
\begin{equation}\label{eq:1}
T+\sigma_{k-1} + \varphi_k(c)+ c \neq \bar b + \Th(M,\bar b)+p\om 
\end{equation}
is consistent, where $c$ is a new constant symbol. Since $p(\bar x)$ is a limit in 
$T+\Th(M,\bar b)$ it is also a limit in 
\begin{equation}\label{eq:2}
T+\sigma_{k-1}+\varphi_k(c)+c \neq \bar b + \Th(M,\bar b).
\end{equation}
By the argument preceding the theorem let $S\in\scott{X}$ be any completion of \eqref{eq:2} such that
$S+p\om$ is consistent, and let
\[
q(x) = \set{ \psi(x) \in \La(\bar b) | \psi(c) \in S}.
\]
The type $q(x) \in\scott{X}$, so it is realized by $ d \in M$ where $d \neq \bar b$. 
Let $\sigma$ be $\varphi(d)$ and expand $\bar b$ to include $d$.

In all cases we have that
\[
T+\sigma_{k-1}+\sigma+\Th(M,\bar b)+p\om
\]
is consistent.

Let $N$ witness this, i.e., 
\[
N \models T+\sigma_{k-1}+ \sigma+p\om+\Th(M,\bar b).
\] 
For all $\bar m \subseteq \bar b$ let $\psi_{\bar m}(\bar x) \in p(\bar x)$ be such that $N \models \neg \psi_{\bar m}(\bar m)$. Finally let $\sigma_k$ be the conjunction of $\sigma_{k-1}$, $\sigma$ and all $\neg \psi_{\bar m}(\bar m)$.

We have to check that $\sigma_k$ satisfies all the properties it is supposed to: \eqref{a:ett} is 
clear since $\sigma_{k-1}$ is one of the conjuncts of $\sigma_k$. Property $(2')$ is also easily seen to be true since all the conjuncts of $\sigma_k$ is true in the model $N$ above. The other three, \eqref{a:tre},
\eqref{a:fyra} and \eqref{a:fem} are all obviously true.
\end{proof}

\begin{cor} Let $\scott{X}$ be an arithmetically closed Scott set,
$M$ a countable $\scott{X}\!$-saturated model in a recursive language $\La$,
$\bar a \in M$, and $\La^+$ a recursive extensions of $\La(\bar a)$. For all 
$T,p(\bar x) \in \scott{X}$ in $\La^+$ such that
$p(\bar x)$ is a limit in $T+\Th(M,\bar a,\bar m)$ for all $\bar m \in M$
there exists an expansion of $M$ satisfying $T+p\om$.
\end{cor}
\begin{proof}
The corollary follows from the strong homogeneity of $M$ in the same way as in Proposition
\ref{prop:trans.imp.param.trans}.
\end{proof}

Let $\La^+$ be $\La_A$, the language of arithmetic, with one unary function symbol, $g$, added; and
let $\sigma_{g \in \Aut}$ stand for the sentence expressing that $g$ is an automorphism for the language $\La_A$.
Define the type $p(x)$ to be
\[
\set{ g(x)=x \land x \neq t | \text{$t$ is a closed Skolem term for the language $\La_A$}}.
\]
A model $(M,f) \models \PA+\sigma_{g \in \Aut}+p\om$ iff $M \models \PA$ and $f \in Aut(M)$ is such that
the fixed points of $f$ are exactly the definable points of $M$, i.e., $\mathop\mathrm{fix}(f)=M_0$, where
$M_0$ is the least elementary submodel of $M$. Any such automorphism $f$ is said to be a
\emph{maximal automorphism}. 

\begin{que}
Is there a recursive set of sentences, $S$, in the language $\La_A \cup \set{c,g}$, where $c$ is a new constant symbol and
$g$ is a unary function symbol, satisfying that
for any completion $T_0$ of $\PA$ in the language $\La_A \cup \set{c}$,
$T_0+S$ is consistent and the type $p(\bar x)$ is a limit in $T_0 +S$. 
\end{que}

If the answer is positive we have a converse of \ref{thm:rec.imp.limittr} in the sense that 
a countable model of arithmetic, $M$, is arithmetically saturated iff it satisfies the following property:
For any $T,p(\bar x) \in \SSy(M)$, if $p(\bar x)$ is a limit in $T+\Th(M,m)$ for all $m \in M$, and
$T+\Th(M)$ is consistent then there is $M^+ \expof M$ such that $M^+ \models
T+p\om$. This follows directly from Theorem \ref{thm:k.k.k}.

Observe that there are completions $T_0$, in the language $\La_A \cup \set{c}$, of $\PA$ 
which isolates $p(\bx)$. For example, if $T_0$ is a completion of
$\PA + \set{c \neq t  | \text{$t$ is a Skolem term}}$ then the formula $x=c \land g(c)=c$ isolates 
$p(x)$. 

It should also be mentioned that even if the answer to the question is negative arithmetic saturation
might be strong enough to prove Theorem \ref{thm:rec.imp.limittr} for a slightly larger class of types 
$p(\bx)$,\footnote{For example, types $p(\bx)$ such that $\rk(\Th(M,m),p(\bx))<\omega$ for all $m \in M$,
where $\rk(T,p(\bx))$ is defined in Chapter \ref{ch:proof}.} and
thus, we might still have a converse to the theorem. 

\section{Categorical theories}

We will now study a special sort of theories, which we call categorical,
of the form $T+p\om$, where $T$ is first-order. We
then use the theory of transcendent models to prove an interesting property of them. 
The main, and motivating, example of a categorical theory is $\Kom$.

In model theory we say that a theory is categorical if it only has one model, such a first-order theory does not
exists by the L\"owenheim-Skolem theorems. However, there do exist theories which are $\kappa$-categorical in the 
sense that there is only one model (up to isomorphism) of cardinality $\kappa$ satisfying the theory. The 
categoricity we now define is \emph{over} a model. 

\begin{defin}
A theory $T$, first-order or not, is \emph{categorical over the model} $M$ if there is at most one
expansion of $M$ satisfying $T$. We say that $T$ is categorical over a theory $S$ if $T$ is categorical over 
any model of $S$.  
\end{defin}

For any model $M$, a theory $T$ in the same language as $M$ is categorical over $M$ since
either $M$ satisfies $T$ or not. In either case there is at most one expansion satisfying $T$.
Another, almost as trivial, example is the theory $\Kom$ which is categorical over $\PA$; if $M \models \PA$
there is exactly one expansion of $M$ satisfying $\Kom$, which is $(M,\omega)$. 

Now for the main theorem of this section; 
it says that the expansion to categorical theories $T+p\om$ is well-behaved.

\begin{thm}\label{thm:trans.emb.in.sat}
Assume $M$ is a countable transcendent model in the recursive language $\La$. 
Let $\scott{X}$ be as in the definition of transcendence and 
$T,p(\bar x) \in \scott{X}$ are such that $T+p\om$ is categorical over $M$.
If $\satcon(T+p\om/\Th(M))$ and 
$N \models T+p\om$ witnesses this property, i.e., $N\restrictedto \La$ is $\omega$-saturated, 
then  $M^+ \embin N$, where $M^+$ is the unique expansion of $M$ satisfying $T+p\om$.
\end{thm}
\begin{proof}
Since $N \restrictedto \La$ is $\omega$-saturated we may assume that $M \embin N\restrictedto \La$. 
We prove that the embedding is elementary for $\La^+$. 
Suppose $N \models \varphi(\bar a)$, where $\bar a \in M$ and $\varphi(\bar x)$ is an $\La^+$-formula. Then 
\[
\satcon(T+p\om+\varphi(\bar a)/\Th(M,\bar a))
\]
and therefore, by Theorem \ref{prop:trans.imp.param.trans},  there is an expansion of $(M,\bar a)$ satisfying $T+p\om+ \varphi(\bar a)$.
Since $T+p\om$ is categorical over $M$, and so over $(M,\bar a)$, that expansion has to be $(M^+,\bar a)$. 
Thus $M^+ \models \varphi(\bar a)$.
\end{proof}

\begin{cor}
Let $T+p\om$ be categorical over a complete theory $T_0$ such that 
$\satcon(T+p\om/T_0)$ and let $N_1$ and $N_2$ be witnesses of $\satcon(T+p\om/T_0)$. 
Then $N_1 \equiv N_2$.
\end{cor}
\begin{proof}
Let $M$ be a countable transcendent model of $T_0$ such that $T,p(\bx) \in \scott{X}$, 
where $\scott{X}$ is as in the definition of transcendence. Then $M^+$ is elementary 
embeddable in both $N_1$ and $N_2$, where
$M^+$ is the unique expansion of $M$ satisfying $T+p\om$. Thus; $N_1 \equiv N_2$.
\end{proof}

Indeed, a small modification of the proof shows something stronger, namely that  
the models $N_1$ and $N_2$ are back-and-forth equivalent.

\begin{cor}
Let $M_1 \equiv M_2$ be countable transcendent models, and $\scott{X}_1$ and $\scott{X}_2$ as in the definition 
of transcendence for $M_1$ and $M_2$ respectively. Suppose $T,p(\bar x) \in \scott{X}_1 \cap \scott{X}_2$, that
$T+p\om$ is categorical over both $M_1$ and $M_2$, and that $\satcon(T+p\om/\Th(M_1))$. 
If $M_1^+$ and $M_2^+$ are the unique expansions of $M_1$ and $M_2$, 
respectively, satisfying $T+p\om$, then $M_1^+ \equiv M_2^+$.
\end{cor}
\begin{proof}
Let $N$ witness that $\satcon(T+p\om/\Th(M_1)$, then, by Theorem \ref{thm:trans.emb.in.sat}, 
$M_1^+ \embin N$ and $M_2^+ \embin N$ and so $M_1^+ \equiv M_2^+$.
\end{proof}

We knew before that if $M$ is a transcendent model of $\PA$ then $\SSy(M)$ is a $\beta_\omega$-model, i.e.,
that $\SSy(M) \embin \Nat_2$. 
Theorem \ref{thm:trans.emb.in.sat} says that if $M$ also is countable then $(M,\omega) \embin (N,\omega)$ for 
any $\omega$-saturated model $N \models \Th(M)$, which, by the $K$-translate of second-order formulas, 
is stronger than saying that $\SSy(M) \embin \Nat_2$.

\begin{thm}
Let $M \models \PA$ be transcendent, $\bar a \in M$ and  $T,p(\bar x) \in \SSy(M)$ in an extension $\La^+$ 
of $\La_A(\bar a)$ such that $T+p\om$ is categorical over $(M,\bar a)$. 
If $\satcon(T+p\om / \Th(M,\bar a))$ then there is a completion $T_c \in \SSy(M)$ of $T$ such that
$\satcon(T_c+p\om/\Th(M,\bar a))$. 
\end{thm}
\begin{proof}
Let 
\[
S=\set{ (c)_n \neq0 \ekv \varphi | \text{$\varphi$ an $\La^+$-sentence and $n=\godel{\varphi}$}}
\]
where $c$ is a new constant. Let $N$ witness that $\satcon(T+p\om/\Th(M,\bar a))$ and let $d \in N$
code the theory $\Th(N)$; then $(N,d) \models S$ and therefore 
\[
\satcon(T+S+p\om/\Th(M,\bar a)).
\]
By the transcendence
of $M$ there is an expansion $(M^+,m) \models T+S+p\om$ of $M$, where $m$ is the interpretation of $c$.
Then $\codedset_M(m)=\Th(M^+)$ and by Theorem \ref{thm:trans.emb.in.sat} $M^+ \embin N$ so 
$\Th(M^+)=\Th(N)$. By letting $T_c=\Th(M^+)$ we have a completion $T_c$ in  $\SSy(M)$ of $T$ such that
$\satcon(T_c+p\om/\Th(M,\bar a))$.   
\end{proof}

The last theorem will help us, in the next section, to 
find something like a converse to Theorem \ref{thm:sat.imp.exp}.

\section{Standard recursive saturation}

There is an interesting special case of transcendence for models of $\PA$, 
very much like recursive saturation is a special case of resplendence. If we expand a model
of arithmetic with the standard predicate, $K$, it is easily seen not to be recursively saturated; 
the type
\[ 
\set{ x > n \land K(x) | x \in \omega}
\]
is omitted. However, if we strengthen the consistency assumption of types including the standard predicate to
$\satcon$ we get a new notion of recursive saturation. We call it \emph{recursive standard saturation}.
\begin{defin}
Let $M \models \PA$ and $\bar a \in M$.
A type $q(\bar x,\bar a)$
over the model $(M,\omega)$ is a \emph{standard type over $M$} if there is an
$\omega$-saturated model $N \models \Th(M,\bar a)$ such that $q(\bar x,\bar a)$ is realized in 
$(N,\omega)$.
\end{defin}

In other words, a set of formulas $q(\bar x,\bar a)$ in the language $\La_A(K,\bar a)$, where $\bar a \in M$, 
is a standard type over $M$ iff 
\[
\satcon(\Kom + q\re / \Th(M,\bar a)),
\] 
where $q\re$ is the non first-order sentence expressing that $q(\bar x,\bar a)$ is realized.
Observe that any type over $M$ is a standard type over $M$.
We will often say that $M$ realizes, or omits, a standard type even if we really mean that
$(M,\omega)$ realizes, or omits, the type.

\begin{defin}
A model $M \models \PA$ is \emph{recursively standard saturated}  if it realizes all
recursive standard types over $M$.
\end{defin}

Clearly, recursive standard saturation is stronger than recursive saturation, it says that the the expanded model $(M,\omega)$
is, not recursively saturated, but as much recursively saturated as we could hope for.
Also, any transcendent model is standard recursively saturated.

\begin{lem}
If $M$ is recursively standard saturated then any standard type $q(\bar x,\bar a)$ in $\SSy(M)$ over $M$ is
realized in $(M,\omega)$.
\end{lem}
\begin{proof}
Let $d \in M$ code $q(\bar x,\bar y)$ and define
\begin{multline*}
r(\bar x,\bar a,d) =  \bigl\{(d)_n \neq 0 \imp \varphi(\bar x,\bar a) {}\mathbin{\big|}{}\\ \text{$\varphi(\bar x,\bar y)$ is an $\La_A(K)$-formula and 
$n=\godel{\varphi(\bar x,\bar y)}$}\bigr\}.
\end{multline*}
It is easy to check that $r(\bar x,\bar a,d)$ 
is a recursive standard type over $M$, and that $(M,\omega) \models r\re \ekv q\re$. Therefore,
$q(\bar x,\bar a)$ is realized in $(M,\omega)$.
\end{proof}

\begin{lem}\label{lem:std.comp.type}
Let $q(\bar x,\bar a) \in \SSy(M)$ be a standard type over a recursively standard saturated model $M$.
Then 
there is a complete standard type $r(\bar x,\bar a) \in \SSy(M)$ over $M$ extending $q(\bar x,\bar a)$.
\end{lem}
\begin{proof}
Let 
\begin{multline*}
s(z,\bar x,\bar a) = q(\bar x,\bar a) \cup  \bigl\{(z)_n\neq 0 \ekv \varphi(\bar x,\bar a) \mathbin{\big|} {} \\
\text{$\varphi(\bar x,\bar y)$ is an $\La_A(K)$ formula and $n=\godel{\varphi(\bar x,\bar y)}$}\bigr\}.
\end{multline*}
Let $N \models \Th(M,\bar a)$ be $\omega$-saturated, $\bar b \in N$ such that
$(N,\omega) \models q(\bar b,\bar a)$, and $d \in N$ code the theory $\Th(N,\omega,\bar b,\bar a)$ in $N$. 
It should be clear that $(N,\omega) \models s(d,\bar b,\bar a)$ and therefore that $s(z,\bar x,\bar a)$ is
a standard type over $M$. 
Let $d',\bar b' \in M$ realize $s(z,\bar x,\bar a)$ in $M$; the set coded by $d'$ in $M$ is the theory
$\Th(M,\omega,\bar b',\bar a)$. Let 
\[
r(\bar x,\bar a) = \set{ \varphi(\bar x,\bar a) | \varphi(\bar b',\bar a) \in \Th(M,\omega,\bar b',\bar a)};
\]
then $r(\bar x,\bar a)$ is a complete type, coded in $M$, and extending $p(\bar x,\bar a)$. To see that 
$r(\bar x,\bar a)$ is a standard type we prove that $(M,\omega) \embin (N,\omega)$ which shows that 
$(N,\omega)$ realizes $r(\bar x,\bar a)$ since it is realized in $(M,\omega)$.

Suppose $ \bar a \in M$ and $\varphi(\bar a)$ is an $\La_A(K,\bar a)$-sentence true in $(N,\omega)$. Define the trivial 
recursive standard type 
\[
p(\bar x,\bar a) = \set{\varphi(\bar a),x=x}.
\]
Since $M$ is recursively standard saturated the standard type $p(\bar x,\bar a)$ is realized in $(M,\omega)$ and so 
$(M,\omega)\models \varphi(\bar a)$. Thus $(M,\omega) \embin (N,\omega)$ as we wanted.
\end{proof}

\begin{thm}
Let $M \models \PA$ be countable and recursively saturated; then $M$ is recursively standard saturated iff
for all standard types $q(\bar x,\bar a) \in \SSy(M)$ over $M$ there is a 
complete standard type $r(\bar x,\bar a) \in \SSy(M)$ over $M$ extending $q(\bar x,\bar a)$.
\end{thm}
\begin{proof}
Lemma \ref{lem:std.comp.type} takes care of the left to right direction of the equivalence. 
For the other direction suppose $M \models \PA$ is
countable, recursively saturated and that $\SSy(M)$ satisfies the closedness-condition in the statement of the theorem. Let 
$q(\bar x,\bar a)$ be a recursive standard type over $M$ and let $T=q(\bar c,\bar a)$ where $\bar c$ are some new constants.
We have that $\satcon(\Kom+T/\Th(M,\bar a))$ and we can construct an expansion of $(M,\bar a)$ satisfying $T+\Kom$ in the same
way as Theorem \ref{thm:sat.imp.exp} is proven.

We only have to check that if 
\[
\satcon(\Kom+T+\sigma+\exists x \varphi(x) /\Th(M,\bar b)),
\]
where $\sigma$ and $\varphi(x)$ are
formulas in $\La_A(K,\bar b)$ for some $\bar b \supseteq \bar a$, then there is a complete type $r(x,\bar b)$ 
over $\Th(M,\bar b)$ such that 
\[
\satcon(\Kom+T+\sigma+\varphi(c) /\Th(M,\bar b)+r(c,\bar b)). 
\]
Let 
$s(x,\bar b)$ be the standard type $T+\varphi(x)+\sigma$ over $M$. 
By the assumption on $M$ let $r(x,\bar b)$ be a complete standard
type over $M$ extending $T+\varphi(x)+\sigma$. Thus, 
\[
\satcon(\Kom+T+\sigma+\varphi(c)+r(c,\bar b)/\Th(M,\bar b))
\]
and therefore
\[
\satcon(\Kom+T+\sigma+\varphi(c)/\Th(M,\bar b)+r(c,\bar b))
\]
as we hoped for.

Therefore; there is an expansion of $(M,\bar a)$ satisfying $\Kom+T$. Let $\bar m \in M$ be the interpretation of $\bar c$; then
\[
(M,\omega) \models q(\bar m,\bar a)
\] 
and the arbitrarily chosen recursive standard type $q(\bar x,\bar a)$ is realized in $(M,\omega)$.
\end{proof}

Thus; we have a complete characterisation of recursive standard saturation for countable recursively saturated models
of arithmetic in terms of their standard system. By the proof of Theorem \ref{thm:super.duper.jump} we see that
any recursively standard saturated model has a standard system which is closed under the operation
\[
A \subseteq \omega \mapsto \tp_{\Nat_2}(A).
\]

\begin{que}
Is this condition also sufficient for countable recursively saturated models of arithmetic, i.e., given a countable
recursively saturated model of arithmetic such that for any $A \in \SSy(M)$ the set $\tp_{\Nat_2}(A) \in \SSy(M)$, is 
the model recursively standard saturated?
\end{que}

\begin{que}
Are all countable recursively standard saturated models of arithmetic transcendent?
\end{que}

\section{Subtranscendence}

We now consider a slightly different notion of transcendence. It is much weaker then full transcendence but still quite a lot
stronger than recursive saturation and resplendence. This notion is like transcendence but with the conclusion that there is
an elementary submodel of $M$ with an expansion satisfying $T+p\om$. However, we need to include parameters in the definition which 
makes the definition slightly more complicated.

Before we define this notion let us consider a variant of resplendence, which we call \emph{subresplendence}.
We introduce the following practical short-hand:
\begin{align*}
M\ExpofEmbin N \quad &\text{iff}\quad \exists K ( M \expof K \embin N). 
\end{align*}

\begin{defin}
A model $M$ is \emph{subresplendent} if for every $\bar a \in M$, every
recursive theory $T$ in an extension $\La^+$ of $\La(\bar a)$, where $\La$ is the
language of $M$, such that $T+\Th(M,\bar a)$ is consistent, there is an elementary submodel
$L \embin M$ and an expansion $L^+$ of $L$ satisfying $T$, i.e., there is $L^+ \ExpofEmbin M$ such
that $\bar a \in L^+ \models T$.
\end{defin}

For countable models resplendence and recursive saturation coincide, but not for uncountable
models. Subresplendence, on the other hand, coincide with recursive saturation for arbitrary 
models.

\begin{prop}\label{prop:subresplendency}
A model is recursively saturated iff it is subresplendent.
\end{prop}
\begin{proof}
To prove that a subresplendent model is recursively saturated all we have to do is to observe that it is enough for 
a type to be realized in an elementary submodel of $M$ to be realized in $M$. The other direction of the equivalence 
is proven by a Henkin construction which is a combination of the constructions in the proof of Theorem 
\ref{thm:rec.ekv.resp} and \ref{thm:beta.imp.emb}.
We omit the details here.
\end{proof}

Let us now define the notion of subtranscendence.

\begin{defin}
A recursively saturated model $M$ is \emph{subtranscendent} if there is a 
Scott set $\scott{X}$ for which $M$ is $\scott{X}\!$-saturated, and for all 
$\bar{a} \in M$ and all extensions $\La^+$ of the language 
$\La(\bar a)$ of $(M,\bar{a})$, all $\La^+$-theories $T \in \scott{X}$, and all $\La^+$-types $p(\bx) 
\in \scott{X}$ such that $T+p\om+\Th(M,\bar{a})$ is consistent there is 
$L^+ \ExpofEmbin M$ such that $\bar{a} \in L^+ \models T+p\om$.
\end{defin}

Subtranscendent models exist, as the next theorem shows us.

\begin{thm}\label{thm:beta.imp.emb}  
Every $\beta$-saturated model is subtranscendent.
\end{thm}
\begin{proof}
Let $M$ be a $\beta$-saturated model and $\scott{X}$ a $\beta$-model such that $M$ is
$\scott{X}\!$-saturated.  Let $\La^+$, $\bar{a}$, $T$ and $p(\bar x)$ be as
in the definition of subtranscendence.

We construct the model $L^+ \ExpofEmbin M$ by a similar construction
as the one in the proof of Theorem \ref{thm:sat.imp.exp}.

Let $\set{\exists x \varphi_k(x, \bar{y})}_{k \in
\omega}$ be an enumeration of all $\La^+$-formulas
starting with an existential quantifier such that $\exists x
\varphi_k(x,\bar{y})$ has at most $k$ free variables.  For technical reasons we also assume that
every formula occurs infinitely often in the sequence. We construct a
sequence of sentences $\set{\sigma_k}_{k \in \omega}$, and a sequence
$\set{b_k}_{k \in \omega}$ of elements in $M$, such that if $\bar b_k$ is $b_0, \ldots, b_k$ then

\begin{enumerate}

\item\label{b:ett} $\sigma_{k+1} \prf \sigma_k$, and all elements of $\bar b_k$ and $\bar a$ occur in $\sigma_{k+1}$,

\item\label{b:tva} $T+\sigma_k+p\om+\Th(M,\bar b_k,\bar a)$ is consistent,

\item\label{b:tre} $\sigma_k \prf \exists x \varphi_k(x,\bar b_{k-1})$ or $\sigma_k \prf
\neg \exists x \varphi_k(x,\bar b_{k-1})$,

\item\label{b:fyra} if $\sigma_k \prf \exists x \varphi_k(x,\bar b_{k-1})$ then
$\sigma_k \prf \varphi_k(b_k,\bar b_{k-1})$, and

\item\label{b:fem} if all elements of $\bar m \in M$ occur in $\bar b_k$ or $\bar a$ then there exists
$\psi(\bar x) \in p(\bar x)$ such that $\sigma_k \prf \neg \psi(\bar m)$,

\end{enumerate}
for all $k \in \omega$.

Given such a sequence $\set{\sigma_k}_{k \in \omega}$, the union
$S = \set{\sigma_k| k \in \omega}$ is a complete Henkin theory: To see that $S$ is complete let $\varphi(\bar a,\bar b)$
be a sentence in 
\[
\La^+(\bar a,b_0,b_1,\ldots)
\] 
and let $k \in \omega$ be so large that all elements in $\bar b$
are among $b_0,b_1,\ldots,b_k$ and $\varphi_k(x,\bar b_{k-1})$ is $\varphi(\bar a,\bar b) \land x=x$. By \eqref{b:tre} either
$\sigma_k \prf \varphi(\bar a,\bar b)$ or $\sigma_k \prf \lnot \varphi(\bar a,\bar b)$. The argument to prove that $S$ is Henkin
is similar.

Let
$L^+$ be the term model of $S$. We can identify the domain of
$L^+$ with the set $\set{ b_k | k \in \omega} \subseteq M$, since if
$t(\bar a,\bar b)$ is a term in $\La^+(\bar a,b_0,b_1,\ldots)$ 
then there is a $k \in \omega$ such that
$\varphi_k(x,\bar b_{k-1})$ is $t(\bar a,\bar b)=x$ and so $\sigma_k \prf t(\bar a,\bar b)=b_{k}$.

Furthermore $\bar a \in L^+ \embin_\La M$ since
\[
\Th(L^+\restrictedto\La ,a)_{a \in L^+} \subseteq \ThM.
\]
Clearly $L^+ \models T$, and if $\bar m \in L^+$ then there is $k$ such that
all elements in $\bar m$ appears in $\bar b_k$ and so there is $\psi(\bar x) \in p(\bar x)$
such that $\sigma_k \prf \lnot \psi(\bar m)$ by condition \eqref{b:fem}.

We have to construct such a sequence $\set{\sigma_{k}}_{k\in\omega}$.
For the sake of uniformity let $\sigma_{-1}$ be $\exists x (x=x)$ and assume
$\varphi_0(x)$ to be $x=x$; this makes the base case of the construction trivial.

Inductively define $\sigma_{k}$, given $\sigma_{k-1}$, as follows; if 
\[
T+\lnot \exists x \varphi_k(x,\bar b_{k-1})+\Th(M,\bar a,\bar b_{k-1})+ p\om
\]
is consistent let $\sigma$ be $\lnot \exists x \varphi_k(c,\bar b_{k-1})$ and $b_k=b_{k-1}$. 
Observe that in this case $k\geq 1$ and so $b_{k-1}$ really exists.

Otherwise; if there exists a parameter $d$ in $\bar b_{k-1}$ or in $\bar a$ such that
\[
T+\sigma_{k-1}+\varphi(d)+\Th(M,\bar a,\bar b_{k-1}) +p\om
\]
is consistent let $\sigma$ be $\varphi(d)$ and $b_k=d$.

For the last case let $c$ be a new constant symbol and $S \in \scott{X}$ satisfy
\begin{multline}\label{eq:S}
\text{``} \varphi_k(c,\bar b_{k-1}) \in S \land {}\text{$p(\bar x)$ is
      not isolated in $S$} {} \land {}\\ \text{$S$ is a complete $\La^+(\bar b_{k-1},c)$-theory 
      including} \\\Th(M,\bar a,\bar b_{k-1})
      +T+ \sigma_{k-1}+c \neq \bar a+ c \neq \bar b_{k-1}\text{'',}
\end{multline}
where $c \neq \bar a$ is a shorthand for $c \neq a_0 \land 
\ldots \land c \neq a_{l}$ and similar
for $c \neq \bar b_{k-1}$.
Such a theory $S$ can be found since $\scott{X}$ is a $\beta$-model;
there is an $S \in \power(\omega)$ satisfying \eqref{eq:S} since
\[  
\Th(M,\bar{a},\bar b_{k-1}) + \sigma_{k-1}  +\varphi_k(c,\bar b_{k-1})  + c \neq \bar a + 
c \neq \bar b_{k-1} + p\om
\]
is consistent and hence $p(\bar x)$ is not isolated in the theory of such a model.

Let
\[
q(x) = \set{ \varphi(x) \in \La(\bar a,\bar b_{k-1}) | \varphi(c) \in S(c)};
\]
since $S \in \scott{X}$ we have $q(x) \in \scott{X}$, and since
\[
q(c) + \Th(M,\bar a,\bar b_{k-1})
\]
is consistent, there is an $e \in M$ realizing $q(x)$.

We know that $e$ does not occur in $\bar a$ or $\bar b_{k-1}$ since
$x \neq \bar a \land x \neq \bar b_{k-1} \in q(x)$. If we let $b_{k} = e$ then
\[
T+\sigma_{k-1}+\varphi_k(b_k,\bar b_{k-1})+\Th(M,\bar a,\bar b_k)+p\om
\]
is consistent since, by the omitting types theorem, $S+p\om$ is consistent. 
Let $\sigma$ be $\varphi_k(b_k,\bar b_{k-1})$.

In all three cases we have that
\[
T+\sigma_{k-1}+\sigma+p\om+\Th(M,\bar a,\bar b_k)
\]
is consistent; let $N$ be a model of it.
For all 
$\bar b \subseteq \bar a \cup \bar b_k$ let $\psi_{\bar b}(\bar x) \in 
p(\bar x)$ be such that $N \models \neg \psi_{\bar b}(\bar b)$. 
Finally let $\sigma_k$ be the conjunction of $\sigma_{k-1}$, 
$\sigma$ and all sentences of the form $\neg \psi_{\bar b}(\bar b)$.

We claim that $\sigma_{k}$ satisfies the five properties above.
\end{proof}


\begin{thm}\label{thm:emb.imp.beta}
  If $M \models \PA$ is subtranscendent then $\SSy(M)$ is a $\beta$-model.
\end{thm}
\begin{proof}
Let $\Theta(X,\bar A)$ be an arithmetic formula with set-parameters
$\bar A$ from $\SSy(M)$, such that $\Nat_2 \models
\exists X \Theta(X,\bar{A})$. We want to find $B \in \SSy(M)$
such that $\Nat \models \Theta(B,\bar A)$.

Let $T+p\om$ be $\exists x \Theta^K(x,\bar a) + \Kom$, where $\bar a$ codes $\bar A$.
To see that $\Th(M,\bar a) + T + p\om$ is consistent, take a model $N$ of
$\Th(M,\bar a)$ such that $N$ is $\beta$-saturated, then 
$(N,\omega) \models \Th(M,\bar a) + T+p\om $.

By the assumption that $M$ is subtranscendent there is a model $L^+$ such that
\[
\bar a \in L^+ \embin_\La M \; \text{and}\; L^+ \models T+p\om.
\]

Thus, if $b \in L^+$ is such that $L^+ \models \Theta^K(b,\bar a)$ then
$\Nat \models \Theta(\codedset_{L^+}(b),\bar A)$ and since $L^+$ is elementary
embedded in $M$ 
the set $B=\codedset_{L^+}(b)$ is in $\SSy(M)$.
This completes the proof.
\end{proof}

\begin{cor}
  If $M \models \PA$ is subtranscendent then it is $\beta$-saturated.
\end{cor}
\begin{proof}
  This follows from the fact that a recursively saturated model $M$ of $\PA$ is
  $\SSy(M)$-saturated.
\end{proof}

This characterises the subtranscendent models of $\PA$ as those which are $\beta$-saturated.

\begin{cor}
If $M \models \PA$ is transcendent then it is subtranscendent.
\end{cor}
\begin{proof}
If $M \models \PA$ is transcendent then $M$ is $\beta_\omega$-saturated by
Theorem \ref{thm:exp.imp.ana} and so by Theorem \ref{thm:beta.imp.emb}
$M$ is subtranscendent.
\end{proof}

We also get an interesting corollary about $\beta$-models.

\begin{cor}\label{cor:beta.models.and.completions}
A Scott set $\scott{X}$ is a $\beta$-model iff for every $T,p(\bx) \in \scott{X}$ such that  
$T+p\om$ is consistent there is a completion $T_c \in \scott{X}$ of $T$ such that $T_c+p\om$ 
is consistent.
\end{cor}
\begin{proof}
Assume that $\scott{X}$ is a $\beta$-model, and that $T,p(\bar x) \in \scott{X}$ are such that
$T+p\om$ is consistent. Let $\Theta(X,T,p(\bar x))$ be an arithmetic formula expressing
\[
\text{``$X$ is a complete theory} \land \text{$p(\bar x)$ is a limit in $X$} \land
T \subseteq X \text{''}.
\]
Since there is $X \in \power(\omega)$ satisfying $\Theta(X,T,p(\bar x))$ and $\scott{X}$ is a $\beta$-model,
there is $T_c \in \scott{X}$ such that $\Nat_2 \models \Theta(T_c,T,p(\bar x)$. By the omitting 
types theorem $T_c+p\om$ is consistent.

For the other direction let $\scott{X}$ be such and $M \models \PA$ a countable $\scott{X}$-saturated model.
The proof of Theorem \ref{thm:beta.imp.emb} goes through since it uses only that $\scott{X}$ is closed under
such completions and no other properties of  $\beta$-models, thus $M$ is subtranscendent. Theorem 
\ref{thm:emb.imp.beta} then says that $\SSy(M)=\scott{X}$ is a $\beta$-model.
\end{proof}

%
%
%
%
\def\ro#1{\hbox to 0pt {#1 \hskip 0 pt minus 1 fil}} 
\begin{sidewaysfigure}
\begin{equation*}
\xymatrix{ 
*+[F-:<3pt>]{\mathstrut\text{$\satcon$-saturation}}
\ar@<0ex>@{->}^{\text{countable}}[rrrd]\ar@<0ex>@{->}[dd]& 
  &\hskip 10ex &\\ 
& && *+[F-:<3pt>]{\mathstrut\text{transcendence}}\ar@<0ex>@{->}[d]\ar@<-0.5ex>@{->}_{\PA}[llld] \\ 
*+[F-:<3pt>]{\mathstrut\text{closure under completing standard types}} \ar@<0ex>@{->}[d]\ar@<0ex>@{->}[dr]\ar@<0.5ex>@{->}[rrr]^{\text{countable}}& &&
*+[F-:<3pt>]{\mathstrut\text{standard recursive saturation}}\ar@<0ex>@{->}[dd]\ar@<0.5ex>@{->}[lll]^{\PA}\\
*+[F-:<3pt>]{\mathstrut A \mapsto \tp_{\Nat_2}(A)}\ar@<0ex>@{->}[r]^{\PD \mathbin\lor \VisL} \ar@<0ex>@{->}[d]&    
*+[F-:<3pt>]{\mathstrut\text{$\beta_\omega$-saturation}} \ar@<0ex>@{->}[dl]  \\
*+[F-:<3pt>]{\mathstrut\text{$\beta$-saturation}}\ar@<0ex>@{->}[d]\ar@<0.5ex>@{->}[rrr]& &&
*+[F-:<3pt>]{\mathstrut\text{subtranscendence}}\ar@<0ex>@{->}[d]\ar@<0.5ex>@{->}[lll]^{\PA} \\
*+[F-:<3pt>]{\mathstrut\text{arithmetic saturation}}\ar@<0.5ex>@{->}[rrr]^{\text{countable}}& &&
*+[F-:<3pt>]{\mathstrut\text{limit transcendence}}} 
\end{equation*}
\caption{This is a summary of the results in this chapter. On the left hand side are saturation properties of a model $M$,
on the right hand model theoretic expandability properties of $M$. ``Closed under completing standard types'' means that there is a
Scott set $\scott{X}$ for which $M$ is $\scott{X}\!$-saturated and which is closed under completing standard types.
``$A \mapsto \tp_{\Nat_2}(A)$'' means that there is such an $\scott{X}$ closed under the operation which takes a set $A$
to its complete type in $\Nat_2$. ``Limit transcendence'' refers to the property in Theorem \ref{thm:rec.imp.limittr}. The arrows correspond to implications, and labelled arrows to implications which hold under
some extra assumption (either on $M$ or on the set theoretic universe).}
\end{sidewaysfigure}


\chapter{Proof theory of omitting a type}\label{ch:proof}

The property of 
omitting a type from a proof theoretic point of view is considered in this chapter. 
First, we recapitulate the definition of an isolated type. 

\begin{defin}
  The type $p(\bx)$ is \emph{isolated in $T$} if there exists a formula $\varphi(\bx)$
  such that $T + \exists \bx\varphi(\bx)$ is consistent and $T
  \models \forall \bx \bigl(\varphi(\bx) \imp \psi(\bx)\bigr)$ for
  all $\psi(\bar x) \in p(\bx)$. It is called \emph{strongly isolated in $T$} if, in addition, $T \models
  \exists \bx \varphi(\bx)$. If $p(\bar x)$ is not isolated in $T$ then we say that it is a \emph{limit
  in $T$}.
\end{defin}

For complete theories $T$ a type $p(\bar x)$ is isolated in $T$ iff it is strongly
isolated in $T$, iff $T \models p\re$. Remember that $p\re$ stands for the sentence expressing
that $p(\bx)$ is realized.

Let $S_k(T)$ be the Stone space of complete $k$-types in $T$, i.e., the space
with complete $k$-types in $T$ as points and 
\[
[\varphi(\bx)]=\set{q(\bx) \in S_k(T) | \varphi(\bx) \in q(\bx)}
\]
as basic open sets. These spaces are compact and Hausdorff.
A complete type $p(\bar x) \in S_k(T)$ is isolated in $T$ iff it is isolated as a point
in the topological space $S_k(T)$. 

\begin{thm}[The omitting types theorem]
If $p(\bx)$ is a limit in $T$ then there is a model of $T+p\om$.
\end{thm}

Thus; if $p(\bx)$ is a limit in $T$ then $\Th(T)=\Th(T+p\om)$, where $\Th(T+p\om)$ is the set
of all sentences true in all models of $T+p\om$.
The omitting types theorem reduces the proof theory of $T+p\om$, when $p(\bar x)$ is a limit in $T$,
to the first-order proof theory of $T$.
On the other hand if $p(\bx)$ is isolated in $T$ then $\Th(T)\neq \Th(T+p\om)$.
In this chapter we will investigate the theory $\Th(T+p\om)$ when $p(\bar x)$ is isolated in $T$. 
First we will give a syntactical description of it.

\section{Syntactic characterisation of $\Th(T+p\om)$}

If nothing else is said, the language we are working with will be recursive, and therefore countable. Thus; any type
is countable and we may use a standard enumeration $\set{p_i(\bar x)}_{i\in \omega}$ of the type $p(\bar x)$.

Given a type $p(\bar x)$ we will define an extension of first-order logic by adding a new
inference rule to the ordinary rules. You should think about this rule, schematically, as 
\begin{equation}
  \tag{$p$-rule}\label{p-rule}
  \inference{\lnot \exists \bar x \varphi(\bar x)}{\inferforall{i \in \omega}{\forall \bar x (\varphi(\bar x) \imp
  p_i(\bar x))}}
\end{equation}
i.e., we may deduce $\lnot \exists \bar x \varphi(\bar x)$ if we can deduce $\forall \bar x (\varphi(\bar x) \imp
p_i(\bar x))$ for all $i \in \omega$. However, instead of discussing proof-trees we will define theories $[T]^p_\alpha$,
where $[T]^p_0$ is the first-order closure of $T$, and $[T]^p_{\lambda+1}$ are all sentences provable from $[T]^p_{\lambda}$
with at most one application of the $p$-rule. Let us make this more precise.

When $T$ is a first-order theory, $\Th(T)$ is the logical closure of $T$, i.e., the set
of all sentences provable, in first-order logic, from $T$.
Define
\begin{align*}
(T)^p &= \set{\neg \exists \bx \varphi(\bx) | T \models
\forall \bx \bigl(\varphi(\bx) \imp \psi(\bx)\bigr) \; \text{for all} \; \psi(\bx) \in p(\bx)},\\
[T]^p &= \Th\bigl(T+(T)^p\bigr), 
\end{align*}
and, by recursion,
\begin{align*}
[T]_0^p &= \Th(T), \\
[T]_{\alpha+1}^p &= \bigl[[T]_\alpha^p\bigr]^p, \text{ and} \\
[T]_\lambda^p &= \bigcup_{\alpha < \lambda} [T]^p_\alpha 
\end{align*}
for limit ordinals $\lambda$. Let finally
\[
[T]^p_\infty = \bigcup_{\alpha \in \mathrm{Ord}} [T]^p_\alpha, 
\]
where $\mathrm{Ord}$ is the class of all ordinals.

\begin{lem}\label{lem:tp}
Suppose that 
$\neg \exists \bar x \varphi_0(\bar x)$, \dots, $\neg \exists \bar x 
\varphi_{k-1}(\bar x) \in (T)^p$, then 
\[
\neg \exists \bar x \bigl(\varphi_0(\bar x) 
\lor \ldots \lor \varphi_{k-1}(\bar x)\bigr) \in (T)^p.
\]
\end{lem}
\begin{proof}
If $T \models \forall \bar x \bigl(\varphi_i (\bar x) \imp p_j(\bar x)\bigr)$ for all $i < k, 
j \in \omega$, then 
\[
T \models \forall \bar x \bigl(\varphi_0(\bar x) \lor {}\ldots 
{}\lor \varphi_{k-1}(\bar x) \imp p_j(\bar x)\bigr)
\]
for all $j \in \omega$. By the definition of $(T)^p$ we have that
$\neg \exists \bar x \bigl(\varphi_0(\bar x) 
\lor \ldots \lor \varphi_{k-1}(\bar x)\bigr) \in (T)^p$. This proves the lemma.
\end{proof}

Next some basic properties of the theories $[T]^p$ and $[T]^p_\infty$.

\begin{prop}\label{prop:when.isolated}
Let $p(\bx)$ be a type over $T$, then
\begin{enumerate}
\item $p(\bx)$ is strongly isolated in $T$ iff $[T]^p$ is inconsistent,
\item\label{e:tva} $p(\bx)$ is isolated in $T$ iff $\Th(T) \neq [T]^p$, and
\item $p(\bx)$ is a limit in $[T]_\infty^p$.
\end{enumerate}
\end{prop}
\begin{proof}
\begin{enumerate}
\item If $p(\bx)$ is strongly isolated in $T$ then there is $\varphi(\bx)$ such that $T 
\models \forall \bar x \bigl(\varphi(\bar x) \imp p_i(\bar x)\bigr)$ for all $i \in \omega$ 
and $T \models \exists \bar x \varphi (\bar x)$. Therefore $\neg \exists \bar x 
\varphi(\bar x) \in [T]^p$ and so $[T]^p \models \perp$. 
On the other hand; if $[T]^p \models \perp$ then $T+(T)^p \models \perp$ and, by Lemma \ref{lem:tp}, there 
is $\neg \exists \bar x \varphi(\bar x) \in (T)^p$ such that $T\models \exists \bar 
x \varphi(\bar x)$. Thus; $p(\bx)$ is strongly isolated in $T$ by $\varphi(\bar x)$.

\item If $p(\bx)$ is isolated in $T$ the set $(T)^p$ is non empty. Let $\neg\exists \bar x \varphi (\bar x) 
\in (T)^p$. $T \nmodels \neg\exists \bar x \varphi (\bar x)$ since $T+\exists \bar 
x \varphi (\bar x)$ is consistent, so $[T]^p \neq \Th(T)$. 
On the other hand if $\Th(T) \neq [T]^p$ then there is $\neg \exists \bar x \varphi (\bar 
x) \in (T)^p$ such that $T \nmodels \neg \exists \bar x \varphi (\bar x)$ and so $T+ 
\exists \bar x \varphi (\bar x)$ is consistent, thus $p(\bx)$ is isolated in $T$ 
by $\varphi(\bar x)$.

\item Assume $p(\bx)$ is isolated in $[T]_\infty^p$. Let $\alpha$ be an ordinal such that 
$[T]_\alpha^p=[T]_{\alpha+1}^p$. Then $[T]_\alpha^p = [T]_\infty^p$, and, 
by \eqref{e:tva}, we have 
\[
[T]_\alpha^p = \Th([T]_\alpha^p) \neq \bigl[[T]_\alpha^p\bigr]^p = [T]_{\alpha+1}^p,
\]
which is a contradiction.\qedhere
\end{enumerate}
\end{proof}

The next proposition characterises $[T]_\infty^p$ as the smallest theory including $T$ and 
in which $p(\bar x)$ is a limit.

\begin{prop}\label{prop:least.theory}
Let $p(\bx)$ be a type over $T$. Then
$[T]_\infty^p$ is the least theory closed under first-order provability,
including  $T$, and such that $p(\bx)$ is a limit in $T$.
\end{prop}
\begin{proof}
Clearly $[T]_\infty^p$ is closed under 
first-order provability, including $T$ and, by Proposition \ref{prop:when.isolated}, $p(\bx)$ is a limit in $[T]_\infty^p$. 
Suppose $S$ is another such theory, i.e., $S$ is closed under 
first-order provability, $p(\bx)$ is a limit in $S$, and $T \subseteq S$.

By induction assume that $[T]_\alpha^p \subseteq S$ for all $\alpha < 
\beta$. If $\beta$ is a limit ordinal then clearly $[T]_\beta^p \subseteq S$, if 
not then $\beta=\gamma+1$ and $[T]_\gamma^p \subseteq S$. We need to prove that 
$([T]_\gamma^p)^p \subseteq S$. Assume $\neg \exists \bar x \varphi(\bar x) \in 
([T]_\gamma^p)^p$, i.e., $[T]_\gamma^p \models \forall \bx \bigl(\varphi(\bar x) \imp p_j(\bar 
x)\bigr)$ for all $j \in \omega$. Since $[T]_\gamma^p \subseteq S$ we also have
 $S \models \forall \bar x \bigl(\varphi(\bar x) \imp p_j(\bar 
x)\bigr)$ for all $j \in \omega$. By the assumption on $S$, $p(\bx)$ is a limit in $S$, we have $S 
\models \neg \exists \bar x \varphi(\bar x)$. Therefore $[T]_\beta^p \subseteq S$ and, by induction,
$[T]_\infty^p \subseteq S$. 
\end{proof}

\begin{prop}  
Let $p(\bx)$ be a type over $T$, then
$T+p\om$ is consistent iff $[T]^p_\infty$ is consistent. In 
fact; any model of $T+p\om$ is a model of $[T]^p_\infty$, i.e., $T+p\om \models [T]_\infty^p$.
\end{prop}
\begin{proof}
Suppose $M \models T+p\om$ and 
assume, by induction, that $M \models [T]^p_\alpha$ for all $\alpha < \beta$. 
If $\beta$ is a limit ordinal we get directly that $M \models [T]^p_\beta$; if 
$\beta$ is a successor ordinal, $\beta=\gamma+1$, we need to check that if $M 
\models \forall \bar x \bigl(\varphi(\bar x) \imp p_i(\bar x)\bigr)$ for all $i \in 
\omega$ then $M \models \neg \exists \bar x \varphi(\bar x)$. 
If not then any $\bar a \in M \models \varphi(\bar a)$ would 
realize $p(\bar x)$. Therefore; $M \models [T]_\beta^p$ and by induction $M \models [T]_\infty^p$.

If $[T]_\infty^p$ is consistent then by the omitting types 
theorem there is a model of $[T]_\infty^p+p\om$ since, by 
Proposition \ref{prop:when.isolated}, $p(\bx)$ is 
not isolated in $[T]^p_\infty$. That model is, of course, also a model of 
$T+p\om$.
\end{proof}

The other way around is, in general, not true, i.e., in general not every model 
of $[T]^p_\infty$ is a model of $T+p\om$. A trivial example of this would be to take a type $p(\bar x)$ over $T$,
which is a limit in $T$; then $[T]^p_\infty=T$ and so $[T]^p_\infty+p\re$ is consistent. 

We are now in a position where we can prove the deduction theorem for this kind of proof system.

\begin{prop}\label{prop:deduction}
If $p(\bx)$ is a type over $T$ then
$(\varphi \imp \sigma) \in [T]^p_\infty$ iff $\sigma \in [T+\varphi]^p_\infty$.
\end{prop}
\begin{proof}
If  $(\varphi \imp \sigma) \in [T]_\infty^p$ then clearly $\sigma \in 
[T+\varphi]_\infty^p$ since $[T+\varphi]_\infty^p$ is closed under first-order 
provability, $[T]_\infty^p \subseteq [T+\varphi]_\infty^p$ and $\varphi \in 
[T+\varphi]_\infty^p$.

For the other direction we will, by induction on $\beta$, prove that, 
for all $\sigma$ and $\varphi$, if $\sigma \in [T+\varphi]_\beta^p$ then $\varphi \imp \sigma \in [T]_\beta^p$. 
If $\beta$ is zero or a limit ordinal this is trivial. Therefore; 
assume $\beta=\gamma+1$ and 
that $\sigma \in [T+\varphi]_{\gamma+1}^p$. By Lemma \ref{lem:tp} there is $\neg 
\exists \bar x\psi(\bar x) \in \bigl([T+\varphi]_\gamma^p\bigr)^p$ such that 
\[
[T+\varphi]_\gamma^p \models \neg \exists \bar x \psi(\bar x) \imp \sigma. 
\]
By the induction hypothesis we have
\begin{equation}\label{eq:proof}
[T]_\gamma^p \models \varphi \imp \bigl(\neg \exists \bar x \psi(\bar x) \imp \sigma\bigr).
\end{equation}
We show that 
$[T]_{\gamma+1}^p +\varphi \models \neg \exists \bar x \psi(\bar x)$ since then we 
would have $[T]_{\gamma+1}^p + \varphi \models \sigma$. We know that 
\[
[T+\varphi]_\gamma^p \models \forall \bar x \bigl(\psi(\bar x) \imp p_i(\bar x)\bigr)
\] 
for all $i \in \omega$, and so, by the induction hypothesis again, we have 
\[
[T]_\gamma^p \models \varphi \imp \forall \bar x\bigl(\psi(\bar x) \imp p_i(\bar x)\bigr).
\]
Thus;
\[
[T]_\gamma^p \models \forall \bar x\bigl(\varphi \land \psi(\bar x) \imp p_i(\bar x)\bigr)
\]
for all $i \in \omega$, which implies that 
\[
[T]_{\gamma+1}^p \models \neg \exists \bar x \bigl(\varphi \land \psi(\bar x)\bigr),
\] 
i.e., $[T]_{\gamma+1}^p +\varphi \models \neg \exists \bar x \psi(\bar x)$.
By using \eqref{eq:proof}, we conclude that $[T]_{\gamma+1}^p \models \varphi \imp \sigma$.
\end{proof}

Finally we can prove that $\Th(T+p\om)=[T]^p_\infty$.

\begin{cor}
Let $p(\bx)$ be a type over $T$; then $\Th(T+p\om)=[T]_\infty^p$.
\end{cor}
\begin{proof}
By Proposition \ref{prop:when.isolated} the theory $[T]_\infty^p$ is true in 
any model of $T+p\om$. If $\varphi \notin [T]_\infty^p$ it means, by 
Proposition \ref{prop:deduction}, that $[T+\neg\varphi]_\infty^p$ is consistent. By 
Proposition \ref{prop:when.isolated}, again, there is a model of $T+\neg \varphi + p\om$, i.e., 
$\varphi$ is not true in all models of $T+p\om$, and so $\varphi \notin \Th(T+p\om)$. 
\end{proof}

To follow the practice of logic we write $T+p\om \prf \varphi$ instead of $\varphi \in \Th(T+p\om)$, 
which then, as the corollary shows, is the same as $\varphi \in [T]^p_\infty$.

Observe that up to this point all results are valid for arbitrarily large languages; the assumption that the languages are
recursive was purely for typographic laziness. 
We will now, however, use the recursiveness of the language and investigate the complexity of $\Th(T+p\om)$.

\section{The complexity of $\Th(T+p\om)$}

Let us pin down where the set $\Th(T+p\om)$ is in the analytic hierarchy. It turns out that, for recursive 
$T$ and $p(\bx)$, it is both $\Pi_1^1$ and implicit $\Pi_1^1$. 

\begin{prop}
The set $\Th(T+p\om)$ is uniformly $\Pi_1^{1,T,p(\bx)}$, i.e., there exists a 
$\Pi_1^1$-formula $\Theta(X,Y,z)$ such that for all theories $T$, all types $p(\bx)$ and theories $T$, and all sentences
$\varphi$ we have $\Nat_2 \models  \Theta(T,p(\bx),\varphi)$ iff $\varphi \in \Th(T+p\om)$.
\end{prop}
\begin{proof}
Let $\Theta(T,p(\bar x),\varphi)$ be
\[ 
\forall S \bigl([S]^p \subseteq S \land T \subseteq S \imp \varphi \in S\bigr).
\]
The formula $\Theta(X,Y,z)$ is $\Pi_1^1$ since $\sigma \in [S]^p$ is equivalent to the following 
arithmetic expression:
\begin{multline*}
\leftfix \exists \tau \bigl(\text{$\tau(\bar x)$ is a formula with free variables $\bar x$} \land (S+\lnot \exists \bar x
\tau(\bar x) \prf \sigma) \land {} \\ \forall i [S \prf \forall \bar x (\tau(\bar x) \imp p_i(\bar x))]\bigl)
\end{multline*}
By definition $[T]_\infty^p$ is the smallest theory including $T$, and closed under the operator
$X \mapsto [X]^p$. Thus; $\Nat_2 \models  \Theta(T,p(\bx),\varphi)$ iff $\varphi \in \Th(T+p\om)$ as we
wanted.
\end{proof}

A set $A \subseteq \omega$ is \emph{implicit} $\Pi^{1,\bar B}_1$ if there is a $\Pi_1^1$-formula
$\Theta(X,\bar B)$, with $\bar B$ as parameters, such that 
\[
\Nat_2 \models \exists ! X \Theta(X,\bar B) \land \Theta(A,\bar B).
\]   

\begin{prop} 
The set $\Th(T+p\om)$ is uniformly implicit $\Pi^{1,T,p(\bar x)}_1$, i.e., there 
is a $\Pi_1^1$-formula $\Theta(X,Y,Z)$ such that
\[
\Nat_2 \models \exists ! Z \Theta(T,p(\bx),Z) \land \Theta(T,p(\bx),[T]^p_\infty)
\] 
for all theories $T$ and all types $p(\bx)$.
\end{prop}
\begin{proof}  
Let $\Theta(T,p(\bx),Z)$ express 
\begin{quote}
for all $S$ such that $T \subseteq  S$, $S$ is closed under first-order provability, 
and if $p(\bx)$ is a limit in $S$ then $Z \subseteq T$.
\end{quote}
It should be clear that $\Theta$ can be taken to be $\Pi_1^1$ and that it satisfies the proposition.
\end{proof}

We now investigate a proof theoretic measure of complexity, which indeed is strongly connected to
the analytic hierarchy. 

\begin{defin}
Let 
$\rk(T,p(\bx))$ be the least ordinal $\alpha$ such that $[T]_\alpha^p = [T]_\infty^p$. 
\end{defin}

If $\delta: \power(\omega) \to \power(\omega)$ is an operator on $\power(\omega)$ then we say that
$\delta$ is $\Gamma$, where $\Gamma$ is some complexity class, e.g., $\Pi^{1,\bar A}_k$, 
if the relation $n \in \delta(X)$ is $\Gamma$, i.e., if the set 
$\set{ \langle n,X \rangle | n \in \delta(X) }$ is $\Gamma$. 
  
For any $T$ and $p(\bx)$ let 
\[
\Psi_{T,p} : X \mapsto [X+T]^p
\]
be the operator taking a theory $X$ to the theory $[X+T]^p$. 

The countable ordinal $\omega_1^\mathrm{CK}$ is the least ordinal not order 
isomorphic to any recursive well-ordering.

By adopting Theorem IV.2.15 in \cite{Hinman:78} to our setting we get:

\begin{prop}\label{prop:hinman.ett}
If $\Psi_{T,p}$ is $\Pi^1_1$ then $\rk(T,p)\leq \omega_1^\mathrm{CK}$.
\end{prop}

\begin{prop}\label{prop:rank}
Let $p(\bx)$ be a type over a theory $T$, both in a recursive language $\La$.
If $T$ and $p(\bx)$ both are $\Pi^1_1$  
then $\rk(T,p(\bx))\leq\omega_1^{\mathrm{CK}}$.
\end{prop}
\begin{proof}
Follows directly from the proposition above since if both $T$ and $p(\bar x)$ are $\Pi^1_1$
then so is $\Psi_{T,p}$.
\end{proof}

A set $A$ is \emph{hyperarithmetic} in $\bar B$ if $A$ is $\Delta_1^{1,\bar B}$, i.e., if there are
a $\Sigma_1^1$-formula $\Theta(x,\bar Y)$ and a $\Pi^1_1$-formula $\Psi(x,\bar Y)$ such that
\[
\Nat_2 \models \forall x \bigl(\Theta(x,\bar B) \ekv \Psi(x,\bar B) \ekv x \in A\bigr).
\]

We will now use the theory of subtranscendent models to prove that Proposition \ref{prop:rank} is the best
possible result, at least for hyperarithmetic $T$ and $p(\bx)$. 
First we need a proposition which is sort of a reverse to Proposition \ref{prop:hinman.ett}:

\begin{prop}
If $\rk(T,p(\bx)) < \omega_1^\mathrm{CK}$ then $\Th(T+p\om)$ is $\Delta^{1,T,p(\bx)}_1$.
\end{prop}
\begin{proof}
By Corollary III.3.12 in \cite{Hinman:78}, the least fixed point of $\Psi_{T,p}$ is $\Delta_1^{1,T,p(\bx)}$.
$\Th(T+p\om)$ is the least fixed point of $\Psi_{T,p}$.
\end{proof}

\begin{prop}
There are hyperarithmetic $T$ and $p(\bx)$ (in a recursive language) such that 
$\rk(T,p(x))=\omega_1^\mathrm{CK}$.
\end{prop}
\begin{proof}
If not then, by the proposition above,
$\Th(T+p\om)$ is $\Delta_1^1$ for all hyperarithmetic $T$ and $p(\bx)$.
By using Corollary \ref{cor:beta.models.and.completions} we will prove that the set of hyperarithmetic
sets $\mathrm{HYP}$ is a $\beta$-model, which is a contradiction.

Given $T,p(\bx) \in \mathrm{HYP}$ in a language $\La$, where $p(\bx)$ is a type and $T$ a theory such that 
$T+p\om$ is consistent, 
we will find a completion of $T$ in $\mathrm{HYP}$, which is consistent with $p\om$, 
in much the same way as the omitting types theorem is proven. 

We first observe that the relation $T+p\om \prf \varphi$ is hyperarithmetic, since it is equivalent to
$\varphi \in \Th(T+p\om)$. 

For simplicity we will assume that the type $p(\bx)$ only has one free variable and thus can be written as $p(x)$. 
Let $\exists x \varphi_i(x)$ be an
enumeration of all $\La(C)$-sentences, where $C$ is a set of infinitely 
many new constant symbols. 
 
We will construct sentences $\sigma_i$ for all $i \in \omega$. Start off by letting $\sigma_0$ be $\neg \exists x (x \neq x)$.
If $T+p\om \nprf \sigma_i \imp \exists x \varphi_i(x)$ let $\sigma_{i+1}$ be
\[
\sigma_i \land \neg\exists x\varphi_i(x).
\]
Otherwise let $\psi(x) \in p(x)$ be such that 
\[
T +p\om \nprf \sigma_i \imp \forall x (\varphi_i(x) \imp \psi(x)) 
\]
and let $\sigma_{i+1}$ be
\[
\sigma_i \land \varphi_i(c) \land \neg \psi(c)
\]
for some $c \in C$ not occurring in $\sigma_i$ or $\varphi_i(x)$.

It should be clear that $T+\sigma_i+p\om$ is consistent for each $i \in \omega$ and that the 
theory $T_c = \set{ \varphi \in \La | \exists i \in \omega (\varphi \in T_i) }$ is complete, consistent with $p\om$, and
includes $T$.

The construction of $T_c$ is arithmetic using $\Th(T+p\om)$, $T$, and $p(\bx)$ as oracles. Thus;
$T_c$ is $\Delta^1_1$, or in other words, $T_c \in \mathrm{HYP}$.

By Corollary \ref{cor:beta.models.and.completions} $\mathrm{HYP}$ is a $\beta$-model 
which is a contradiction
(see, for example, \cite[p. 39]{Simpson:99} or \cite[Corollary III.4.8]{Hinman:78}).
\end{proof}

%
%
%
%

\section{Some theories with high rank}

We have seen that there exists hyperarithmetic $T$ and $p(\bx)$ such that $\rk(T,p(\bx))=\omega_1^\mathrm{CK}$,
we will now build some concrete examples of recursive types and theories with high ranks. More exactly, given 
$\alpha < \epsilon_0$ we construct
recursive $T$ and $p(\bx)$ with rank $\alpha$. It should be noted that $\epsilon_0$, the least $\epsilon$ such that
$\omega^\epsilon=\epsilon$, is much less than $\omega_1^\mathrm{CK}$.
 
We start off with an example where $\rk(T,p(x))=2$ which is taken 
from \cite{Casanovas.Farre:96}. 

Let the language $\La$ be $\set{P,Q_i,U_i}_{i \in 
\omega}$, where all symbols are unary predicate symbols. Let the theory $T$ be the set 
of the axioms:
\begin{gather*}
\exists xP(x) \\
\forall y \lnot Q_i(y) \imp \forall x\bigl( P(x) \imp U_i(x)\bigr) \\
\forall x \bigl(Q_i(x) \imp U_j(x)\bigr)\\
\exists^{\geq i} x\bigwedge_{k \leq i} U_k(x)
\end{gather*}
for all $i,j \in \omega$;
and let 
\[ 
p(x)=\set{U_i(x)}_{i \in \omega}.
\]
Clearly $T + p\om$ is inconsistent, since if $M \models T$ then either there is $i$ such that
$M \models \exists x Q_i(x)$ in which case such an $x \in M$ would realize $p(x)$, otherwise 
$M \models \forall x (P(x) \imp U_i(x))$ for all $i$ and any $x \in M$ satisfying $P^M$ will
realize $p(x)$. 
 
Therefore; $[T]^p_\infty$ is the inconsistent theory. In \cite{Casanovas.Farre:96} 
it is shown that $p(x)$ is not strongly isolated in $T$, which means that $[T]_1^p$ is 
consistent; thus $[T]^p_\infty \neq [T]^p_1$. It is easy to see that $[T]_2^p=[T]_\infty^p$, and so we have that
$\rk(T,p)=2$.

\newcommand{\emptyseq}{\epsilon}\newcommand{\conc}{\smallfrown}

We will now generalise this example and find theories $T_\alpha$  such that 
$\rk(T,p(x)) = \alpha$, where $p(x)$ is the same type as above, for each ordinal $\alpha < \omega^2$. However, we first
need some 
general theory about well-founded trees. 

To talk about trees we need to talk about sequences. We denote the empty sequence by  $\emptyseq$ and
$\langle s_0,s_1,\ldots,s_{k-1}\rangle$ the sequence of length $k$ with the $i$th element $s_{i-1}$.
Given a sequence $s=\langle s_0,s_1,\ldots,s_{k-1}\rangle$, the sequence $s \restrictedto l$ is $s$ if $l \geq k$, and
$\langle s_0,s_1,\ldots,s_{l-1}\rangle$ otherwise.

\begin{defin}
A \emph{tree} $\tau$ is a 
subset of $\omega^{<\omega}$ such that for all $s \in \tau$ $s \restrictedto k 
\in \tau$ for every $k\in\omega$. An \emph{infinite branch} in a tree $\tau$ is an $f 
\in\omega^\omega$ such that $\langle f(0),f(1),\ldots,f(k)\rangle \in\tau$ for 
every $k \in \omega$. A tree is \emph{well-founded} if it has no infinite 
branches.
\end{defin}

The reader should be warned that 
the natural numbers $i \in \omega$ will play a dual role, both as a number and as the singleton 
sequence $\langle i \rangle$. Moreover;
if $s,t \in \omega^{<\omega}$ we define $s \conc t \in \omega^{<\omega}$ to be the finite
sequence starting with $s$ and ending with $t$, i.e., the concatenation of $s$ and $t$.

\begin{defin}  
If $\tau$ is a well-founded tree and $s \in  
\omega^{\mathord<\omega}$ we define $\rk_\tau(s)$, the rank of $s$ in $\tau$, 
to be either $-1$ or an ordinal, by recursion:
\[ 
\rk_\tau(s) = \begin{cases}\sup_{i \in \omega} (\rk_\tau(s\conc i) +1) & 
\text{if $s \in \tau$}\\
-1 & \text{otherwise,}
\end{cases}
\]
and let $\rk(\tau)=\rk_\tau(\emptyseq)$.
\end{defin}

The function $\rk_\tau: \omega^{<\omega} \to \omega$ is well-defined
since $\tau$ is well-founded: if not then for some node $s \in\tau$
$\rk_\tau(s)$ would not be defined and so there is a child of $t$ where $\rk_\tau$ is not defined, and so on; this would
define an
infinite path through $\tau$ and violate the well-foundedness of $\tau$.
Observe that $\rk_{\tau}(s)=0$ iff $s$ is a terminating node in the tree $\tau$.
And that $\rk_\tau(s)$ only depends on the subtree of $\tau$ below $s$, i.e., if $\tau\restrictedto s=\set{t | s \conc t \in \tau}$ 
then $\rk_\tau(s)=\rk(\tau\restrictedto s)$.
 
For every countable limit ordinal $\lambda$ fix a strictly increasing sequence of ordinals cofinal in $\lambda$, 
i.e., let $\{\lambda\}: \omega
\to \lambda$ be such that $\sup_{i \in \omega} \{\lambda\}(i) = \lambda$ and $\{\lambda\}(i+1)>\{\lambda\}(i)$ for all $i \in \omega$.

By transfinite recursion define trees $\tau_\alpha$, for countable ordinals $\alpha$: Let
\begin{align*}
\tau_0&=\set{\emptyseq}, \\
\tau_{\alpha+1} &= \set{\emptyseq} \cup \set{ i \conc s | i \in \omega, s
  \in \tau_\alpha }, \text{ and} \\
\tau_\lambda &= \set{\emptyseq} \cup \set{i \conc s | i \in \omega, s \in
  \tau_{\{\lambda\}(i)}},
\end{align*}
for limit ordinals $\lambda$.

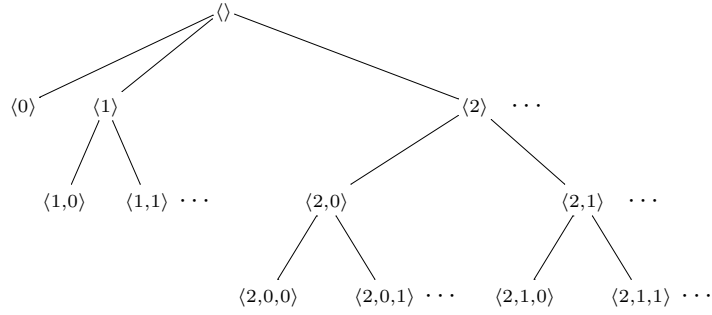
\begin{figure}
\def\binar{\ar@{-}[dl]\ar@{-}[dr]}
\def\seq#1{\scriptstyle\langle#1\rangle}
\[
\xymatrix@1@W=0pt@C=0pt@H=0pt@W=0pt{
   &   &         &   &      &\seq{}\ar@{-}[dlllll]\ar@{-}[dlll]\ar@{-}[drrrrr] \\
\seq{0}&   &\seq{1}\binar&   &      &  &   &         &   &      &\seq{2}\ar@{-}[dlll]\ar@{-}[drr] & \cdots \\
   &\seq{1,0}&         &\seq{1,1}&\cdots&  &   &\seq{2,0}\binar&   &      &  &   &\seq{2,1}\binar&\cdots \\
   &   &         &   &      &  &\seq{2,0,0}&         &\seq{2,0,1}&\cdots&  &\seq{2,1,0}&         
   &\seq{2,1,1}&\cdots 
}
\]
\caption{A picture of $\tau_\omega$ given that $\{\omega\}(k)=k$ for all $k \in \omega$.}
\end{figure}
By transfinite induction on $\alpha$ it is easy to see that $\rk(\tau_\alpha) = \alpha$
for every countable ordinal $\alpha$.

Let us now return to the problem of finding $T$ and $p(\bx)$ with high rank. 
Given a well-founded tree $\tau$ we define a theory $T_\tau$ in the language
\[
\La_\tau=\set{P_s | s \in \tau, s \neq \emptyseq} \cup \set{U_i | i \in \omega}, 
\]
where all predicate symbols are unary, as the set of axioms
\begin{gather*}
\neg \exists x P_{s \conc i}(x) \imp \forall x \bigl(P_s(x)
  \imp U_i(x)\bigr), \text{ for } s\conc i \in \tau, s \neq \emptyseq \text{ and } i \in \omega, \\
\forall x \bigl(P_s(x) \imp U_i(x)\bigr), \text{ for } i \in \omega \text{ and } \rk_\tau(s)=0, \text{ and}\\
\exists^{\geq i} x \bigwedge_{j \leq i} U_j(x), \text{ for } i \in \omega,
\end{gather*}
where $\exists^{\geq i} x \varphi(x)$ is a short-hand for 
\[
\exists x_0, x_1, \ldots, x_{i-1} \bigwedge_{\stackrel{\scriptstyle k,j < i}{k\neq j}} \bigl(x_k \neq x_j \land
\varphi(x_j)\bigr).
\]
As before, let $p(x)$ be the type $\set{U_i(x) | i \in \omega}$.

For convenience we will write $T_\alpha$ for the theory $T_{\tau_\alpha}$; we will also write  $[T]_\beta$ when we really
mean $[T]^p_\beta$ (the type $p(x)$ does not depend on $\beta$ so this should not introduce any ambiguities).

Sometimes the exposition may benefit from thinking of $T_\alpha$ as formulated in the full language 
\[ 
\La_{2^{<\omega}} = \set{P_s | s \in 2^{<\omega}, s \neq \emptyseq} \cup \set{U_i | i \in \omega} 
\]
with the extra axioms
\[
\neg \exists x P_s(x), \text{ for } s \notin \tau_\alpha.
\] 

\begin{lem}\label{lem:strong.theory}
If $\beta \leq \alpha$ are countable ordinals then $[T_\alpha]_{\beta} \prf \lnot \exists x P_s(x)$
for all $s \in \tau_\alpha$ such that $\rk_{\tau_\alpha}(s) < \beta$.
\end{lem}
\begin{proof}
The ordinal $\alpha$ will be fixed, so we simplify things by letting 
$\tau=\tau_\alpha$, $\rk=\rk_\tau$, and $T=T_\alpha$.

We proceed by induction on the ordinal $\beta$. For $\beta=0$ the statement holds trivially. The case $\beta=1$ needs
special treatment. Let $\rk(s)=0$, then $[T]_0 \prf \forall x \bigl(P_s(x) \imp U_i(x)\bigr)$ for all $i \in \omega$. Thus;
$[T]_1 \prf \neg \exists x P_s(x)$ as we hoped for.

For the induction step assume $\beta=\gamma +1 > 1$. Since $[T]_\gamma \prf \neg \exists x P_s(x)$ for
all $\rk(s) < \gamma$, we have $[T]_\gamma \prf \forall x\bigl(P_s(x)
\imp U_i(x)\bigr)$ if $\rk(s \conc i) < \gamma$. Therefore, $[T]_\beta \prf
\neg \exists x P_s (x)$ for all $\rk(s) < \beta$, since if $\rk(s) <
\beta$ then $\rk(s \conc i) < \gamma$.

If $\beta=\lambda$ is a limit ordinal, then $[T]_\lambda=\cup_{\gamma <
  \lambda} [T]_\gamma$ and $[T]_\gamma \prf \neg \exists x P_s(x)$
  for all $\rk(s) < \gamma<\lambda$. Therefore $[T]_\lambda \prf \neg \exists x P_s(x)$ for all $\rk(s) < \lambda$.
\end{proof}

\begin{lem}\label{lem:t.equivalent}
If $\beta \leq \alpha < \omega^2$ then 
\[ 
[T_\alpha]_\beta \equiv T_\alpha + \set{ \neg \exists x P_s(x)| \rk_{\tau_\alpha}(s)<\beta}. 
\]
\end{lem}
\begin{proof}
As before, let $\tau=\tau_\alpha$, $\rk=\rk_\tau$ and $T=T_\alpha$.

We prove the lemma by induction on $\beta$. For $\beta=0$ it is clear.
Assume $\beta=\gamma+1$. We prove that if
\[ 
[T]_\gamma \prf \forall x \bigl(\varphi(x) \imp U_i(x)\bigr)
\]
for every $i \in \omega$, then
\[ 
[T]_\gamma \prf \forall x \bigl(\varphi(x) \imp \bigvee_{s \in t}
P_s(x)\bigr)
\]
for some finite set $t \subseteq \tau$ such that $\rk(s)=\gamma$ for
all $s \in t$. If not, then the theory
\[
[T]_\gamma + \varphi(a) + \set{\neg  P_s(a) | \rk(s)=\gamma }
\]
is consistent by compactness; and since $\exists^{\geq i} x (U_0(x) \land U_1(x)\land \ldots \land U_i(x))$ 
is provable in $[T]_\gamma$ we can conclude that 
\begin{equation}\label{eq:gamma}
[T]_\gamma + \varphi(a) + \set{\neg  P_s(a) | \rk(s)=\gamma } + a\neq
b + \set{U_i(b) | i \in \omega}
\end{equation}
is consistent, again by compactness. Let $M$ be a model of \eqref{eq:gamma}. 

We will
find $i_0 \in \omega$ such that no $P_{s\conc i_0}$ nor any $U_{i_0}$
occur in $\varphi$, and such that if
$\rk(s\conc i_0) < \gamma$ then $\rk(s) \leq \gamma$. Such an $i_0$ can
be found since there are only a finite number of limit ordinals $<
\alpha$, and therefore we can define $i_0$ such that $\{\lambda\}(i_0)
\geq \gamma$ for every limit ordinal $\gamma < \lambda < \alpha$:
Suppose that $\rk(s\conc i_0) < \gamma$, if $\rk(s)$ is not a
limit then $\rk(s) = \rk(s\conc i_0) + 1 \leq \gamma$; otherwise, $\rk(s)$ is
a limit and $\rk(s\conc i_0) = \{\rk(s)\}(i_0) < \gamma$. Thus;
$\rk(s) \leq \gamma$ by the choice of $i_0$.

To recapitulate, $i_0$ is such that $U_{i_0}$ is not in $\varphi$, neither are
any $P_{s \conc i_0}$; furthermore, if $\rk(s \conc i_0) < \gamma$ then
$\rk(s) \leq \gamma$.
Define a model $N$ just like $M$ except that
\begin{align*}
P_{s\conc i_0}^N&=P_{s \conc i_0}^\M \cup \set{b^M}\quad\text{and}\\
U_{i_0}^N&=U_{i_0}^M \setminus \set{a^M}
\end{align*}
for every $s\conc i_0 \in \tau$ such that $\rk(s \conc i_0)\geq\gamma$.
By the choice of $i_0$ $N$ satisfies $\varphi(a)$ and $\lnot U_{i_0}(a)$; also 
if $\rk(s) < \gamma$ then $N$ satisfies $\lnot \exists x P_s(x)$ since, by Lemma \ref{lem:strong.theory}, 
$M$ satisfies it and
$P_s^N=P_s^M$ for all $s$ such that $\rk(s) < \gamma$. Therefore
\[ 
\N \models \set{\lnot \exists x P_s(x) | \rk(s) < \gamma} + \varphi(a) +
\neg U_{i_0}(a). 
\]

We will now prove that $\N \models T$, because then it would follow, by the induction 
hypothesis, that
\[ 
\N \models [T]_\gamma + \varphi(a) + \lnot U_{i_0}(a)
\]
and therefore $[T]_\gamma \nprf \forall x \bigl(\varphi(x) \imp U_{i_0}(x)\bigr)$.

In checking that $N \models T$ we first observe 
that all sentences of the form $\exists^{\geq i}
x \bigwedge_{j \leq i} U_j(x)$ hold in $N$ since 
\[
M \models \exists^{\geq i+1} x \bigwedge_{j \leq i} U_j(x).
\]
Also; the sentences of
the form $\forall x \bigl(P_s(x) \imp U_i(x)\bigl)$, where $s$ is a terminating node, hold in $N$. This follows by the induction
hypothesis trivially if $\gamma > 0$,  since then $M \models \lnot \exists P_s(x)$. Assume $\gamma=0$; then
if $P_s^N(c)$, where $c \in N$, we have $P_s^M(c)$ or $c=b$, in either case $c$ satisfies $U_i^M$. We know that 
$M \models \lnot P_s(a)$ since $\rk(s)=\gamma=0$ and so $c \neq a$ and therefore $c$ satisfies $U_i^N$. 

We need to check that the sentences in $T$ of the form
\[
\lnot \exists x P_{s \conc i}(x) \imp \forall x \bigl(P_s(x) \imp 
U_{i}(x)\bigr)
\]
hold in $N$.
There are two cases which both need careful checking:
\begin{itemize}
\item If $i=i_0$ and $\N \models \lnot \exists x P_{s\conc i_0}(x)$ then 
  $M \models \lnot \exists x P_{s\conc i_0}(x)$ and so
  $M \models \forall x \bigl(P_s(x) \imp U_{i_0}(x)\bigr)$. We also have
  that $\rk(s \conc i_0) < \gamma$ and so, by the choice of $i_0$
  $\rk(s) \leq \gamma$.
  Therefore $\N \models \lnot P_s(a)$.
  If $c \in N$ satisfies $P_s^N$ then either $c$ satisfies $P_s^M$ or $c=b$. 
  In either case it satisfies $U_{i_0}^M$, and, since $c \neq a$, it also
  satisfies $U_{i_0}^N$. Thus $\N
  \models \forall x \bigl(P_s(x) \imp U_i(x)\bigr)$. 
\item If $i \neq i_0$ and $\N \models \lnot \exists x P_{s \conc
  i}(x)$ we have $\M \models \lnot \exists x P_{s \conc i}(x)$ and so
  $\M \models \forall x \bigl(P_s(x) \imp U_i(x)\bigr)$. Therefore $\N \models
  \forall x \bigl(P_s(x) \imp U_i(x)\bigr)$, since the only possible difference in the
  interpretations of $P_s$ is that in $\N$ it may happen that $\N
  \models P_s(b)$ but $U_i^M$ was chosen in such a way that $U_i^M(b)$ and therefore 
  $U_i^N(b)$.
\end{itemize}

All this proves that $[T]_{\gamma +1}=[T]_\beta$ is weaker than 
\[
[T]_\gamma +\set{\lnot \exists x P_s(x) | \rk(s) = \gamma }
\] 
and so, by the induction
hypothesis, is weaker than 
\[
T + \set{\lnot \exists x P_s(x) | \rk(s) \leq \gamma}.
\]
By the preceding lemma it is at least as strong as that theory;
therefore they coincide.

If $\beta=\lambda$ is a limit ordinal, the result follows
immediately by the compactness theorem.
\end{proof}

\begin{thm}
The rank of $T_\alpha$ and $p(x)$ is $\alpha$, i.e., $\rk(T_\alpha,p(x))=\alpha$, for any ordinal $\alpha < \omega^2$.
\end{thm}
\begin{proof}
By the preceding lemma we have
\[ [T_\alpha]_\alpha \equiv \set{ \lnot\exists x P_s(x) |
s \in \tau_\alpha, s \neq \emptyseq} + \Bigl\{ \exists^{\geq i} x
  \bigwedge_{j \leq i} U_j(x) | i \in \omega\Bigr\},
\]
since for $s \in \tau_\alpha$, $\rk_{\tau_\alpha}(s)<\alpha$ iff $s \neq \emptyseq$.  

We prove that $p(x)$ is a limit in $T=[T_\alpha]_\alpha$.
Suppose $T + \exists x \varphi(x)$ is
consistent. Let $\M \models T + \varphi(a)$ and let $i_0$ be such that
$U_{i_0}$ does not occur in $\varphi$. Let $\N$ be like $\M$
except that $U_{i_0}^\N=U_{i_0}^\M \setminus \set{a}$; then $\N \models
T + \varphi(a) + \lnot U_{i_0}(a)$ and so $T \nprf \forall x
\bigl(\varphi(x) \imp U_{i_0}(x)\bigr)$. This proves that $p(x)$ is
not isolated in $T$ and so $\rk(T_\alpha,p(x))\leq \alpha$.

If $\beta < \alpha$ then $p(x)$ is isolated in $[T_\alpha]_\beta$: If
$\rk_{\tau_\alpha}(s)=\beta$ then $p(x)$ is isolated by $P_s(x)$ since, by the Lemma \ref{lem:t.equivalent},
$[T_\alpha]_\beta \models \lnot \exists x P_{s \conc i}(x)$ for all $i \in \omega$. 
Also
$[T_\alpha]_\beta + \exists x P_s(x)$ is consistent: Let 
$M$ be a model with domain $\omega$, and with the predicates 
$U_i^M=\omega$ for all $i \in \omega$, $P_s^M=\omega$, and $P_t^M=\emptyset$ if $t\neq s$; then 
$M \models [T_\alpha]_\beta + \exists x P_s(x)$.

Thus $\rk(T_\alpha,p(x)) > \beta$ for every $\beta < \alpha$ and so $\rk(T_\alpha,p(x))=\alpha$.
\end{proof}

Observe that the proof, and therefore the theorem, does not depend on
the choice of the functions $\{\lambda\} : \omega \to \lambda$ which
was used in the construction of the well-founded trees $\tau_\alpha$.
By choosing these functions explicitly we can extend the result to all
ordinals strictly less than $\epsilon_0$, the least solution of
$\omega^\epsilon=\epsilon$.

\begin{thm}
There are functions $\{\lambda\} : \omega \to \lambda$ such that 
$\rk(T_\alpha,p(x)) = \alpha$ for all countable ordinals $\alpha < \epsilon_0$.
\end{thm}
\begin{proof}
Given any ordinal $\lambda<\epsilon_0$ we can write (uniquely)
\[
\lambda = n_1 \omega^{\lambda_1} + \ldots + n_k \omega^{\lambda_k} 
\]
where $n_i \in \omega$, $n_i > 0$, and $\lambda >\lambda_1 > \ldots > \lambda_k \geq
0$.\footnote{By $\alpha \beta$ we mean $\beta$ taken $\alpha$ times, which, for us, seems more natural than the more
common definition of $\alpha \beta$ as $\alpha$ taken $\beta$ times.} 
This is called the Cantor normal form of $\lambda$ \cite[Theorem 2.26]{Jech:03}. Observe that $\lambda_k > 0$ iff
$\lambda$ is a limit ordinal. 

Define the functions $\{\lambda\} : \omega \to \lambda$, for $\lambda < \epsilon_0$, by transfinite
recursion:
\begin{multline*}
  \leftfix\{n_1 \omega^{\lambda_1} + \ldots + n_k \omega^{\lambda_k}\}(i) =\\
  \begin{cases}
   n_1 \omega^{\lambda_1} + \ldots + (n_k-1) \omega^{\lambda_k} +
  \omega^{\{\lambda_k\}(i)} & \text{if $\lambda_k$ is limit}  \\ n_1
  \omega^{\lambda_1} + \ldots + (n_k-1) \omega^{\lambda_k} + i
  \omega^\xi & \text{if $\lambda_k=\xi+1$}
\end{cases}
\end{multline*}

Fix $\alpha<\epsilon_0$ and $\gamma < \alpha$, and let $\rk=\rk_{\tau_\alpha}$. 
If we can find an $i_0$ such that $\rk(s) \leq \gamma$ when $\rk(s
\conc i_0) < \gamma$ for all $s \in \tau_\alpha$ 
then the proof of Lemma \ref{lem:t.equivalent} will go through as it
stands and $\rk(T_\alpha,p(x))=\alpha$.

Therefore; we need to find an $i_0$ such that for all limit ordinals $\lambda$, if
$\gamma < \lambda < \alpha$ then $\{\lambda\}(i_0) \geq \gamma$.
Fix such a limit ordinal $\lambda$; let 
\[ 
\lambda = n_1 \omega^{\lambda_1} + \ldots  + n_k \omega^{\lambda_k},
\]
and  
\[ 
\gamma = m_1 \omega^{\gamma_1} + \ldots + m_{l}\omega^{\gamma_{l}},
\]
be the Cantor normal forms of $\gamma$ and $\lambda$ respectively 
(observe that $\lambda_k > 0$ but $\gamma_l$ might be $0$). We may assume that $\{\lambda\}(0)<\gamma$.

If $\lambda_k$ is a limit then 
\begin{multline*}
\leftfix\{\lambda\}(0)=n_1 \omega^{\lambda_1} + \ldots + (n_k-1)
\omega^{\lambda_k} + \omega^{\{\lambda_k\}(0)}< m_1 \omega^{\gamma_1}+
\ldots + m_l \omega^{\gamma_{l}} = \gamma \\< n_1 \omega^{\lambda_1} + \ldots + n_k
\omega^{\lambda_k}=\lambda.
\end{multline*}
We can easily see that $l \geq k$, $n_i=m_i$ for $i < k$, $n_k=m_k+1$, and $\lambda_i =
\gamma_i$ for $i \leq k$. Therefore there are at most $l$ different such
$\lambda$, each corresponding to different values of $k$, or in more
loose terms, where to ``chop off'' $\gamma$.

More or less the same argument works if $\gamma_k$ is a successor ordinal, then
\begin{multline*}
\leftfix\{\lambda\}(0)=n_1 \omega^{\lambda_1} + \ldots + (n_k-1)
\omega^{\lambda_k} < m_1 \omega^{\gamma_1}+
\ldots + m_l \omega^{\gamma_{l}} =\gamma \\< n_1 \omega^{\lambda_1} + \ldots + n_k
\omega^{\lambda_k}=\lambda,
\end{multline*}
and $l \geq k$, $n_i=m_i$ for $i < k$, $n_k=m_k+1$, and $\lambda_i =
\gamma_i$ for $i \leq k$.

Thus, in both cases, we can choose $i_0$ such that for each limit $\lambda$ satisfying 
$\gamma < \lambda < \alpha$ we have $\{\lambda\}(i_0) > \gamma$.
\end{proof}


\chapter{Scott's problem}\label{ch:scott}

Scott's problem is to characterise the standard systems for models of first-order
arithmetic. For countable models Dana Scott showed in \cite{Scott:62} that the standard systems are exactly the countable 
Scott sets, i.e., countable boolean algebras of sets of natural numbers closed under relative recursion and K\"onig's lemma. 
It follows, by a union of chains argument, see Theorem \ref{thm:scott.sets.are.ssy}, 
that for models of cardinality $\aleph_1$ the standard systems are 
the Scott sets of cardinality $\aleph_1$.

If the continuum hypothesis holds this settles the problem.
However, if it fails then very little is known about 
standard systems of models of cardinality strictly greater than $\aleph_1$, although it is
easy to see that any standard system of any model is a Scott set. 

Given a Scott set, $\scott{X}\!$, closed under jump, realizing the countable chain condition, and of cardinality 
strictly less than $2^{\aleph_0}$
we will, assuming Martin's axiom, 
construct a model of arithmetic with $\scott{X}$ as its standard system. Any countable Scott set closed under
jump satisfies these conditions. However, we do not know if there
exists any such uncountable Scott sets.

The construction in this chapter is strongly inspired by one of Kanovei's constructions in \cite{Kanovei:96} where he, given a 
countable arithmetically closed set $\scott{X}$,
constructs a model $M$ of true arithmetic with $\SSy(M)=\scott{X}$ and such that a set $A \subseteq \omega$ is representable (without parameters) over $(M,\omega)$ by a $\Sigma_k$-formula iff it is definable (without parameters) 
by a $\Sigma^1_k$ formula over $\scott{X}$.

\section{Definitions}

Let $(P,<)$ be a partial order. Two elements $x,y \in P$ are said to be
\emph{compatible} if there is $z \in P$ such that $z\leq x$ and 
$z\leq y$. 

The partial order $(P,<)$ is said to have the \emph{countable
chain condition} (c.c.c. for short) if for every uncountable set $A \subseteq P$ 
there are $x,y \in A$ such that $x$ and $y$ are compatible. 

Recall that a set $\scott{X} \subseteq \power(\omega)$, where $\power(\omega)$ is the power set of $\omega$,
 is called a \emph{Scott set} if it is a boolean algebra closed under
relative recursion and such that if $\tau \in \scott{X}$ is an infinite binary tree (coded with a suitable G\"odel numbering)
then there is an infinite path in $\scott{X}$ through $\tau$. Any arithmetically closed, i.e., closed under 
relative recursion and the jump operator, set $\scott{X} \subseteq \power(\omega)$ is a Scott set.

A Scott set $\scott{X} \subseteq \power(\omega)$ is said to have the c.c.c. if
the partial order $(\scott{X}_{\mathrm{inf}},\subset)$ has the c.c.c.,
where $\scott{X}_\mathrm{inf}$ is the collection of all infinite sets in $\scott{X}$.

A \emph{filter} $F$ on a partial order $P$ is an up-wards closed subset of $P$ such that if 
$x,y \in F$ then there is
$z \in F$ satisfying $z \leq x$ and $z \leq y$. If $P$ is a boolean algebra a filter $F$ is an \emph{ultrafilter}
if for all $p \in P$ either $p \in F$ or $\lnot p \in F$.

A set $D$ is \emph{dense in $P$} if for every $p \in P$ there is $q\in D$ such that $q\leq p$. If $\mathscr{D}$
is a collection of dense sets we say that a filter $F$ on $P$ is $\mathscr{D}$-\emph{generic} if for 
every $D \in \mathscr{D}$ we have $D \cap F \neq \emptyset$. 

\emph{Martin's axiom}, MA for short, says that for any partial order $P$ with the c.c.c. and any collection
$\mathscr{D}$ of dense sets of cardinality $< 2^{\aleph_0}$ there is a $\mathscr{D}$-generic filter on $P$. Clearly,
$\mathrm{ZFC}+\mathrm{CH} \prf \mathrm{MA}$, but it is also the case that if $\mathrm{ZFC}$ is consistent then so is
$\mathrm{ZFC}+\mathrm{MA}+\lnot\mathrm{CH}$. 
In fact, if $\mathrm{ZFC}$ is consistent and $\kappa \geq \omega_1$ is regular such that
$2^{<\kappa}=\kappa$, then  $\mathrm{ZFC}+\mathrm{MA}+ 2^{\aleph_0}=\kappa$ is consistent, see
\cite[Theorem VIII.6.3]{Kunen:80}.  

Finally, the \emph{standard system} of a model $M$ of $\PA$, written $\SSy(M)$, 
is the collection of all sets of natural numbers coded in $M$, i.e.,
$\SSy(M)=\set{\codedset_M(a) | a \in M}$ where $\codedset_M(a)=\set{n \in \omega | M \models (a)_\num{n} \neq 0}$, $\num{n}$ is
the $n$th numeral and
$(a)_x$ is the $x$th element of the sequence coded by $a$. 

The main theorem we will prove is the following.

\begin{thm}
If $\mathrm{MA}$ holds, $|\scott{Y}| < 2^{\aleph_0}$ is an arithmetically closed Scott 
set with the c.c.c., and $T_0 \in \scott{X}$ is some
completion of $\PA$; then there is a $K \models T_0$ such that $\SSy(K)=\scott{Y}$.
\end{thm}

\section{The construction}

Let $\scott{X}$ be a subset of the power set of the natural numbers. We say that a model $M$ is \emph{coded
in $\scott{X}$} if there is a model $N$ isomorphic to $M$ whose domain  
is in $\scott{X}$ and such that the elementary diagram of $N$, 
$\Th(N,a)_{a \in N}$, also is in $\scott{X}$ (we will as usual identify a formula with
its G\"odel number). Since $M$ and $N$ are isomorphic we will usually assume, for simplicity, 
that $M=N$ when saying that $M$ is coded in $\scott{X}$. 
A set $A$ is recursive in $M$ when $A$ is recursive in the elementary diagram of
$M$, $\Th(M,a)_{a \in M}$. It should be pointed out that this is not the standard definition of a set being recursive in
a model, nor of a model being coded.

Let $\scott{X}$ be any Scott set, 
$T \in \scott{X}$ a completion of $\PA$ and $M \models T$ a model coded in $\scott{X}$. 
Such a model $M$ can be found, for any $T \in \scott{X}$, by doing the ordinary
Henkin construction starting with $T$; the resulting complete Henkin theory will be recursive in $T$ and the 
corresponding model $M$ is therefore coded in $\scott{X}$. 
Let $\prod_\scott{X} M$ be the set of all functions $f: \omega \to M$ which are in $\scott{X}$ 
(we identify a function with its graph, and a pair with its G\"odel number).

For any ultrafilter $U$ in 
$\scott{X}$ define $\prod_\scott{X} M / U$ to be the set of equivalence classes of the equivalence relation
$\equiv_U$ defined on $\prod_\scott{X} M$ by 
\[
f \equiv_U g \quad\text{iff}\quad \set{n | f(n)=g(n)} \in U.
\]
The set $\prod_\scott{X} M / U$ 
can be interpreted as a structure in the language of $\PA$ by interpreting the function symbols pointwise. 
That $\scott{X}$ is a Scott set guarantees that $\prod_\scott{X} M$ is closed under addition and multiplication.
By the canonical map $c \mapsto \bar c$, where
$\bar c$ is the constant function $n \mapsto c$, $M$ is a submodel of $\prod_\scott{X} M /U$. 

Let $K$ be the model $\prod_\scott{X} M /U$, where $U$ is some fixed ultrafilter on $\scott{X}\!$. 

Let $\sigma$ be a sentence in the language $\La_\PA(\prod_\scott{X} M)$. 
The $\La_A(M)$-sentence we get by replacing all occurrences of functions $f$ by the value $f(i)$ will be 
denoted by $\sigma[i]$. By $[\sigma]$ we mean the $\La_A(\prod_\scott{X} M / U)$-sentence we get by replacing all functions
$f$ by the equivalence class $[f]$ of $f$.

The ordinary \L os theorem follows:

\begin{lem}
For any sentence $\sigma$ of $\La_\PA(\prod_\scott{X} M)$ we have
\[ 
K \models [\sigma] \quad\text{iff}\quad \set{i | M \models \sigma[i]} \in U.
\]
\end{lem}
\begin{proof}
This is proved by induction on $\sigma$. 

First, we prove that for a term $t([f_1],\ldots,[f_k])$, if $f \in \prod_\scott{X} M$ is the function
defined by $M \models f(i)= t(f_1(i),\ldots,f_k(i))$ then 
\[
K \models t([f_1],\ldots,[f_k])=[f].
\]

For the base case $t$ is $[f_1]$ and clearly $K \models [f_1]=[f]$ if $f_1=f$. 
If $t$ is 
\[
t_1([f_1],\ldots,[f_k])+t_2([f_1],\ldots,[f_k])
\]
then by the induction hypothesis if 
\[
M \models t_j(f_1(i),\ldots,f_k(i))=g_j(i)
\]
for $j=1,2$ then 
\[
K \models t_j([f_1],\ldots,[f_k])=[g_j]
\]
and so 
\[
K \models t([f_1],\ldots,[f_k])=[g_1]+[g_2]=[g_1+g_2],
\]
and clearly $M \models t(f_1(i),\ldots,f_k(i))=g_1(i)+g_2(i)=(g_1+g_2)(i)$. Similar for the case with
multiplication.

This takes care of the case when $\sigma$ is atomic. The induction step for $\lnot$ and $\vee$ are easy; they
only need that $\scott{X}$ is a boolean algebra.
 
Let us focus on the induction step for the existential quantifier; 
when $\sigma$ is of the form $\exists x \varphi(x,f_1,\ldots,f_k)$.

Then, $K \models [\sigma]$ iff there is $f \in \prod_\scott{X} M$ such that 
\[
K \models \varphi([f],[f_1],\ldots,[f_k]),
\]
iff there is $f \in \prod_\scott{X} M$ such that 
\begin{equation}\label{eq:scott.1}
\set{i | M \models \varphi(f(i),f_1(i),\ldots,f_k(i))} \in U.
\end{equation}
Thus, if $K \models [\sigma]$ there is $f$ satisfying \eqref{eq:scott.1} and so the larger set 
\[
\set{i | M \models \sigma[i]}
\]
is in $U$ since it is recursive using $\ThM$, $f_1$, \dots, $f_k$ as oracles.  

On the other hand if $A=\set{i | M \models \sigma[i]} \in U$ then the function defined by
\[
f(i)=
\begin{cases}
(\mu x) M \models \varphi(x,f_1(i),\ldots,f_k(i)) & \text{ if } i \in A \\
0^M & \text{ otherwise}
\end{cases}
\]
is recursive using $\ThM$, $A$, $f_1$, \dots, $f_k$ as oracles, and so is in $\prod_\scott{X} M$. 
Therefore 
\[
K \models \varphi([f],[f_1],\ldots,[f_k])
\]
and we get that $K \models [\sigma]$.
\end{proof}

Observe that in proving \L os theorem all that is needed is that $\scott{X}$ is a boolean algebra closed under
relative recursion.

\L os theorem gives us that $M \embin K$. To see this
let $K \models \varphi(\bar{c}_1,\ldots,\bar{c}_k)$, i.e., 
\[
\set{i | M \models \varphi(c_1,\ldots,c_k)} \in U.
\]
Since $\emptyset \notin U$ we must have that $M \models \varphi(c_1,\ldots,c_k)$.

Let us now return to the proof of the main theorem.

\begin{proof}[Proof of Theorem 1]
Let $U$ be a non-principal ultrafilter on $\scott{Y}$. That $U$ is non-principal 
implies that all sets in $U$ are infinite and
that if $A \bigtriangleup B =(A \setminus B) \cup (B \setminus A)$ 
is finite then $A \in U$ iff $B \in U$. 
Let $K_U$ be the \emph{proper} elementary extension $\prod_\scott{Y} M_0 / U$ of $M_0$. 

Given $X \in \scott{Y}$ we will find $[f] \in K_U$ such that $[f]$ codes $X$, i.e.,
\[
K_U \models ([f])_\num{n} \neq 0 \quad\text{iff}\quad n \in X,
\]
where $\num{n}$ is the $n$th numeral.
Let $f(k)$ be the code in $M_0$ of the finite set $X \cap \set{0,1,\ldots,k-1}$. 
The definition of $f$ is recursive 
in $M_0$ and $X$, and therefore $f \in \scott{Y}$. For any $n,k \in \omega$ we have
$M_0 \models (f(k))_\num{n}  \neq 0$ iff  $n < k$ and $n \in X$, and so
\begin{align*}
n \in X &\quad\Rightarrow\quad \set{ k | M_0 \models (f(k))_\num{n} \neq 0} =\omega \setminus  \set{0,1,\ldots,n}, \\
n \notin X &\quad\Rightarrow\quad \set{ k | M_0 \models (f(k))_\num{n} \neq 0}  = \emptyset.
\end{align*}

Thus, $n \in X$ iff $\set{ k | M_0 \models (f(k))_\num{n} \neq 0} \in U$, i.e., iff $K_U \models ([f])_\num{n} \neq 0$. 

So far we have neither used the fact that $\scott{Y}$ is arithmetically closed nor that $\scott{Y}$ has the c.c.c. To prove that
all sets coded are in $K_U$ we need to find a generic ultrafilter $U$, for that we need these extra conditions on $\scott{Y}$.

Given any filter $F$ on $\scott{Y}_\mathrm{inf}$, let $F'$ be $F$ with all cofinite sets added. The collection $F'$ of sets 
has the finite intersection property in $\scott{Y}$ 
since if $X$ is infinite and $Y$ is cofinite then $X \cap Y$ is infinite. 
In $\scott{Y}$ $F'$ can be extended to an ultrafilter
$U$ and since $U$ includes all cofinite sets it is non-principal. Thus, every filter on 
$\scott{Y}_\mathrm{inf}$ can be extended to a non-principal ultrafilter on $\scott{Y}\!$. Therefore 
MA gives us, for any collection  $\mathscr{D}$ of
dense subsets of $\scott{Y}_\mathrm{inf}$ of cardinality strictly less than $2^{\aleph_0}$, a non-principal 
$\mathscr{D}$-generic ultrafilter on $\scott{Y}$.  

We say that a $\Delta_0(\prod_\scott{Y}M_0)$ sentence $\sigma$ is \emph{true on} $X \subseteq \omega$ if
$M_0 \models \sigma[k]$ for almost all $k \in X$ (i.e., for all but a finite number). Observe that, since we are assuming that
$U$ is a non-principal ultrafilter,
$K_U \models [\sigma]$ iff there is $X \in U$ such that 
$\sigma$ is true on $X$. The set $X$ \emph{decides} $\sigma$ if either $\sigma$ or $\lnot \sigma$
is true on $X$. 

If $f \in \prod_\scott{X} M_0$ let
\[
D_f=\set{V | \forall n \bigl(V \text{ decides } (f)_\num{n}=0\bigr)}
\]
and 
$\mathscr{D}=\set{D_f | f \in \prod_\scott{X} M_0}$.
We prove that $\mathscr{D}$ is a collection of dense sets, i.e.,  
given $f \in \prod_\scott{Y} M_0$ and $X \in \scott{Y}_\mathrm{inf}$ there is
$V \in \scott{Y}$, $V \subseteq X$, deciding all sentences of the form $(f)_\num{n}=0$.

Let such $f$ and $X$ be given, 
 $V_0=X$, and let $V_{n+1} \subseteq V_n$ be either of
\[
\set{k \in V_n | M_0 \models (f(k))_\num{n} = 0} \quad\text{or}\quad 
\set{k \in V_n | M_0 \models (f(k))_\num{n} \neq 0},
\]
whichever is infinite.

Define $V=\set{i_n}_{n \in \omega}$ where $i_n \in V_n$ and $i_{n+1}$ is chosen to be the least element of $V_{n+1}$
strictly larger than $i_n$;  then 
\[
V \setminus V_n \subseteq \set{i_0,i_1,\ldots i_{n-1}}. 
\]
Hence 
$V$ decides all sentences $(f)_\num{n}=0$.   
Since $\scott{Y}$ is arithmetically closed it is easy to see that
$V \in \scott{Y}$ and so $V \subseteq X$ and $V \in D_f$. 

Since $\mathscr{D}$ is a collection of dense sets and $|\mathscr{D}|< 2^{\aleph_0}$, by 
MA there is a $\mathscr{D}$-generic non-principal ultrafilter on $\scott{Y}$. 
Let $G$ be such a non-principal $\mathscr{D}$-generic ultrafilter on $\scott{Y}$ and let
$K$ be $K_G$.

To prove that any $X \subseteq \omega$ coded by some $[f] \in K$ is in $\scott{Y}$ we 
use the fact that if $f \in \prod_\scott{Y} M_0$ then there
is a $V \in G$ deciding all formulas $(f)_\num{n} = 0$. Since
\[
K \models ([f])_\num{n} = 0 \quad\text{iff}\quad \set{ k | (f(k))_\num{n} = 0 } \in G
\]
we get
$K \models ([f])_\num{n} = 0$ iff $(f)_\num{n}=0$ is true on $V$.
Thus the set $X$ is arithmetic in $f$ and $V$ and therefore $X \in \scott{Y}$.
This ends the proof of the equality $\SSy(K)=\scott{Y}$.
\end{proof}

\begin{thm}
Assume $\scott{X}$ is a Scott set, $T \in \scott{X}$ is any completion of $\PA$ and
$M \models T$ is coded in $\scott{X}$. Let also $U$ be a non-principal ultrafilter on
$\scott{X}\!$; then the elementary extension $K=\prod_\scott{X} M/U$ of $M$ is recursively
saturated.
\end{thm}
\begin{proof}
Let $p(x,[g])=\set{\varphi_i(x,[g])}_{i \in \omega}$ be a recursive type over $K$ with, 
for simplicity, $[g]$ as the only parameter. We may assume that
\[
\models \forall x \bigl(\varphi_{i+1}(x,[g]) \rightarrow \varphi_i(x,[g])\bigr).
\]
Then
\[
K\models \exists x \varphi_i(x,[g])
\]
for all $i \in \omega$; let 
\[
V_i = \set{k | M \models \exists x \varphi_i(x,g(k)) } \in U.
\]
Define $f \in \prod_\scott{X} M$ as follows: Given $k \in \omega$ find the largest $i\leq k$, if it exists, 
such that $k \in V_i$, i.e., $M \models \exists x \varphi_i(x,g(k))$, 
and let $f(k)$ be such that $M \models \varphi_i(f(k),g(k))$. If no such $i$ exists let
$f(k)$ be $0^M$. The function $f$ is recursive in $M$ and satisfies the property 
\[ 
\text{if } k \in V_i \text{ and } i \leq k \text{ then } M \models \varphi_i(f(k),g(k)).
\]
Fix $i \in \omega$ and consider the set
\[ 
V_i'=\set{ k | M \models \varphi_i(f(k),g(k))}.
\]
If $k \geq i$ then $k \in V_i$ iff $k \in V_i'$ so $V_i \bigtriangleup V_i'$ is finite and therefore 
$V_i' \in U$. Thus;
\[
K \models \varphi_i([f],[g])
\]
for all $i \in \omega$. Any recursive type over $K$ is realized in $K$ and so $K$ is recursively saturated.
\end{proof}

Thus, given any arithmetically closed Scott set $\scott{X}$ with the c.c.c. and
of cardinality strictly less than the continuum, 
and a completion $T \in \scott{X}$ of $\PA$, we can
construct an $\scott{X}\!$-saturated model $M \models T$.

However, if interesting such Scott sets exist is unclear to us.

\begin{que}
Does there exist uncountable Scott sets closed under jump with the c.c.c.?
\end{que}

Clearly any countable Scott has the c.c.c. By a rather simple argument, see for example
\cite[Chapter II: Theorem 1.2]{Kunen:80}, it can be
seen that the full power set $\power(\omega)$ 
does not have the c.c.c.
Nothing else concerning Scott sets and the c.c.c. is known to us.

\backmatter

\bibliography{thesis}
\bibliographystyle{halpha}
\end{document}

%% file: head.tex
\newcommand{\rep}{\mathrm{rep}}

\newcommand{\WKL}{\mathsf{WKL}}
\newcommand{\RCA}{\mathsf{RCA}}
\newcommand{\ACA}{\mathsf{ACA}}
\newcommand{\M}{M}
\newcommand{\N}{N}

\newcommand{\LOW}{\mathrm{LOW}}

\newcommand{\scott}[1]{\mathscr{#1}}
\newcommand{\ThM}{\mathrm{Th}(M,a)_{a \in M}}

\newcommand{\PA}{\mathrm{PA}}
\newcommand{\SSy}{\mathop\mathrm{SSy}}
\newcommand{\Aut}{\mathop\mathrm{Aut}}
\newcommand{\tp}{\mathop\mathrm{tp}\nolimits}
\newcommand{\dcl}{\mathop\mathrm{dcl}}

\newcommand{\Th}{\mathop\mathrm{Th}}
\newcommand{\om}{\mathord\uparrow}
\newcommand{\re}{\mathord\downarrow}
\newcommand{\La}{\mathscr{L}}

\newcommand{\VisL}{V\mathord=L}
\newcommand{\PD}{\text{PD}}
\newcommand{\Kom}{T_{K=\omega}}

\newcommand{\nmodels}{\nvDash}
\newcommand{\prf}{\vdash}
\newcommand{\nprf}{\nvdash}
\newcommand{\imp}{\rightarrow}

\newcommand{\restrictedto}{\mathord\upharpoonright}
\newcommand{\trans}{\satcon}

\DeclareFontFamily{U}{feifs}{}
\DeclareFontShape{U}{feifs}{m}{n}{<5> <6> <7> <8> <9> <10> 
   <11> <12> feifs10}{}
\DeclareSymbolFont{fs}{U}{feifs}{m}{n}
\DeclareMathSymbol{\expof}{\mathrel}{fs}{3}
\DeclareMathSymbol{\redof}{\mathrel}{fs}{4}
\DeclareMathSymbol{\ExpofEmbin}{\mathrel}{fs}{11}
\DeclareMathSymbol{\RedofEmbin}{\mathrel}{fs}{12}
\DeclareMathSymbol{\ExpofExtin}{\mathrel}{fs}{21}
\DeclareMathSymbol{\RedofExtin}{\mathrel}{fs}{22}

\DeclareMathSymbol{\ExtofRedof}{\mathrel}{fs}{13}
\DeclareMathSymbol{\ExtofExpof}{\mathrel}{fs}{14}
\DeclareMathSymbol{\EmbinRedof}{\mathrel}{fs}{23}
\DeclareMathSymbol{\EmbinExpof}{\mathrel}{fs}{24}

\newcommand{\embin}{\prec}
\newcommand{\extof}{\succ}

\newcommand{\Nat}{\mathbb{N}}
\renewcommand{\set}[1]{\Set{#1}}
\newcommand{\power}{\mathop\mathscr{P}}
\newcommand{\num}[1]{{\underline{#1}}}
\newcommand{\ekv}{\leftrightarrow}

\newcommand{\rk}{\mathrm{rk}}
\newcommand{\bx}{\bar{x}}
\newcommand{\Ba}{\bar{a}}

\newcommand{\Prf}{\mathrm{Prf}}
\newcommand{\con}{\mathrm{Con}}

\newcommand{\codedset}{\mathop\mathrm{set}\nolimits}

\newcommand\satcon{\mathrm{SatCon}}

\newcommand{\leftfix}{}

 \def\godel#1{\,\godelalt{#1}\,}

 \def\godelalt#1{\setbox0=\hbox{$#1$}
 \dimen0=.3pt \dimen1=3.5pt 
 \dimen2=1.5pt \dimen3=1.5pt
 \dimen4=\dimen3 \advance\dimen4 by \ht0%
 \raise\dimen4\leftcorner %
 \dimen5=\dimen1 \advance\dimen5 by -\dimen2 %
 \kern-\dimen5 %
 #1 %
 \kern-\dimen5 %
 \raise\dimen4\rightcorner}
 
 \def\hlineseg{\vbox to \dimen0{\hrule width\dimen1 depth0pt height\dimen0}}
 
 \def\vlineseg{\hbox to \dimen0{\advance\dimen1 by -\dimen0 \vrule
     height\dimen0 depth\dimen1 width\dimen0}}
 
 \def\leftcorner{\hbox{\vlineseg \kern-\dimen0 \hlineseg}}
 
 \def\rightcorner{\hbox{\hlineseg \kern -\dimen0 \vlineseg}}

%% file: thesis.bbl
\begin{thebibliography}{KKK91}

\bibitem[AM74]{Apt.Marek:73}
Krzysztof~R. Apt and Wiktor Marek.
\newblock Second order arithmetic and related topics.
\newblock \emph{Ann. Math. Logic}, 6:177--229, 1973/74.

\bibitem[BS76]{Barwise.Schlipf:76}
Jon Barwise and John Schlipf.
\newblock An introduction to recursively saturated and resplendent models.
\newblock \emph{J. Symbolic Logic}, 41(2):531--536, 1976.

\bibitem[CF96]{Casanovas.Farre:96}
Enrique Casanovas and Rafel Farr{\'e}.
\newblock Omitting types in incomplete theories.
\newblock \emph{J. Symbolic Logic}, 61(1):236--245, 1996.

\bibitem[Eng02]{Engstrom:02}
Fredrik Engstr\"om.
\newblock \emph{Satisfaction classes in nonstandard models of arithmetic}.
\newblock Licentiate thesis, Chalmers University of Technology, 2002,
  math/0209408.

\bibitem[Fef61]{Feferman:60}
Solomon Feferman.
\newblock Arithmetization of metamathematics in a general setting.
\newblock \emph{Fund. Math.}, 49:35--92, 1960/1961.

\bibitem[Fri73]{Friedman:73}
Harvey Friedman.
\newblock Countable models of set theories.
\newblock In \emph{Cambridge Summer School in Mathematical Logic (Cambridge,
  1971)}, pages 539--573. Lecture Notes in Math., Vol. 337. Springer, Berlin,
  1973.

\bibitem[Hin78]{Hinman:78}
Peter~G. Hinman.
\newblock \emph{Recursion-theoretic hierarchies}.
\newblock Springer-Verlag, Berlin, 1978.
\newblock Perspectives in Mathematical Logic.

\bibitem[Jec03]{Jech:03}
Thomas Jech.
\newblock \emph{Set theory}.
\newblock Springer Monographs in Mathematics. Springer-Verlag, Berlin, 2003.
\newblock The third millennium edition, revised and expanded.

\bibitem[Kan96]{Kanovei:96}
Vladimir Kanovei.
\newblock On external {S}cott algebras in nonstandard models of {P}eano
  arithmetic.
\newblock \emph{J. Symbolic Logic}, 61(2):586--607, 1996.

\bibitem[Kay91]{Kaye:91}
Richard Kaye.
\newblock A generalization of {S}pecker's theorem on typical ambiguity.
\newblock \emph{J. Symbolic Logic}, 56(2):458--466, 1991.

\bibitem[KKK91]{Kaye.Kossak.ea:91}
Richard Kaye, Roman Kossak, and Henryk Kotlarski.
\newblock Automorphisms of recursively saturated models of arithmetic.
\newblock \emph{Ann. Pure Appl. Logic}, 55(1):67--99, 1991.

\bibitem[Kun80]{Kunen:80}
Kenneth Kunen.
\newblock \emph{Set theory}, volume 102 of \emph{Studies in Logic and the
  Foundations of Mathematics}.
\newblock North-Holland Publishing Co., Amsterdam, 1980.
\newblock An introduction to independence proofs.

\bibitem[Mos61]{Mostowski:61}
Andrzej Mostowski.
\newblock Formal system of analysis based on an infinitistic rule of proof.
\newblock In \emph{Infinitistic Methods (Proc. Sympos. Foundations of Math.,
  Warsaw, 1959)}, pages 141--166. Pergamon, Oxford, 1961.

\bibitem[MS04]{Mummert.Simpson:04}
Carl Mummert and Stephen~G. Simpson.
\newblock An incompleteness theorem for {$\beta\sb n$}-models.
\newblock \emph{J. Symbolic Logic}, 69(2):612--616, 2004.

\bibitem[Rog87]{Rogers:87}
Hartley Rogers, Jr.
\newblock \emph{Theory of recursive functions and effective computability}.
\newblock MIT Press, Cambridge, MA, second edition, 1987.

\bibitem[Sco62]{Scott:62}
Dana Scott.
\newblock Algebras of sets binumerable in complete extensions of arithmetic.
\newblock In \emph{Proc. Sympos. Pure Math., Vol. V}, pages 117--121. American
  Mathematical Society, Providence, R.I., 1962.

\bibitem[Sco65]{Scott:65}
Dana Scott.
\newblock Logic with denumerably long formulas and finite strings of
  quantifiers.
\newblock In \emph{Theory of Models (Proc. 1963 Internat. Sympos. Berkeley)},
  pages 329--341. North-Holland, Amsterdam, 1965.

\bibitem[Sim99]{Simpson:99}
Stephen~G. Simpson.
\newblock \emph{Subsystems of second order arithmetic}.
\newblock Springer-Verlag, Berlin, 1999.

\bibitem[Smo81]{Smorynski:81*4}
Craig Smory{\'n}ski.
\newblock Recursively saturated nonstandard models of arithmetic.
\newblock \emph{J. Symbolic Logic}, 46(2):259--286, 1981.

\bibitem[Smo82]{Smorynski:82}
Craig Smory{\'n}ski.
\newblock Addendum: ``{R}ecursively saturated nonstandard models of
  arithmetic''.
\newblock \emph{J. Symbolic Logic}, 47(3):493--494, 1982.

\bibitem[Wil75]{Wilmers:75}
George Wilmers.
\newblock \emph{Some problems in set theory: non-standard models and their
  applications to model theory}.
\newblock PhD thesis, University of Oxford, 1975.

\end{thebibliography}
